\documentclass[11pt]{article}

\usepackage[T1]{fontenc}
\usepackage[utf8]{inputenc}
\usepackage{lmodern}
\usepackage{microtype}

\usepackage{amsmath,amssymb,amsfonts,amsthm}
\usepackage{mathtools}

\usepackage[a4paper,margin=3cm]{geometry}
\usepackage{enumitem}
\usepackage{xcolor}
\usepackage[
    colorlinks=true,
    linkcolor=blue,
    citecolor=blue,
    urlcolor=blue
]{hyperref}

\usepackage[numbers,sort&compress]{natbib}

\usepackage{comment}

\theoremstyle{plain}
\newtheorem{theorem}{Theorem}[section]
\newtheorem{lemma}[theorem]{Lemma}
\newtheorem{proposition}[theorem]{Proposition}
\newtheorem{corollary}[theorem]{Corollary}

\theoremstyle{definition}
\newtheorem{definition}[theorem]{Definition}
\newtheorem{assumption}[theorem]{Assumption}
\newtheorem{example}[theorem]{Example}

\theoremstyle{remark}

\newcommand{\R}{\mathbb{R}}
\newcommand{\N}{\mathbb{N}}
\newcommand{\M}{\mathcal{M}}
\newcommand{\Scal}{\mathcal{S}}
\newcommand{\Leb}{\mathcal{L}}

\DeclareMathOperator{\supp}{supp}
\DeclareMathOperator{\sgn}{sgn}

\numberwithin{equation}{section}

\begin{document}

\title{\bfseries Existence and uniqueness of entropy solutions for nonautonomous pressureless gas dynamics}

\author{ Richard De~la~cruz\thanks{School of Mathematics and Statistics, Universidad Pedagógica y Tecnológica de Colombia, Tunja 150003, Colombia.
E-mail: \texttt{richard.delacruz@uptc.edu.co}} \and Rakib Mondal\thanks{Department of Mathematics, Indian Institute of Technology Bombay, Powai, Mumbai, Maharashtra 400076, India.
E-mail: \texttt{rkbmon90@gmail.com}} \and Wladimir Neves\thanks{Instituto de Matem\'atica, Universidade Federal do Rio de Janeiro, Cidade Universit\'aria, Rio de Janeiro 21945-970, Brazil.
E-mail: \texttt{wladimir@im.ufrj.br}} }

\date{}

\maketitle

\begin{abstract}
We study the Cauchy problem for one-dimensional pressureless gas dynamics with a continuous positive transport coefficient $\alpha(t)$ and a continuous nonnegative linear damping coefficient $\sigma(t)$, for every nonnegative locally finite Radon measure $\rho_0$ as initial density and every $u_0\in L^\infty(d\rho_0)$ as initial velocity.
A central contribution of this paper is a new and direct proof of uniqueness under the original assumptions of Huang and Wang.
It is carried out on the generalized potentials, and the energy condition enters only through the vanishing of the velocity dispersion at the initial time. The proof applies on every finite time strip and requires neither finite mass nor finite moments.
We then show that time-dependent transport and damping act on entropy solutions only through an effective time, so that the nonautonomous system is the homogeneous one observed on a time interval $[0,A_\infty)$. 
Together with uniqueness on strips, this yields a global entropy solution, unique in density and momentum at every time, for all such coefficients.
We also show that the energy condition cannot be dropped, and that its defect equals the mass of an atom times the variance of the velocities into which it splits.
As consequences, we obtain an explicit representation of the solution by a generalized potential, the Riemann solutions with delta shocks, and the vanishing-damping limit.\\
\end{abstract}

\bigskip


\medskip


\section{Introduction}
\label{intro}

Hyperbolic systems of conservation and balance laws provide the mathematical language for a wide range of phenomena in gas dynamics, traffic flow, sedimentation, chromatography, shallow water theory, and continuum mechanics \cite{Dafermos2016,Bressan2000,LeVeque2002,Smoller1994}.
Among them, the pressureless gas dynamics system occupies a distinctive place. It is weakly, but not strictly, hyperbolic. Compressive configurations carrying positive mass can produce point concentrations in finite time, even for smooth initial profiles. In one space dimension, the resulting Dirac masses, measures of the form $c\,\delta_{x_0}$ with $c>0$, move along curves, the so-called {\em delta-shocks}, whereas expansive data produce vacuum.
In either case the solution leaves the class of bounded or $BV$ entropy solutions of the theory of strictly hyperbolic systems and must be sought among Radon measures. Standard tools of the theory, such as $BV$ estimates and Kruzhkov entropy inequalities, do not apply in this setting, and both the notion of solution and the admissibility conditions must be formulated for measures.
For the pressureless system, all the terms of the equations, including the momentum flux $\rho u^2$, are measures, because the velocity is bounded almost everywhere with respect to the density. Delta-shocks are therefore weak solutions in the sense of measures. This distinguishes them from the singular shocks of Keyfitz and Kranzer \cite{KeyfitzKranzer1995}, which are defined as limits of approximate solutions and do not satisfy the equations in the sense of distributions \cite{Keyfitz2011}, and from the $\delta'$-shocks of Panov and Shelkovich \cite{PanovShelkovich2006}, which involve derivatives of Dirac measures.

The subject of this paper is the Cauchy problem for the nonautonomous pressureless system
\begin{equation}\label{NA}
\begin{cases}
\partial_t\rho+\alpha(t)\,\partial_x(\rho u)=0,\\[1mm]
\partial_t(\rho u)+\alpha(t)\,\partial_x(\rho u^2)=-\sigma(t)\rho u,\\[1mm]
(\rho,u)|_{t=0}=(\rho_0,u_0),
\end{cases}
\end{equation}
for $x\in\R$ and $t\ge0$, where $\rho=\rho(x,t)\ge0$ is the density, $u=u(x,t)$ is the velocity, and $\alpha,\sigma\in C([0,\infty))$, with $\alpha>0$ and $\sigma\ge0$, are chosen independently of each other. The initial density $\rho_0$ is a nonnegative measure in $\M_{\mathrm{loc}}(\R)$, the space of locally finite signed Radon measures on $\R$, and the initial velocity belongs to $L^\infty(\R,d\rho_0)$, the space of bounded functions that are measurable with respect to $\rho_0$, identified when they coincide $\rho_0$-almost everywhere. The coefficient $\alpha$ modulates the transport, and $\sigma$ is a time-dependent linear damping. For $\alpha\equiv1$ and $\sigma\equiv0$, \eqref{NA} is the homogeneous pressureless system. The present work continues De~la~cruz and Neves \cite{DelaCruzNeves2025}, where the Riemann problem was solved for a class of such systems with time-gradually-degenerate damping, and addresses the Cauchy problem with general Radon measure data.

The homogeneous system has been studied intensively over the last three decades, largely because of its role as a model of sticky-particle dynamics and of the formation of large-scale structures, such as galaxies, filaments, and voids, from pressureless matter in an expanding universe \cite{BrenierGrenier,ERykovSinai}.
Approaches based on measure-valued transport with discontinuous or one-sided Lipschitz coefficients were developed by Bouchut \cite{Bouchut1994}, Bouchut and James \cite{BouchutJames1998}, and Poupaud and Rascle \cite{PoupaudRascle1997}. Explicit measure solutions of the Cauchy problem for the associated transport equations were constructed by Wang, Huang, and Ding \cite{WangHuangDing1997}, Wang and Ding \cite{WangDing1997}, and Cheng, Li, and Zhang \cite{ChengLiZhang1997}, and the Riemann problem was solved by Sheng and Zhang \cite{ShengZhang1999}.
The decisive step for the Cauchy problem with general data was taken by Huang and Wang \cite{HuangWang2001}, who constructed entropy solutions for arbitrary nonnegative Radon measure initial density and bounded measurable initial velocity. They showed that uniqueness in this class requires, besides a one-sided estimate of Oleinik type on the velocity, familiar from scalar conservation laws, an energy condition, namely the weak continuity of the kinetic energy $u^2\rho$ at the initial time, and they proved by an explicit counterexample that this condition cannot be dropped \cite[see also][Example~2.1]{BouchutJames1999}.
Nonhomogeneous variants of the autonomous system were treated by Ding and Huang \cite{DingHuang2004} and Huang \cite{Huang2005}.
The Cucker--Smale model describes particles whose velocities tend to align through interactions. Motivated by its hydrodynamic description \cite{CuckerSmale2007,HaTadmor2008,HaLiu2009}, Jin \cite{Jin2016} studied a linear relaxation toward zero velocity, a simplification of alignment dissipation, and extended this theory to the system with the linear damping $-\rho u$ by reconstructing the generalized potential in the damped setting, after Ha, Huang, and Wang \cite{HaHuangWang2014} had treated the same system for locally integrable initial densities that are positive almost everywhere.
More recently, Tang \cite{Tang2025} studied the zero-relaxation and vanishing-damping limits of the system with constant linear damping within the same framework.

Uniqueness is the central question in this theory, because once mass concentrates the Oleinik condition alone no longer selects a solution. An essential contribution of this paper is a new proof of uniqueness of entropy solutions of the homogeneous system, under the same assumptions on the data as in \cite{HuangWang2001}, namely a nonnegative locally finite Radon measure as initial density and a bounded initial velocity. The proof applies on every time strip $[0,T)$ with $T\in(0,\infty]$, for velocities that are bounded only on each finite strip, and it requires neither a finite total mass nor finite moments of the density. The new step is the passage to the initial time, which we carry out on the generalized potentials themselves and which does not use the size of any exceptional set of initial points. This step needs a new argument. In \cite{HuangWang2001} it is based on the Lebesgue nullity of the set of initial points whose particles belong to a concentration at every positive time, and we give bounded initial data for which this set has full measure in an interval. This concerns only an auxiliary set of the proof, and the uniqueness theorem itself is contained in our result. Uniqueness on finite strips is, moreover, the form of uniqueness that the nonautonomous problem requires, because after the change of time used below the evolution on $[0,\infty)$ may correspond to a homogeneous evolution on a bounded time interval.

Compared with the homogeneous case, the nonautonomous problem has remained much less developed.
Two structural properties of the homogeneous system are used throughout the autonomous theory, namely that free particles move along straight lines and that the system is invariant under the scaling $(x,t)\mapsto(\lambda x,\lambda t)$, $\lambda>0$.
The generalized potential of the homogeneous theory is built from the straight-line characteristics, and the self-similar structure of the Riemann solution comes from the scaling invariance.
For general coefficients, neither property needs to hold in physical time.
To the best of our knowledge, previous work on the nonautonomous system has treated Riemann data.
De~la~cruz \cite{Delacruz2020} and De~la~cruz and Neves \cite{DelaCruzNeves2025} constructed Riemann solutions, with and without delta shocks, for damped pressureless systems, the latter for a class of time-dependent coefficients that contains the one considered here, by means of a generalized similarity variable and a vanishing-viscosity argument adapted to piecewise-constant data.
That work already uses the velocity rescaling and the effective time employed below.

Systems of the form \eqref{NA} arise naturally in cosmology, where pressureless matter moves in an expanding background.
To see how the coefficients arise, consider a one-dimensional, homogeneous, and expanding background, described by a function $a\in C^2([0,\infty))$, called the scale factor, such that two points that move only with the expansion of the background, at distance $d$ at time $0$, are at distance $a(t)d$ at time $t$. We assume that $a(t)>0$ and $a'(t)\ge0$ for every $t\ge0$ and that $a(0)=1$. The second derivative of $a$ is used only for the particle acceleration below, and for the coefficients themselves $a\in C^1([0,\infty))$ suffices.
Throughout the paper, a prime denotes differentiation with respect to the only variable of a function of one variable, here the time $t$.
If $R(t)$ is the physical position of a particle at time $t$, it is natural to describe it by its comoving position $X(t):=R(t)/a(t)$, and to measure its motion relative to the background by its peculiar velocity
$$
U(t):=R'(t)-\frac{a'(t)}{a(t)}\,R(t)=a(t)\,X'(t),
$$
where the ratio $a'/a$ is the relative expansion rate of the background, known in cosmology as the Hubble parameter \cite[equation~(22)]{MartelShapiro1998}.
A particle that moves only with the background has constant comoving position, hence $R(t)=a(t)X(0)$ and $R''=(a''/a)R$.
Suppose now that the gravitational field generated by the departures from homogeneity is neglected, so that every particle undergoes only this background acceleration, $R''=(a''/a)R$. Differentiating $U=R'-(a'/a)R$ and using $R'=a'X+aX'$, we obtain
$$
X'(t)=\frac{U(t)}{a(t)} \qquad \mbox{ and } \qquad U'(t)=-\frac{a'(t)}{a(t)}\,U(t).
$$
Thus, in comoving variables, a free particle is transported with speed $U/a$, and its peculiar velocity is damped at the rate $a'/a$.
For a continuum of such particles that stick together upon collision \cite{BrenierGrenier,ERykovSinai}, described by the comoving density $\rho(x,t)$, that is, mass per unit comoving length, and by the peculiar velocity field $u(x,t)$, so that $U(t)=u(X(t),t)$ along each particle path, the conservation of mass and momentum then gives the first two equations of \eqref{NA} with $\alpha(t)=1/a(t)$ and $\sigma(t)=a'(t)/a(t)$.
The passage to comoving coordinates produces the time-dependent coefficient in front of the fluxes, and the expansion acts on the momentum as a time-dependent linear damping.
In this cosmological model the two coefficients are not independent, since both are determined by the single function $a$ and $\sigma=a'/a=-\alpha'/\alpha$.
The cosmological coefficients are a particular case of those considered here, because the conditions $a>0$ and $a'\ge0$ give $\alpha>0$ and $\sigma\ge0$.

The starting point of our analysis is the transformation
$$
\tau=A(t) \qquad \mbox{ and } \qquad
w=e^{S(t)}u,
$$
where $S(t)=\int_0^t\sigma(s)\,ds$ and $A(t)=\int_0^t\alpha(s)e^{-S(s)}\,ds$.
The effective time $A(t)$ is the displacement at time $t$ of a free particle of unit initial velocity governed by \eqref{NA}, and $e^{-S(t)}$ is the factor by which its velocity has decayed.
For classical solutions, a direct computation shows that this change of variables turns \eqref{NA} into the homogeneous pressureless system, observed through a clock that runs at the rate $\alpha(t)e^{-S(t)}$.
We do not claim that this rescaling is new. For the cosmological coefficients it reduces, up to constant reference scales, to the supercomoving variables of Martel and Shapiro \cite{MartelShapiro1998}.
The content of the paper is to make the reduction rigorous for measure-valued solutions with general data, where it is not a formal consequence of the classical computation.

For hyperbolic systems, an entropy solution is usually a weak solution that satisfies $\partial_t\eta+\partial_xq\le0$ in the sense of distributions for every convex entropy $\eta$ with entropy flux $q$. For the homogeneous pressureless system, the pairs $\eta=\rho g(u)$ and $q=\rho u\,g(u)$, with $g$ convex, lead to the entropic solutions of Bouchut \cite{Bouchut1994}. The admissibility used here treats the velocity and the density separately. Where $\rho>0$, classical solutions of \eqref{NA} satisfy $\partial_tu+\alpha(t)u\,\partial_xu=-\sigma(t)u$, an equation of Burgers type for the velocity alone, and the first condition is the one-sided bound of Oleinik type
$$
\frac{u(x_2,t)-u(x_1,t)}{x_2-x_1}\le\frac{e^{-S(t)}}{A(t)}\qquad(x_1<x_2)
$$
for almost every $t>0$. The second condition is the energy condition, the weak continuity of $u^2\rho$ at $t=0$, where $u^2\rho$ is the entropy $\rho g(u)$ with $g(u)=u^2$. We call these two requirements the \emph{Oleinik--energy entropy condition}, and a weak solution that satisfies them an \emph{entropy solution}. No entropy inequality is imposed at positive times. For $\alpha\equiv1$ and $\sigma\equiv0$, the bound becomes $1/t$, and the notion reduces to the entropy solutions of \cite[Definition~2]{HuangWang2001}.

In this paper we prove that, for arbitrary continuous coefficients $\alpha>0$ and $\sigma\ge0$, arbitrary nonnegative Radon measure initial density, and bounded initial velocity, the Cauchy problem \eqref{NA} has a global entropy solution. It is unique in the sense that its density and its momentum are determined by the data at every time, so that its velocity is unique $\rho_t$-almost everywhere. Our main result reads as follows.

\begin{theorem}[Main theorem]\label{main-introduction}
Let $\alpha,\sigma\in C([0,\infty))$ with $\alpha>0$ and $\sigma\ge0$, let $S(t)=\int_0^t\sigma(s)\,ds$ and $A(t)=\int_0^t\alpha(s)e^{-S(s)}\,ds$, and put $A_\infty:=\lim\limits_{t\to\infty}A(t)\in(0,\infty]$.
Let $\rho_0\in\M_{\mathrm{loc}}(\R)$ be nonnegative and let $u_0\in L^\infty(\R,d\rho_0)$.
Then the following hold.
\begin{enumerate}[label=(\roman*)]
\item The Cauchy problem \eqref{NA} has a global entropy solution (in sense of Definition~\ref{entropy-solution}). If $(\rho^1,u^1)$ and $(\rho^2,u^2)$ are entropy solutions with the same initial data, then $\rho^1(\cdot,t)=\rho^2(\cdot,t)$ and $u^1(\cdot,t)\rho^1(\cdot,t)=u^2(\cdot,t)\rho^2(\cdot,t)$ as measures for every $t\ge0$. In particular, the velocity is unique $\rho(\cdot,t)$-almost everywhere.
\item For a pair $(\rho,u)$ on $\R\times[0,\infty)$, let $\mathcal T(\rho,u):=(\tilde\rho,w)$ be the pair on $\R\times[0,A_\infty)$ given by
$$
\tilde\rho(\cdot,\tau):=\rho\big(\cdot,A^{-1}(\tau)\big),
\qquad
w(\cdot,\tau):=e^{S(A^{-1}(\tau))}\,u\big(\cdot,A^{-1}(\tau)\big).
$$
Then $\mathcal T$ is a bijection from the set of entropy solutions of \eqref{NA} on $[0,\infty)$ onto the set of entropy solutions with the same initial data of the homogeneous pressureless system on $[0,A_\infty)$, which is \eqref{NA} with $\alpha\equiv1$, $\sigma\equiv0$, and $\tau$ in place of $t$. Under $\mathcal T$, the Oleinik--energy entropy condition of \eqref{NA} corresponds to that of the homogeneous system.
\item Let $\bar a\in\R$ with $\rho_0(\{\bar a\})=0$, and let
$$
F(y;x,t):=\int_{\bar a}^y\big(A(t)u_0(\eta)+\eta-x\big)\,d\rho_0(\eta),
\qquad
\nu(x,t):=\inf_{y\in\R}F(y;x,t),
$$
where, for $a,y\in\R$, a measure $\mu$, and a locally $\mu$-integrable function $f$, the oriented integral is $\int_a^yf\,d\mu:=\int_{[a,y)}f\,d\mu$ if $y\ge a$ and $\int_a^yf\,d\mu:=-\int_{[y,a)}f\,d\mu$ if $y<a$. Then the entropy solution of \eqref{NA} satisfies, for every $t>0$,
$$
\rho(\cdot,t)=-\partial_x^2\nu(\cdot,t)\qquad\text{in }\mathcal D'(\R),
$$
and $\alpha(t)\,u\rho=\partial_x\partial_t\nu$ in $\mathcal D'(\R\times(0,\infty))$.
\end{enumerate}
\end{theorem}

Part (i) is Theorem~\ref{existence-uniqueness}, part (ii) is Theorem~\ref{entropy-correspondence} with $T=\infty$, and part (iii) is Corollary~\ref{potential-solution}. The proof of the theorem requires three results, which we regard as the main contributions of the paper.

The first is the uniqueness of entropy solutions of the homogeneous system on strips, Theorem~\ref{strip-uniqueness}.
The effective time may remain bounded, and then the whole evolution of \eqref{NA} on $[0,\infty)$ corresponds to the homogeneous evolution on the bounded interval $[0,A_\infty)$. Every concentration that the homogeneous dynamics would produce only after the time $A_\infty$ never forms.
For $\alpha\equiv1$, this happens for every constant damping and, for $\sigma(t)=\sigma_0/(1+t)$, exactly when $\sigma_0>1$.
It also happens for the cosmological coefficients of an Einstein--de Sitter background \cite{EinsteinDeSitter1932}, in which the scale factor grows like the power $2/3$ of the time elapsed since the origin of the expansion. If $t_0>0$ is the time elapsed at the initial instant, then $a(t)=(1+t/t_0)^{2/3}$ and $A_\infty=3t_0$.
A bounded effective time already occurs for constant damping, where $A(t)=1-e^{-t}$ enters the generalized potential of Jin \cite{Jin2016}.
In this situation uniqueness is needed on a strip $[0,T)$ with $T$ finite, for velocities bounded only on each finite strip, and it does not follow from uniqueness for all positive times, because an entropy solution on a strip need not be the restriction of an entropy solution on the half-line. In the proof of Theorem~\ref{main-introduction}, the results of \cite{HuangWang2001} are used only for existence.
Our proof follows the three steps of the homogeneous theory, namely a description of the solution at positive times, its representation through a generalized potential built from its state at a positive time $s$, and the passage $s\to0^+$. The first two steps are obtained from a Lagrangian description in mass coordinates, in which every concentration is an interval of mass coordinates and the Oleinik condition enters only through the forward uniqueness of characteristics.
The passage to the initial time is carried out directly on the generalized potentials. The initial traces of the mass, the momentum, and the kinetic energy imply that the joint distributions of position and velocity converge, as $s\to0^+$, to a measure concentrated on the graph of the initial velocity, Lemma~\ref{zero-variance}. This controls the contribution of the atoms of the initial density, so that the generalized potentials converge to the generalized potential of the data, Theorem~\ref{convergence}, and the density and the momentum are recovered by differentiating the limit in the sense of distributions.
This argument does not use the size of any exceptional set of initial points. In Section~\ref{immediate-concentration} we give bounded initial data for which the particles issued from every point of $(0,1)$ outside a countable set belong to a concentration at every positive time, both for the generalized characteristics and for the limits of the regularized characteristics, although at every positive time the solution consists of finitely many atoms.

The second is the transfer of the homogeneous theory to \eqref{NA} through an intrinsic solution class.
Our weak and entropy solutions are defined in terms of $(\rho,u)$ and the coefficients only, without a cumulative mass potential and without a velocity defined at every point.
We construct the cumulative potential of every weak solution from the weak form of the mass equation and the continuity in time of the density and the momentum, in Proposition~\ref{intrinsic-potential}, and in Lemma~\ref{joint-velocity} we construct a jointly Borel velocity that satisfies the Oleinik inequality at every pair of points, including points of vacuum, for almost every time, and defines the given momentum at every time.
For $\alpha\equiv1$ and $\sigma\equiv0$, the intrinsic class coincides, up to the choice of a representative of the velocity, with the class of \cite{HuangWang2001}, which is formulated for a potential $m$ with $\partial_xm=\rho$ and $\partial_\tau m+w\,\partial_xm=0$. Hence the homogeneous theory itself can be formulated without a potential.
In this class, the velocity rescaling and the effective time identify the entropy solutions of \eqref{NA} with those of the homogeneous system, so that existence and uniqueness transfer to independent time-dependent transport and damping coefficients.

The third is a quantitative description of the role of the energy condition.
A Dirac mass may split into particles whose velocities have the same mean, and every such splitting yields a weak solution satisfying the Oleinik condition.
We show in Proposition~\ref{nonuniqueness} that the defect in the initial trace of $u^2\rho$ equals the mass of the atom times the variance of the normalized velocity distribution, so that the energy condition selects exactly the solution without splitting.
This places the examples of Bouchut and James \cite{BouchutJames1999} and Huang and Wang \cite{HuangWang2001} in a single family. The same quantity governs the proof of uniqueness, where the energy condition forces the limit of the joint distributions of position and velocity at the initial time to carry no velocity dispersion.

A further technical point is that $A$ is in general only of class $C^1$, so the test functions of \eqref{NA} and of the homogeneous system correspond to each other only in the class $C^1$. We therefore extend the weak formulation to $C^1$ test functions, using a uniform local mass bound.

The remaining results are consequences of the correspondence.
The theory of Jin \cite{Jin2016} for the damping $-\rho u$ is recovered as the case $\alpha\equiv1$ and $\sigma\equiv1$, which explains the similarity between his generalized potential and that of the homogeneous theory.
The vanishing-damping limit for the density and the momentum, locally uniformly in time, follows from the continuity in time of the homogeneous solution alone, without the analysis of minimizers used by Tang \cite{Tang2025}.
The generalized potential transfers through the change of variables and determines the solution explicitly from the data. The generalized characteristics transfer as well, and their speed is $\alpha(t)$ times the velocity. The Riemann solutions of De~la~cruz and Neves \cite{DelaCruzNeves2025} are recovered and shown to be the unique entropy solutions in the class of measure-valued solutions.

Related equivalence questions are treated by Carrillo and Galtung \cite{CarrilloGaltung2026}, who prove the equivalence of entropy and Lagrangian solutions of pressureless Euler systems on the line for probability measures with finite second moment and velocities of finite kinetic energy, and who include a constant linear damping. Their velocities need not be bounded. The present framework requires bounded velocities, but it allows arbitrary locally finite density measures, without a global mass or moment assumption, together with time-dependent transport and damping coefficients.

The paper is organized as follows.
In Section~\ref{preliminaries} we fix the hypotheses on the coefficients and establish the properties of the primitive $S$ of $\sigma$ and of the effective time $A$ that are used throughout.
Section~\ref{classical-solutions} proves the transformation theorem for classical solutions, where the cancellation of the damping term and the reduction of $\alpha$ and $\sigma$ to the single factor $A'$ appear in their simplest form.
Section~\ref{weak-solutions} contains the analytical core of the paper. There we introduce the intrinsic notion of weak solution, show that it provides a velocity that is jointly Borel in space and time, construct the cumulative mass potential from every weak solution, extend the weak formulation to $C^1$ test functions, and prove the equivalence of weak solutions.
In Section~\ref{entropy-cauchy-section} we formulate the intrinsic Oleinik and energy conditions, show by an explicit family of solutions that the energy condition cannot be dropped, and prove the equivalence of the two Cauchy problems.
Section~\ref{global-existence-uniqueness} proves global existence and uniqueness. There we recall the formulation of the homogeneous theory and obtain existence from its existence theorem, construct a velocity that satisfies the Oleinik inequality at every pair of points, give a direct proof of uniqueness on finite time strips through the convergence of generalized potentials, and show that the intrinsic solution class coincides with the homogeneous one.
Section~\ref{lagrangian} describes entropy solutions in mass coordinates. There every concentration is an interval of mass coordinates, the velocity of each portion of mass is identified, and the representation of the solution at positive times used in Section~\ref{global-existence-uniqueness} is proved.
In Section~\ref{characteristics} the generalized potential, the representation formulas for the solution, and the generalized characteristics are transferred to \eqref{NA}, together with the characteristic differential inclusion and the speed of the forward characteristics.
Section~\ref{examples} computes the effective time for constant, time-gradually-degenerate, and cosmological coefficients, proves the vanishing-damping limit, compares the transformed Riemann solutions with those of De~la~cruz and Neves \cite{DelaCruzNeves2025}, and gives bounded data for which the particles issued from almost every initial point belong to a concentration at every positive time.

\section{Preliminaries}
\label{preliminaries}

The coefficients $\alpha$ and $\sigma$ enter the transformation only through two auxiliary functions, the primitive $S$ of $\sigma$ and the effective time $A$, both defined below.
In this section we recall the problem, fix the standing hypotheses on the coefficients, introduce these two functions, and record their basic properties.
We consider the Cauchy problem \eqref{NA} with data $\rho_0$ and $u_0$ as in the introduction. Measures in $\M_{\mathrm{loc}}(\R)$ act on $C_c(\R)$ by integration. Throughout, weak-$*$ convergence means convergence of these pairings against every function in $C_c(\R)$. No finiteness of the total mass on $\R$ is assumed.
Since $u_0$ is determined only $\rho_0$-almost everywhere, its values on regions without mass play no role.
For the moment \eqref{NA} is regarded formally, and the notions of classical, weak, and entropy solution are introduced once the transformation is available.
System~\eqref{NA} makes sense formally for essentially arbitrary coefficients, but the analysis below requires the following standing assumptions.

\begin{assumption}\label{coefficients}
The coefficients $\alpha,\sigma:[0,\infty)\to\R$ satisfy
\begin{enumerate}[label=(H\arabic*)]
\item\label{alpha} $\alpha\in C([0,\infty))$ and $\alpha(t)>0$ for every $t\ge0$.
\item\label{sigma} $\sigma\in C([0,\infty))$ and $\sigma(t)\ge0$ for every $t\ge0$.
\end{enumerate}
\end{assumption}

To identify the two auxiliary functions, consider a single free particle governed by \eqref{NA}. In the absence of interactions, its position $X(t)$ and its velocity $U(t)$ satisfy $X'(t)=\alpha(t)U(t)$ and $U'(t)=-\sigma(t)U(t)$, with $X(0)=x_0$ and $U(0)=v_0$.
Let $S(t):=\int_0^t\sigma(s)\,ds$.
The second equation gives $U(t)=v_0e^{-S(t)}$, so that $e^{-S(t)}$ is the factor by which the velocity of a free particle has decayed at time $t$.
The first equation then gives $X(t)=x_0+v_0A(t)$, where
\begin{equation}\label{effective-time}
A(t):=\int_0^t\alpha(s)e^{-S(s)}\,ds.
\end{equation}
We call $A$ the effective time.
Although $S$ and $A$ are obtained from the motion of a single free particle, the change of variables defined with them, which takes $A(t)$ as the new time variable and multiplies the velocity by $e^{S(t)}$, is applied below to arbitrary solutions of \eqref{NA}.

The hypotheses of Assumption~\ref{coefficients} enter as follows.
Strict positivity of $\alpha$ makes $A$ strictly increasing, hence invertible, which is what the substitution $\tau=A(t)$ requires.
Continuity of $\sigma$ makes $S$, and hence $e^{\pm S}$, continuously differentiable. This is needed to transport the test functions of the momentum equation, which carry the factor $e^{S(t)}$, within the class of $C^1$ functions. The sign condition $\sigma\ge0$ expresses the interpretation of $\sigma$ as a friction coefficient and yields $S\ge0$, so that $0<e^{-S}\le1$.
Assumption~\ref{coefficients} includes time-gradually-degenerate damping coefficients such as $\sigma(t)=\sigma_0(1+t)^{-\ell}$, with $\sigma_0>0$ and $\ell\ge0$, of the type considered by De~la~cruz and Neves \cite{DelaCruzNeves2025}, and the cosmological coefficients $\alpha=1/a$ and $\sigma=a'/a$ of the introduction whenever $a\in C^1([0,\infty))$ is positive and nondecreasing.

The name effective time reflects the following property.
A free particle of \eqref{NA} with initial position $x_0$ and initial velocity $v_0$ occupies at time $t$ the position $x_0+v_0A(t)$, which a free particle of the homogeneous system, that is, of \eqref{NA} with $\alpha\equiv1$ and $\sigma\equiv0$, reaches at time $\tau=A(t)$.
Thus $A(t)$ is the time that the homogeneous dynamics needs to produce the configuration that \eqref{NA} reaches at the physical time $t$. In particular, all events of the dynamics are synchronized by $A$. For instance, two particles with $x_1<x_2$ and $u_1>u_2$ collide when $A(t)=(x_2-x_1)/(u_1-u_2)$, which is their collision time for the homogeneous system.
Since $A$ is increasing, the limit $A_\infty:=\lim\limits_{t\to\infty}A(t)$ exists in $(0,\infty]$. If $A_\infty<\infty$, only the collisions whose homogeneous collision time is strictly smaller than $A_\infty$ take place at finite physical time. A homogeneous event at exactly $A_\infty$ is not reached at any finite physical time.

\begin{lemma}\label{properties-A}
Under Assumption~\ref{coefficients}, the following hold.
\begin{enumerate}[label=(\roman*)]
\item $S\in C^1([0,\infty))$, $S'=\sigma$, $S(0)=0$, and $0<e^{-S(t)}\le1$ for every $t\ge0$.
\item $A\in C^1([0,\infty))$, $A(0)=0$, and $A'(t)=\alpha(t)e^{-S(t)}>0$ for every $t\ge0$.
\item $A$ is strictly increasing and maps $[0,\infty)$ bijectively onto $[0,A_\infty)$. Moreover, $A^{-1}\in C^1([0,A_\infty))$, and both $A$ and $A^{-1}$ are Lipschitz on every compact subinterval of their domains.
\end{enumerate}
\end{lemma}

\begin{proof}
(i) Since $\sigma$ is continuous, the fundamental theorem of calculus gives $S\in C^1([0,\infty))$ with $S'=\sigma$ at every point, and $S(0)=0$. Since $\sigma\ge0$, we have $S\ge0$, hence $0<e^{-S}\le1$.

(ii) The integrand $\alpha e^{-S}$ in \eqref{effective-time} is continuous by \ref{alpha} and (i), so the fundamental theorem of calculus gives $A\in C^1([0,\infty))$ with $A'=\alpha e^{-S}$, and positivity follows from \ref{alpha}.

(iii) By (ii), $A$ is continuous and strictly increasing with $A(0)=0$, hence a bijection of $[0,\infty)$ onto $[0,A_\infty)$. Since $A'$ is continuous and positive, the inverse function theorem gives $A^{-1}\in C^1([0,A_\infty))$ with $(A^{-1})'=1/(A'\circ A^{-1})$. On a compact subinterval, $A'$ and $(A^{-1})'$ are continuous, hence bounded, which gives the Lipschitz property.
\end{proof}

Lemma~\ref{properties-A}(iii) implies, in particular, that $A$ and $A^{-1}$ map Lebesgue null sets onto Lebesgue null sets and Lebesgue measurable sets onto Lebesgue measurable sets, a fact used repeatedly to transfer statements valid for almost every time.

With $S$ and $A$ at hand, we can state the system to which \eqref{NA} will be reduced, namely the homogeneous pressureless gas dynamics system
\begin{equation}\label{H}
\begin{cases}
\partial_\tau\rho+\partial_x(\rho w)=0,\\[1mm]
\partial_\tau(\rho w)+\partial_x(\rho w^2)=0,
\end{cases}
\end{equation}
studied by Huang and Wang \cite{HuangWang2001} for nonnegative $\rho_0\in\M_{\mathrm{loc}}(\R)$ and $u_0\in L^\infty(\R,d\rho_0)$.

\subsection{Transformation of classical solutions}
\label{classical-solutions}

We begin with classical solutions, for which the effect of the transformation is a short computation. By a classical solution of \eqref{NA} on $[0,T)$ we mean a pair $(\rho,u)\in C^1(\R\times[0,T))$ that satisfies \eqref{NA} pointwise, and classical solutions of \eqref{H} are defined in the same way.
This shows why the transformation works before we address the additional care required for Radon measure data.
The cancellations that appear below, namely the disappearance of the damping term from the momentum equation and the collapse of $\alpha(t)$ and $\sigma(t)$ into the single factor $A'(t)$, reappear in the weak formulation, where only the functional setting changes.

The free-particle computation of Section~\ref{preliminaries} suggests the change of variables that replaces $t$ by the effective time $A(t)$ and removes from the velocity the factor $e^{-S(t)}$ by which the velocity of a free particle decays.
Applied to a general pair $(\rho,u)$, this substitution reads
$\tau=A(t)$,
and, for $(x,\tau)\in\R\times[0,A_\infty)$ and $t=A^{-1}(\tau)$,
\begin{equation}\label{w}
\tilde\rho(x,\tau)=\rho(x,t),
\qquad \mbox{ and } \qquad
w(x,\tau)=e^{S(t)}u(x,t).
\end{equation}
Since $A$ is a bijection of $[0,\infty)$ onto $[0,A_\infty)$, this substitution is invertible, with inverse given by
\begin{equation}\label{Tinverse}
\rho(x,t)=\tilde\rho(x,A(t))
\qquad \mbox{ and } \qquad
u(x,t)=e^{-S(t)}w(x,A(t)).
\end{equation}
The transformation leaves the density unchanged, apart from the relabeling of time, and modifies only the time variable and the velocity.
The following theorem shows that, for classical solutions, it turns \eqref{NA} into \eqref{H} exactly.

\begin{theorem}\label{transformation}
Let Assumption~\ref{coefficients} hold.
A pair $(\rho,u)\in C^1(\R\times[0,\infty))$ satisfies the first two equations of \eqref{NA} if and only if the pair $(\tilde\rho,w)$ defined by \eqref{w}, which belongs to $C^1(\R\times[0,A_\infty))$, satisfies \eqref{H}.
In particular, the substitution \eqref{w} is a bijection between classical solutions of \eqref{NA} on $\R\times[0,\infty)$ and classical solutions of \eqref{H} on $\R\times[0,A_\infty)$ with the same initial data.
\end{theorem}

\begin{proof}
By Lemma~\ref{properties-A}, $S$, $A$, and $A^{-1}$ are of class $C^1$, so $(\rho,u)$ is of class $C^1$ if and only if $(\tilde\rho,w)$ is. For a $C^1$ function $f(x,\tau)$, the function $F(x,t):=f(x,A(t))$ satisfies, by the chain rule,
\begin{equation}\label{rho-chain}
\partial_tF(x,t)=A'(t)\,\partial_\tau f(x,A(t)),\qquad \partial_xF(x,t)=\partial_xf(x,A(t)).
\end{equation}
By \eqref{Tinverse}, $\rho=\tilde\rho(\cdot,A(t))$, $\rho u=e^{-S(t)}(\tilde\rho w)(\cdot,A(t))$, and $\rho u^2=e^{-2S(t)}(\tilde\rho w^2)(\cdot,A(t))$.
For the continuity equation, \eqref{rho-chain} and the identity $\alpha(t)e^{-S(t)}=A'(t)$ give
$$
\partial_t\rho+\alpha(t)\,\partial_x(\rho u)
=A'(t)\,\partial_\tau\tilde\rho+\alpha(t)e^{-S(t)}\partial_x(\tilde\rho w)
=A'(t)\big[\partial_\tau\tilde\rho+\partial_x(\tilde\rho w)\big],
$$
where the bracket is evaluated at $(x,A(t))$.
For the momentum equation, $S'=\sigma$ and \eqref{rho-chain} give
$$
\partial_t(\rho u)
=-\sigma(t)e^{-S(t)}\tilde\rho w+e^{-S(t)}A'(t)\,\partial_\tau(\tilde\rho w)
=-\sigma(t)\rho u+e^{-S(t)}A'(t)\,\partial_\tau(\tilde\rho w),
$$
and $\alpha(t)\,\partial_x(\rho u^2)=\alpha(t)e^{-2S(t)}\partial_x(\tilde\rho w^2)=e^{-S(t)}A'(t)\,\partial_x(\tilde\rho w^2)$.
Hence
$$
\partial_t(\rho u)+\alpha(t)\,\partial_x(\rho u^2)+\sigma(t)\rho u
=e^{-S(t)}A'(t)\big[\partial_\tau(\tilde\rho w)+\partial_x(\tilde\rho w^2)\big],
$$
the damping terms having cancelled exactly.

Since $A'(t)>0$ and $e^{-S(t)}A'(t)=\alpha(t)e^{-2S(t)}>0$, the two equations of \eqref{NA} hold at $(x,t)$ if and only if the two equations of \eqref{H} hold at $(x,A(t))$.
As $A$ maps $[0,\infty)$ onto $[0,A_\infty)$, this proves the equivalence.
Finally, $S(0)=0$ and $A(0)=0$ give $(\tilde\rho,w)(\cdot,0)=(\rho,u)(\cdot,0)$, so the initial data are the same.
\end{proof}

Classical solutions, however, are the exception for the pressureless system, since smooth compressive data generically produce concentrations in finite time, and the relevant solution class is that of Radon measures.
The next section shows that the same computation survives in that class.

\section{Weak solutions}
\label{weak-solutions}

This section carries out, in the measure-valued setting, the computation of Section~\ref{classical-solutions}. The algebra is the same, but two points have no counterpart in the classical case.
First, the theory of Huang and Wang \cite{HuangWang2001} is formulated in terms of a cumulative mass potential, whereas our definition of weak solution refers only to $(\rho,u)$.
We show that every weak solution possesses such a potential, which will allow us to identify the two solution classes. Second, the test functions of the weak formulation must be transported through $\tau=A(t)$, and, since $A$ is in general only $C^1$, this transport does not preserve the class $C_c^\infty$ in which the weak formulation is stated.

Since the density may concentrate in finite time, we work in a measure-theoretic framework. Throughout, $\M_{\mathrm{loc}}(\R)$ is endowed with the weak-$*$ topology defined by its pairings with $C_c(\R)$, and a family $(\mu_t)$ is weak-$*$ continuous if $t\mapsto\int_\R\varphi\,d\mu_t$ is continuous for every $\varphi\in C_c(\R)$.

In this setting the substitution \eqref{w} acts on families of measures rather than on functions, and we now define it as a map between explicit sets.
For $T\in(0,\infty]$, let $\mathcal P_T$ be the set of pairs $(\rho,u)$ formed by a family $\rho=(\rho_t)_{t\in[0,T)}$ of nonnegative measures in $\M_{\mathrm{loc}}(\R)$ and a family $u=(u(\cdot,t))_{t\in[0,T)}$ of bounded functions such that $u(\cdot,t)$ is $\rho_t$-measurable for every $t\in[0,T)$.
With the convention $A(\infty):=A_\infty$, we define $\mathcal T_T:\mathcal P_T\longrightarrow\mathcal P_{A(T)}$ as $\mathcal T_T(\rho,u):=(\tilde\rho,w)$, where, for every $\tau\in[0,A(T))$,
\begin{equation*}
\tilde\rho_\tau:=\rho_{A^{-1}(\tau)}
\qquad \mbox{ and } \qquad
w(\cdot,\tau):=e^{S(A^{-1}(\tau))}\,u(\cdot,A^{-1}(\tau)).
\end{equation*}
By Lemma~\ref{properties-A}(iii), $A$ maps $[0,T)$ bijectively onto $[0,A(T))$, so $\mathcal T_T$ is a bijection, and its inverse $\mathcal T_T^{-1}:\mathcal P_{A(T)}\to\mathcal P_T$ is given by
$$
\rho_t=\tilde\rho_{A(t)}
\qquad \mbox{ and } \qquad
u(\cdot,t)=e^{-S(t)}\,w(\cdot,A(t)),
\qquad t\in[0,T).
$$
If $\rho_t=f(\cdot,t)\,dx$ for a function $f$, then $\tilde\rho_\tau=f(\cdot,A^{-1}(\tau))\,dx$, so on pairs given by functions $\mathcal T_T$ reduces to the substitution \eqref{w} of Section~\ref{classical-solutions}.
The results of this section and of the next one assert that $\mathcal T_T$ maps the weak solutions, and then the entropy solutions, of \eqref{NA} on $[0,T)$ onto those of \eqref{H} on $[0,A(T))$.
When the value of $T$ is clear from the context, we write $\mathcal T$ instead of $\mathcal T_T$.

\subsection{Weak formulation and cumulative potential}

We first fix the notion of solution. It is stated intrinsically, that is, in terms of $(\rho,u)$, $\alpha$, and $\sigma$ only. In particular, it does not involve the transformation $\mathcal T$, and, unlike the formulation of the homogeneous theory, which is stated for a pair formed by a cumulative mass potential and a velocity, it does not require the construction of a potential.
Whether a given pair is a weak solution can therefore be decided by checking conditions on the pair itself. Measurability with respect to a Radon measure means measurability for its completion. In measure formulas, $\rho_t$ denotes the measure at time $t$, whereas $\partial_t\rho$ denotes a distributional time derivative.

\begin{definition}[Weak solution]\label{weak-solution}
Let $T\in(0,\infty]$. A pair $(\rho,u)$ is called a weak solution of \eqref{NA} on $[0,T)$ if the following conditions hold:
\begin{enumerate}[label=(\alph*)]
\item $t\mapsto\rho_t$ and $t\mapsto u(\cdot,t)\rho_t$ are weak-$*$ continuous maps from $[0,T)$ into $\M_{\mathrm{loc}}(\R)$, and $\rho_t\ge0$ for every $t\in[0,T)$.
\item for every $t\in[0,T)$, $u(\cdot,t)$ is bounded and measurable with respect to $\rho_t$,
\begin{equation}\label{local-uniform-u}
\sup_{t\in I}\|u(\cdot,t)\|_{L^\infty(d\rho_t)}<\infty
\qquad\text{for every compact }I\subset[0,T).
\end{equation}
\item for all $\Xi,\Psi\in C_c^\infty(\R\times(0,T))$,
\begin{align}
\int_0^T\!\!\int_{\R}
\partial_t\Xi\,d\rho_t\,dt
+\int_0^T\!\!\int_{\R}
\alpha(t)\,\partial_x\Xi\,u\,d\rho_t\,dt
&=0,\label{weak1-rho}\\
\int_0^T\!\!\int_{\R}
\left(
\partial_t\Psi\,u
+\alpha(t)\partial_x\Psi\,u^2
-\sigma(t)\Psi\,u
\right)\,d\rho_t\,dt
&=0.
\label{weak2}
\end{align}
\item the initial data are attained in the sense that, as $t\to0^+$, $\rho_t\stackrel{*}{\rightharpoonup}\rho_0$ and $u(\cdot,t)\,\rho_t\stackrel{*}{\rightharpoonup} u_0\,\rho_0$ in $\M_{\mathrm{loc}}(\R)$.
\end{enumerate}
When $T=\infty$, we simply speak of a weak solution of \eqref{NA}.
Weak solutions $(\rho,w)$ of \eqref{H} on $[0,T)$ are defined by conditions (a), (b), and (d), with $(\rho,w)$ in place of $(\rho,u)$, together with the identities
\begin{align}
\int_0^T\!\!\int_{\R}\big(\partial_\tau\Phi+w\,\partial_x\Phi\big)\,d\rho_\tau\,d\tau&=0,
\label{weak1-H}\\
\int_0^T\!\!\int_{\R}\big(\partial_\tau Z\,w+\partial_xZ\,w^2\big)\,d\rho_\tau\,d\tau&=0,
\label{weak2-H}
\end{align}
for all $\Phi,Z\in C_c^\infty(\R\times(0,T))$, which are \eqref{weak1-rho} and \eqref{weak2} with $\alpha\equiv1$ and $\sigma\equiv0$.
\end{definition}

By (a), condition (d) states that $\rho_t$ and $u(\cdot,t)\rho_t$ take the values $\rho_0$ and $u_0\rho_0$ at $t=0$.

The weak formulation involves integrals of continuous functions against the measures $\rho_t$, and also the integrals of $u(\cdot,t)^2$, which are not pairings of a weak-$*$ continuous family with a fixed continuous function.
The construction of the cumulative potential below requires, moreover, integrals of functions that are only Borel. The following lemma shows that such integrals depend measurably on the parameters.

\begin{lemma}\label{kernel}
Let $T\in(0,\infty]$, let $(\mu_t)_{t\in[0,T)}$ be a weak-$*$ continuous family in $\M_{\mathrm{loc}}(\R)$, and let $Z$ be a Borel subset of a Euclidean space.
Then the following hold.
\begin{enumerate}[label=(\roman*)]
\item For every bounded Borel set $E\subset\R$, the function $t\mapsto\mu_t(E)$ is Borel on $[0,T)$.
\item Let $f:Z\times[0,T)\times\R\to\R$ be a bounded Borel function, and assume that there is $n\in\N$ such that $f(z,t,y)=0$ whenever $|y|\ge n$.
Then the function $(z,t)\mapsto\int_\R f(z,t,y)\,d\mu_t(y)$ is Borel on $Z\times[0,T)$.
\item Assume that $\mu_t\ge0$ for every $t\in[0,T)$.
For a Borel set $Q\subset\R\times[0,T)$ and $t\in[0,T)$, put $Q_t:=\{y\in\R:(y,t)\in Q\}$.
Then $t\mapsto\mu_t(Q_t)$ is a Borel function with values in $[0,\infty]$, and
$$
\mu(Q):=\int_0^T\mu_t(Q_t)\,dt
$$
defines a measure $\mu$ on the Borel subsets of $\R\times[0,T)$.
Moreover, for every nonnegative Borel function $f$ on $\R\times[0,T)$, the function $t\mapsto\int_\R f(y,t)\,d\mu_t(y)$ is Borel and
$$
\int_{\R\times[0,T)}f\,d\mu=\int_0^T\!\!\int_\R f(y,t)\,d\mu_t(y)\,dt.
$$
In particular, if $N\subset\R\times[0,T)$ is a Borel set with $\mu(N)=0$, then $\mu_t(N_t)=0$ for almost every $t\in[0,T)$.
\end{enumerate}
\end{lemma}

\begin{proof}
For a set $E$, we denote by $\mathbf 1_E$ its characteristic function, and, for $n\in\N$, we put $I_n:=(-n,n)$.
Since each $\mu_t$ is a Radon measure, its total variation $|\mu_t|$ is finite on bounded sets.

(i) Fix $n\in\N$ and let $U\subset I_n$ be open.
There is a sequence $\{\chi_j\}_{j\ge1}\subset C_c(\R)$ with $0\le\chi_j\le\mathbf 1_U$ that increases pointwise to $\mathbf 1_U$ as $j\to\infty$.
Since $|\mu_t|(U)<\infty$, the dominated convergence theorem gives
$$
\mu_t(U)=\lim_{j\to\infty}\int_\R\chi_j\,d\mu_t,
$$
and each map $t\mapsto\int_\R\chi_j\,d\mu_t$ is continuous by the weak-$*$ continuity of the family.
Hence $t\mapsto\mu_t(U)$ is Borel.
Let $\mathcal G_n$ be the family of Borel sets $E\subset I_n$ such that $t\mapsto\mu_t(E)$ is Borel.
It contains $I_n$, it is closed under nested differences because $\mu_t(F\setminus E)=\mu_t(F)-\mu_t(E)$ with finite terms when $E\subset F\subset I_n$, and it is closed under increasing unions by the countable additivity of $\mu_t$ on $I_n$.
Thus $\mathcal G_n$ is a $\lambda$-system that contains the $\pi$-system of the open subsets of $I_n$, which generates the Borel $\sigma$-algebra of $I_n$.
By Dynkin's $\pi$--$\lambda$ theorem, $\mathcal G_n$ contains every Borel subset of $I_n$.
Since every bounded Borel set is contained in some $I_n$, this proves (i).

(ii) Fix $n$ as in the statement, and let $\mathcal H_n$ be the set of bounded Borel functions $F:Z\times[0,T)\times\R\to\R$ such that the function $(z,t)\mapsto\int_\R F(z,t,y)\mathbf 1_{I_n}(y)\,d\mu_t(y)$, which is defined because each section $F(z,t,\cdot)$ is Borel and $|\mu_t|(I_n)<\infty$, is Borel on $Z\times[0,T)$.
Then $\mathcal H_n$ is a vector space and, by the dominated convergence theorem, it is closed under pointwise limits of sequences bounded by a common constant.
If $B\subset Z$, $J\subset[0,T)$, and $E\subset\R$ are Borel sets, then
$$
\int_\R\mathbf 1_B(z)\mathbf 1_J(t)\mathbf 1_E(y)\mathbf 1_{I_n}(y)\,d\mu_t(y)=\mathbf 1_B(z)\mathbf 1_J(t)\,\mu_t(E\cap I_n),
$$
which is Borel in $(z,t)$ by (i).
Let $\mathcal L_n$ be the family of Borel sets $Q\subset Z\times[0,T)\times\R$ such that $\mathbf 1_Q\in\mathcal H_n$.
It contains $Z\times[0,T)\times\R$ by the computation above, it is closed under nested differences because $\mathcal H_n$ is a vector space, and it is closed under increasing unions because $\mathcal H_n$ is closed under bounded pointwise limits.
Hence $\mathcal L_n$ is a $\lambda$-system.
The rectangles $B\times J\times E$ form a $\pi$-system contained in $\mathcal L_n$, since the intersection of two rectangles is a rectangle.
Every open subset of the separable metric space $Z\times[0,T)\times\R$ is a countable union of products of relatively open sets, so these rectangles generate its Borel $\sigma$-algebra, and Dynkin's $\pi$--$\lambda$ theorem gives $\mathbf 1_Q\in\mathcal H_n$ for every Borel set $Q\subset Z\times[0,T)\times\R$.
By linearity, $\mathcal H_n$ contains every Borel simple function.
Every bounded Borel function is the uniform limit of a sequence of Borel simple functions bounded by a common constant, so $\mathcal H_n$ contains every bounded Borel function.
Since $f(z,t,y)=f(z,t,y)\mathbf 1_{I_n}(y)$ for the function $f$ of the statement, this proves (ii).

(iii) Let $Q\subset\R\times[0,T)$ be a Borel set.
For each $n\in\N$, (ii) applied with $Z=\{0\}\subset\R$ to the function $(0,t,y)\mapsto\mathbf 1_Q(y,t)\mathbf 1_{I_n}(y)$ shows that $t\mapsto\mu_t(Q_t\cap I_n)$ is Borel.
As $n\to\infty$, $\mu_t(Q_t\cap I_n)$ increases to $\mu_t(Q_t)$, which is therefore a Borel function of $t$.
Clearly $\mu(\emptyset)=0$.
If $Q$ is the union of pairwise disjoint Borel sets $Q_k$, then $Q_t$ is the union of the pairwise disjoint sets $(Q_k)_t$, so $\mu_t(Q_t)=\sum_k\mu_t((Q_k)_t)$, and the monotone convergence theorem gives $\mu(Q)=\sum_k\mu(Q_k)$.
Thus $\mu$ is a measure.
The integral formula holds for $f=\mathbf 1_Q$ by the definition of $\mu$, hence for nonnegative Borel simple functions by linearity.
A nonnegative Borel function $f$ is the limit of an increasing sequence of nonnegative Borel simple functions, and the monotone convergence theorem, applied to the inner integrals and then to the integrals in $t$, shows that $t\mapsto\int_\R f(y,t)\,d\mu_t(y)$ is Borel and that the formula holds for $f$.
Finally, if $\mu(N)=0$, then $\int_0^T\mu_t(N_t)\,dt=0$ with a nonnegative Borel integrand, so $\mu_t(N_t)=0$ for almost every $t$.
\end{proof}

We also record a bound on the mass of compact sets, uniform on compact time intervals, which is used repeatedly below.

\begin{lemma}[Uniform local mass bound]\label{uniform-bound}
Let $(\rho,u)$ satisfy conditions (a) and (b) of Definition~\ref{weak-solution} on $[0,T)$, and let $K\subset\R$ and $I\subset[0,T)$ be compact. Then $\sup\limits_{t\in I}\rho_t(K)<\infty$ and $\sup\limits_{t\in I}|u(\cdot,t)\rho_t|(K)<\infty$.
\end{lemma}

\begin{proof}
Fix $\chi\in C_c(\R)$ with $\chi\ge\mathbf 1_K$, for instance a piecewise linear function equal to $1$ on the convex hull of $K$. Since $\rho_t\ge0$, $\rho_t(K)\le\int_\R\chi\,d\rho_t$, and the right-hand side is continuous on $[0,T)$ by (a), hence bounded on $I$. The second bound follows from $|u(\cdot,t)\rho_t|(K)\le\|u(\cdot,t)\|_{L^\infty(d\rho_t)}\rho_t(K)$ and \eqref{local-uniform-u}.
\end{proof}

The velocity is only assumed to be measurable with respect to $\rho_t$ at each fixed time. The next lemma shows that conditions (a) and (b) nevertheless provide a velocity that is jointly Borel in space and time and defines the same momentum at every time. It is used for the measurability of the weak formulation and, in Section~\ref{global-existence-uniqueness}, for the comparison with the formulation of the homogeneous theory.

\begin{lemma}[Borel representative of the velocity]\label{borel-velocity}
Let $(\rho,u)$ satisfy conditions (a) and (b) of Definition~\ref{weak-solution} on $[0,T)$.
Then there is a Borel function $h:\R\times[0,T)\to\R$ such that, for every $t\in[0,T)$, $h(\cdot,t)=u(\cdot,t)$ $\rho_t$-almost everywhere and $|h(x,t)|\le\|u(\cdot,t)\|_{L^\infty(d\rho_t)}$ for every $x\in\R$.
Consequently, if $\phi:\R\times[0,T)\to\R$ is a bounded Borel function such that $\phi(x,t)=0$ whenever $|x|\ge n$, for some $n\in\N$, and $k\in\{0,1,2\}$, then the function $t\mapsto\int_\R\phi(\cdot,t)\,u(\cdot,t)^k\,d\rho_t$ is Borel on $[0,T)$ and bounded on every compact subinterval of $[0,T)$.
\end{lemma}

\begin{proof}
For $j\in\N$ and $x\in\R$, put $I_j(x):=[x-1/j,x+1/j]$, and, for $t\in[0,T)$, let $r_j(x,t):=\rho_t(I_j(x))^{-1}\int_{I_j(x)}u(\cdot,t)\,d\rho_t$ if $\rho_t(I_j(x))>0$ and $r_j(x,t):=0$ otherwise.
Fix $n\in\N$. For $|x|<n$, the function $y\mapsto\mathbf 1_{I_j(x)}(y)$ vanishes for $|y|\ge n+1$.
By condition (a), Lemma~\ref{kernel}(ii), applied with $Z=(-n,n)$ to the bounded Borel function $(x,t,y)\mapsto\mathbf 1_{\{|y-x|\le1/j\}}\mathbf 1_{(-n-1,n+1)}(y)$ and to the weak-$*$ continuous families $(\rho_t)$ and $(u(\cdot,t)\rho_t)$, shows that $(x,t)\mapsto\rho_t(I_j(x))$ and $(x,t)\mapsto\int_{I_j(x)}u(\cdot,t)\,d\rho_t$ are Borel on $(-n,n)\times[0,T)$.
Since $n$ is arbitrary, $r_j$ is Borel on $\R\times[0,T)$, and $|r_j(x,t)|\le\|u(\cdot,t)\|_{L^\infty(d\rho_t)}$ for all $x$ and $t$.
Hence $h:=\limsup\limits_{j\to\infty}r_j$ is a real-valued Borel function with the stated bound.

Now fix $t\in[0,T)$. Since $u(\cdot,t)$ is bounded and $\rho_t$-measurable, it belongs to $L^1_{\mathrm{loc}}(\R,d\rho_t)$. By the Lebesgue--Besicovitch differentiation theorem for Radon measures on $\R$ \cite[Chapter~1]{EvansGariepy2015}, for $\rho_t$-almost every $x$ we have $\rho_t(I_j(x))>0$ for every $j$ and $r_j(x,t)\to u(x,t)$ as $j\to\infty$. Therefore $h(\cdot,t)=u(\cdot,t)$ $\rho_t$-almost everywhere. The argument is carried out at each fixed time, so no exceptional set of times arises.

For the last assertion, let $\phi$, $n$, and $k$ be as in the statement. For every $t$, $\int_\R\phi(\cdot,t)u(\cdot,t)^k\,d\rho_t=\int_\R\phi(\cdot,t)h(\cdot,t)^k\,d\rho_t$.
Let $I\subset[0,T)$ be a compact interval and let $U_I$ be the bound in \eqref{local-uniform-u}. On $\R\times I$, $h$ coincides with its truncation $h_I:=\max\{-U_I,\min\{U_I,h\}\}$, and Lemma~\ref{kernel}(ii), applied with $Z=\{0\}$ to the bounded Borel function $(0,t,y)\mapsto\phi(y,t)h_I(y,t)^k\mathbf 1_I(t)$, shows that the restriction of $t\mapsto\int_\R\phi(\cdot,t)u(\cdot,t)^k\,d\rho_t$ to $I$ is Borel. Exhausting $[0,T)$ by countably many such intervals gives the Borel measurability on $[0,T)$. Finally, for $t\in I$, the integral is bounded by $\sup|\phi|\,U_I^k\sup\limits_{s\in I}\rho_s([-n,n])$, which is finite by Lemma~\ref{uniform-bound}.
\end{proof}

In particular, the inner integrals in \eqref{weak1-rho} and \eqref{weak2} are Borel functions of $t$, bounded on compact subintervals of $(0,T)$ and vanishing outside the time support of the test functions, so the left-hand sides of these identities are meaningful. Lemma~\ref{borel-velocity} also shows that the energy pairings $t\mapsto\int_\R\varphi\,u(\cdot,t)^2\,d\rho_t$, $\varphi\in C_c(\R)$, are Borel, which is the reason why no measurability condition on them is included in Definition~\ref{weak-solution}.

Identity \eqref{weak1-rho} is the distributional form of $\partial_t\rho+\alpha(t)\,\partial_x(\rho u)=0$.
The homogeneous theory is formulated instead in terms of a cumulative potential, that is, a function $m$ with $\partial_xm=\rho$ which is transported by the flow.
To compare the two formulations, we show that \eqref{weak1-rho} alone produces such a potential, jointly in space and time. At the initial time, the potential is the cumulative distribution function of the initial density.

\begin{definition}[Initial cumulative distribution function]\label{m}
For nonnegative $\rho_0\in\M_{\mathrm{loc}}(\R)$, define
$$
m_0(x):=
\begin{cases}
\rho_0([0,x)), & x\ge0,\\
-\rho_0([x,0)), & x<0.
\end{cases}
$$
Then $m_0(0)=0$, $m_0$ is nondecreasing and left-continuous, $m_0\in BV_{\mathrm{loc}}(\R)$, $m_0'=\rho_0$ as Radon measures, and $m_0(b)-m_0(a)=\rho_0([a,b))$ for all $a<b$.
\end{definition}

At positive times, the primitive $p(\cdot,t)$ of $\rho_t$ normalized by $p(0,t)=0$ does not in general satisfy $\partial_tp=-\alpha u\rho$, since mass may cross the point $x=0$. In the construction below, this flux is accounted for by subtracting a function $C(t)$ of time alone.
For $T\in(0,\infty]$, we write $m\in L^\infty_{\mathrm{loc}}([0,T);BV_{\mathrm{loc}}(\R))$ if $m$ is a Borel function on $\R\times[0,T)$ such that
$$
\sup_{t\in I}\Big(\|m(\cdot,t)\|_{L^1(K)}+|\partial_xm(\cdot,t)|(K)\Big)<\infty
$$
for all compact sets $K\subset\R$ and $I\subset[0,T)$.

\begin{proposition}[Intrinsic construction of the cumulative potential]
\label{intrinsic-potential}
Let $(\rho,u)$ satisfy conditions (a), (b), and \eqref{weak1-rho} of Definition~\ref{weak-solution} on $[0,T)$. Then there exists a function
$$
m\in L^\infty_{\mathrm{loc}}([0,T);BV_{\mathrm{loc}}(\R))
$$
such that, in $\mathcal D'(\R\times(0,T))$,
\begin{equation}\label{potential-system}
\partial_xm=\rho,
\qquad
\partial_tm=-\alpha(t)u\rho.
\end{equation}
For every $t$, $m(\cdot,t)$ is nondecreasing and left-continuous.
If, in addition, $\rho_t\stackrel{*}{\rightharpoonup}\rho_0$ as $t\to0^+$, the additive constant may be fixed so that
\begin{equation}\label{potential-initial-trace}
m(x,t)\longrightarrow m_0(x)
\quad\text{as }t\to0^+
\end{equation}
at every continuity point of $m_0$, hence for almost every $x\in\R$.
The potential can be chosen in $C([0,T);L^1_{\mathrm{loc}}(\R))$ with $m(\cdot,0)=m_0$ when the initial normalization is imposed. Any two locally integrable potentials satisfying \eqref{potential-system} differ by a constant almost everywhere in space-time. Their normalized representatives in $C([0,T);L^1_{\mathrm{loc}}(\R))$, chosen left-continuous in $x$ at every time, coincide pointwise.
\end{proposition}

\begin{proof}
The proof has four steps. At each fixed time we take the primitive of $\rho_t$ normalized by a fixed weight $\theta$, and we first prove that it is jointly measurable and locally bounded. We then compute, from \eqref{weak1-rho}, the time derivative of its spatial averages.
This derivative differs from the one required by \eqref{potential-system} by a function of time alone, which we remove by subtracting its primitive $C(t)$. Finally, we identify the initial trace and prove uniqueness.

Fix $\theta\in C_c^\infty(\R)$ with $\theta\ge0$ and $\int_\R\theta\,dx=1$.
Define, for each fixed $t$,
$$
p(x,t):=
\begin{cases}
\rho_t([0,x)), &x\ge0,\\
-\rho_t([x,0)), &x<0,
\end{cases}
\qquad
\bar m(x,t):=p(x,t)-\int_\R\theta(z)p(z,t)\,dz.
$$

\medskip \noindent\textit{Step 1. Measurability and local bounds.} We first show that $p$ is jointly Borel.
We have
$$
p(x,t)=\int_\R \kappa(x,y)\,d\rho_t(y),
\qquad
\kappa(x,y):=
\mathbf 1_{\{x\ge0,\ 0\le y<x\}}
-\mathbf 1_{\{x<0,\ x\le y<0\}},
$$
where $\kappa$ is bounded and Borel on $\R^2$.
Fix $n\in\N$.
For $|x|<n$, the function $\kappa(x,\cdot)$ vanishes outside $(-n,n)$, so that $p(x,t)=\int_\R\kappa(x,y)\mathbf 1_{(-n,n)}(y)\,d\rho_t(y)$ for $|x|<n$.
By condition (a), Lemma~\ref{kernel}(ii), applied with $Z=(-n,n)$ and $f(x,t,y):=\kappa(x,y)\mathbf 1_{(-n,n)}(y)$, shows that $p$ is Borel on $(-n,n)\times[0,T)$.
Since $n$ is arbitrary, $p$ is Borel on $\R\times[0,T)$.

Next, choose $r>0$ such that $[-r,r]$ contains $0$ and $\supp\theta$. For $z\in\supp\theta$ we have $|p(z,t)|\le\rho_t([-r,r])<\infty$, so the integral $\int_\R\theta(z)p(z,t)\,dz$ is finite for every $t\in[0,T)$. The Borel $\sigma$-algebra of $\R\times[0,T)$ is the product of the Borel $\sigma$-algebras of $\R$ and $[0,T)$, so Tonelli's theorem, applied to the nonnegative Borel functions $\theta p^+$ and $\theta p^-$, where $p^\pm:=\max\{\pm p,0\}$, shows that this integral is a Borel function of $t$. Hence $\bar m$ is Borel on $\R\times[0,T)$.

For every fixed $t$, $\bar m(\cdot,t)$ is the unique left-continuous primitive of $\rho_t$ normalized by
\begin{equation}\label{smooth-normalization}
\int_\R\theta(x)\bar m(x,t)\,dx=0,
\end{equation}
so that $\partial_x\bar m(\cdot,t)=\rho_t$ for every $t$.

We now bound $\bar m$ locally. Let $K\subset\R$ and $I\subset[0,T)$ be compact, and choose a compact interval $\widehat K$ containing $K$, $0$, and $\supp\theta$.
By Lemma~\ref{uniform-bound}, $M_{\widehat K,I}:=\sup\limits_{t\in I}\rho_t(\widehat K)<\infty$. For $x\in \widehat K$ we have $|p(x,t)|\le M_{\widehat K,I}$, so $\big|\int_\R\theta(z)p(z,t)\,dz\big|\le M_{\widehat K,I}$ and
\begin{equation}\label{m-bar-bound}
\sup_{t\in I}\|\bar m(\cdot,t)\|_{L^1(K)}
\le 2|K|M_{\widehat K,I},
\qquad
\sup_{t\in I}|\partial_x\bar m(\cdot,t)|(K)\le M_{\widehat K,I}.
\end{equation}

\medskip \noindent\textit{Step 2. Time derivative of the spatial averages.}
We first express the spatial averages of $\bar m(\cdot,t)$ through $\rho_t$. For $\varphi\in C_c^\infty(\R)$, put $g:=\varphi-(\int_\R\varphi\,dx)\theta$ and $G(x):=\int_{-\infty}^xg(s)\,ds$. Since $\int_\R g\,dx=0$, $G\in C_c^\infty(\R)$, and an integration by parts at fixed $t$, combined with \eqref{smooth-normalization}, gives
\begin{equation}\label{primitive-pairing}
\int_\R\varphi(x)\bar m(x,t)\,dx=-\int_\R G(x)\,d\rho_t(x)=:F_\varphi(t).
\end{equation}
Since $G$ is fixed, $F_\varphi$ is continuous on $[0,T)$ by (a). In particular, $F_\theta\equiv0$ by \eqref{smooth-normalization}.

To differentiate these averages in time, fix $\varphi\in C_c^\infty(\R)$ and test \eqref{weak1-rho} with $\Xi(x,t):=\zeta(t)\varphi(x)$, where $\zeta\in C_c^\infty(0,T)$ is arbitrary.
Applying \eqref{primitive-pairing} with $\varphi$ replaced by $\varphi'$, and noting that $\int_\R\varphi'\,dx=0$, so that $g=\varphi'$ and $G=\varphi$, we have $F_{\varphi'}(t)=-\int_\R\varphi\,d\rho_t$, and the tested identity becomes
$$
0=-\int_0^T\zeta'(t)F_{\varphi'}(t)\,dt+\int_0^T\alpha(t)\zeta(t)\Big[\int_\R\varphi'u\,d\rho_t\Big]dt
\qquad\text{for every }\zeta\in C_c^\infty(0,T).
$$
Thus the distributional derivative of $F_{\varphi'}$ on $(0,T)$ is the function $t\mapsto-\alpha(t)\int_\R\varphi'u\,d\rho_t$, which is continuous by (a) and \ref{alpha}.
Since a continuous function whose distributional derivative is continuous is of class $C^1$,
\begin{equation}\label{F-derivative}
F_{\varphi'}'(t)=-\alpha(t)\int_\R\varphi'(x)u(x,t)\,d\rho_t(x)
\qquad\text{for every }t\in(0,T).
\end{equation}
Every $g\in C_c^\infty(\R)$ with $\int_\R g\,dx=0$ is the derivative of its compactly supported primitive, so \eqref{F-derivative} holds with $\varphi'$ replaced by any such $g$. For a general $\psi\in C_c^\infty(\R)$, write $\psi=g+(\int_\R\psi\,dx)\theta$ with $g:=\psi-(\int_\R\psi\,dx)\theta$. Since $F_\theta\equiv0$ and $\varphi\mapsto F_\varphi$ is linear, $F_\psi=F_g$, and therefore
\begin{equation}\label{F-derivative-general}
F_\psi'(t)=-\alpha(t)\int_\R\psi\,u\,d\rho_t+\Big(\int_\R\psi\,dx\Big)b(t),
\qquad
b(t):=\alpha(t)\int_\R\theta(x)u(x,t)\,d\rho_t(x),
\end{equation}
for every $t\in(0,T)$. The function $b$ is continuous on $[0,T)$ by (a) and \ref{alpha}.

\medskip \noindent\textit{Step 3. Definition of $m$ and proof of \eqref{potential-system}.}
The last term in \eqref{F-derivative-general} depends on $\psi$ only through $\int_\R\psi\,dx$, and it is removed by subtracting a primitive of $b$.
Let $C_0\in\R$ be an arbitrary constant, put $C(t):=C_0+\int_0^tb(s)\,ds$, so that $C\in C^1([0,T))$ with $C'=b$, and define
\begin{equation}\label{potential-def}
m(x,t):=\bar m(x,t)-C(t).
\end{equation}
For every $\psi\in C_c^\infty(\R)$, the function $t\mapsto\int_\R\psi m\,dx=F_\psi(t)-C(t)\int_\R\psi\,dx$ is of class $C^1$ on $(0,T)$, and by \eqref{F-derivative-general}
$$
\frac{d}{dt}\int_\R\psi(x)m(x,t)\,dx
=F_\psi'(t)-C'(t)\int_\R\psi\,dx
=-\alpha(t)\int_\R\psi\,u\,d\rho_t.
$$
The constant $C_0$ disappears upon differentiation, so its choice does not affect the following computations.
Hence $\partial_xm=\rho$ and $\partial_tm=-\alpha u\rho$ hold when tested against every product $\zeta(t)\psi(x)$ with $\zeta\in C_c^\infty(0,T)$ and $\psi\in C_c^\infty(\R)$.

It remains to pass from products to general test functions. The joint measurability of $\bar m$, the bounds \eqref{m-bar-bound}, and the continuity of $C$ show that $m$ is locally integrable on $\R\times[0,T)$ and belongs to $L^\infty_{\mathrm{loc}}([0,T);BV_{\mathrm{loc}}(\R))$.
Thus, $m$, $\partial_xm$, and $\partial_tm$ are well-defined distributions, and $\rho$ and $\alpha u\rho$ define locally finite measures on $\R\times(0,T)$. For the latter, use the weak-$*$ continuous signed kernel $u(\cdot,t)\rho_t$, Lemma~\ref{kernel}, and the local total-variation bound $\max_I\alpha\,U_I M_{\widehat K,I}$, where $U_I$ is the bound in \eqref{local-uniform-u}. Since both sides of \eqref{potential-system}, tested against $\Xi\in C_c^\infty(\R\times(0,T))$, involve only $\Xi$ and its first-order derivatives, it suffices to approximate $\Xi$ by finite sums of products in this sense.
Let $\supp\Xi\subset(-a,a)\times(c,d)$ with $a>0$ and $[c,d]\subset(0,T)$, and put $Q:=[-a,a]\times[c,d]$.
By the Weierstrass approximation theorem, there are polynomials $q_n$ in $(x,t)$ with $q_n\to\partial_x\partial_t\Xi$ uniformly on $Q$.
Since $\Xi$ and its first-order derivatives vanish on $\{-a\}\times[c,d]$ and on $[-a,a]\times\{c\}$, the polynomials
$$
P_n(x,t):=\int_{-a}^x\!\!\int_c^tq_n(y,s)\,ds\,dy
$$
satisfy $P_n\to\Xi$, $\partial_xP_n\to\partial_x\Xi$, and $\partial_tP_n\to\partial_t\Xi$ uniformly on $Q$.
Choose $\chi\in C_c^\infty((-a,a))$ and $\vartheta\in C_c^\infty((c,d))$ such that $\chi(x)\vartheta(t)=1$ on $\supp\Xi$, and put $\Xi_n(x,t):=\chi(x)\vartheta(t)P_n(x,t)$.
Each $\Xi_n$ is a finite sum of products $\zeta(t)\psi(x)$ as above, $\supp\Xi_n\subset Q$, and, since $\Xi=\chi\vartheta\,\Xi$, the functions $\Xi_n-\Xi=\chi\vartheta\,(P_n-\Xi)$ converge to zero uniformly on $Q$ together with their first-order derivatives.
By the local bounds just mentioned, each side of \eqref{potential-system}, tested against $\Xi_n$, converges to the same expression tested against $\Xi$. This proves \eqref{potential-system} in $\mathcal D'(\R\times(0,T))$.

\medskip \noindent\textit{Step 4. Initial trace and uniqueness.}
Assume $\rho_t\stackrel{*}{\rightharpoonup}\rho_0$ as $t\to0^+$.
Let $\bar m_0$ be the left-continuous primitive of $\rho_0$ normalized by $\int_\R\theta\bar m_0\,dx=0$. By \eqref{primitive-pairing}, applied also at $t=0$, $\int_\R\varphi\,\bar m(\cdot,t)\,dx\to\int_\R\varphi\,\bar m_0\,dx$ as $t\to0^+$, for every $\varphi\in C_c^\infty(\R)$. The functions $\bar m(\cdot,t)$ are nondecreasing and, by \eqref{m-bar-bound}, locally uniformly bounded for $t$ near $0$.
By Helly's selection theorem, every sequence $t_n\to0^+$ has a subsequence along which $\bar m(\cdot,t_n)$ converges pointwise to a nondecreasing function, and the distributional convergence identifies this limit with $\bar m_0$ at every continuity point of $\bar m_0$.
Hence $\bar m(x,t)\to\bar m_0(x)$ as $t\to0^+$ at every continuity point of $\bar m_0$. Since $m_0$ and $\bar m_0$ are left-continuous primitives of the same measure, $m_0-\bar m_0=d_0$ for a constant $d_0$. Since $C_0$ was arbitrary in Step~3, we may choose $C_0:=-d_0$. Then $C(t)\to C_0$ as $t\to0^+$, and $m(x,t)=\bar m(x,t)-C(t)\to\bar m_0(x)+d_0=m_0(x)$ at every continuity point of $m_0$, which is \eqref{potential-initial-trace}.

For uniqueness, let $m^1$ and $m^2$ satisfy \eqref{potential-system}. Then $d:=m^1-m^2$ satisfies $\partial_xd=\partial_td=0$ in $\mathcal D'(\R\times(0,T))$. A distribution on a connected open set whose first-order derivatives vanish is constant, so $m^1-m^2$ is constant almost everywhere in space-time.

For the constructed potential \eqref{potential-def}, \eqref{primitive-pairing} and continuity of $C$ give continuity of $t\mapsto\int\varphi m\,dx$ for every $\varphi\in C_c^\infty(\R)$. Fix $t_*\in[0,T)$ and a sequence $t_j\to t_*$. On each compact spatial interval, the functions $m(\cdot,t_j)$ are uniformly bounded and nondecreasing. Helly's theorem and a diagonal extraction yield a subsequence converging almost everywhere and in $L^1$ on each compact interval. Continuity of the spatial pairings identifies the limit with $m(\cdot,t_*)$ almost everywhere. Every subsequence has the same possible $L^1_{\mathrm{loc}}$ limit, so the entire sequence converges. Thus $m\in C([0,T);L^1_{\mathrm{loc}}(\R))$.
For two potentials in this time-continuous class, their almost-everywhere constant difference holds in $L^1_{\mathrm{loc}}$ at every time by continuity. Their common initial trace makes the constant zero. Left continuity in $x$ upgrades equality almost everywhere in $x$ to equality at every $x$, since a nonzero difference at a point would, by left continuity, persist on an interval of positive length. Without the time-continuity convention, uniqueness is only almost everywhere in space-time.
\end{proof}

For the homogeneous system, Proposition~\ref{intrinsic-potential} takes the form in which it will be compared with the homogeneous formulation in Section~\ref{global-existence-uniqueness}.

\begin{corollary}[Intrinsic and cumulative formulations]
\label{intrinsic-cumulative}
Let $(\rho,w)$ be a weak solution of \eqref{H} on $[0,T)$ in the sense of Definition~\ref{weak-solution}.
Then there is a potential $m$, unique almost everywhere once \eqref{potential-initial-trace} is imposed and unique pointwise with the time-continuous, spatially left-continuous convention of Proposition~\ref{intrinsic-potential}, such that
$$
\partial_xm=\rho,
\qquad
\partial_\tau m+w\, \partial_xm=0
$$
in $\mathcal D'(\R\times(0,T))$, where $w\,\partial_xm$ denotes the measure $w\rho$.
Conversely, if a function $m\in L^\infty_{\mathrm{loc}}([0,T);BV_{\mathrm{loc}}(\R))$ with $\partial_xm=\rho$ satisfies $\partial_\tau m+w\,\partial_xm=0$, then differentiation in $x$ gives \eqref{weak1-H}.
\end{corollary}

\begin{proof}
The first statement is Proposition~\ref{intrinsic-potential} with $\alpha\equiv1$. For the converse, $\partial_\tau m=-w\rho$ in $\mathcal D'(\R\times(0,T))$, and applying $\partial_x$, which commutes with $\partial_\tau$ on distributions, gives $\partial_\tau\rho+\partial_x(w\rho)=0$ in $\mathcal D'(\R\times(0,T))$, which is \eqref{weak1-H}.
\end{proof}

\subsection{Weak formulation with \texorpdfstring{$C^1$}{C1} test functions}

The equivalence of weak solutions is obtained by transporting test functions through the time change $\tau=A(t)$.
If $\Phi$ is a smooth function with compact support in $\R\times(0,A(T))$, the composition $\Phi(x,A(t))$ is of class $C^1$, but in general not of class $C^2$.
Indeed, since $e^{-S}\in C^1$ and $e^{-S}>0$, the derivative $A'=\alpha e^{-S}$ is of class $C^1$ if and only if $\alpha$ is, so $A\in C^2$ if and only if $\alpha\in C^1$, while Assumption~\ref{coefficients} only requires $\alpha$ to be continuous.
Since Definition~\ref{weak-solution} requires \eqref{weak1-rho} and \eqref{weak2} only for smooth test functions, we first show that the weak formulation remains valid for $C^1$ test functions, and then we describe the correspondence of test functions induced by the time change.
The extension uses Lemma~\ref{uniform-bound} and the following continuity property of weak-$*$ continuous families of measures.

\begin{lemma}[Continuity of parametrized integrals]\label{param-continuity}
Let $J\subset\R$ be an interval and let $t\mapsto\mu_t\in\M_{\mathrm{loc}}(\R)$ be weak-$*$ continuous on $J$, with $\sup\limits_{t\in I}|\mu_t|(K)<\infty$ for all compact $K\subset\R$ and $I\subset J$.
Then, for every $\phi\in C_c(\R\times J)$, the function $t\mapsto\int_\R\phi(x,t)\,d\mu_t(x)$ is continuous on $J$.
\end{lemma}

\begin{proof}
Let $\supp\phi\subset K\times I$ with $K$, $I$ compact, and let $s,t\in J$.
Then,
$$
\Big|\int_\R\phi(\cdot,t)\,d\mu_t-\int_\R\phi(\cdot,s)\,d\mu_s\Big|
\le\|\phi(\cdot,t)-\phi(\cdot,s)\|_\infty\,|\mu_t|(K)
+\Big|\int_\R\phi(\cdot,s)\,d(\mu_t-\mu_s)\Big|.
$$
As $t\to s$, the first term tends to zero by the uniform continuity of $\phi$ and the bound on $|\mu_t|(K)$ for $t$ near $s$, and the second by the weak-$*$ continuity of $t\mapsto\mu_t$.
\end{proof}

\begin{lemma}[Extension to $C^1$ test functions]\label{extension}
Let $(\rho,u)$ satisfy conditions (a) and (b) of Definition~\ref{weak-solution} on $[0,T)$.
If \eqref{weak1-rho} and \eqref{weak2} hold for all $\Xi,\Psi\in C_c^\infty(\R\times(0,T))$, then they hold for all $\Xi,\Psi\in C_c^1(\R\times(0,T))$.
\end{lemma}

\begin{proof}
Let $\Xi\in C_c^1(\R\times(0,T))$, with $\supp\Xi\subset K_x\times I$ for compact $K_x\subset\R$ and $I\subset(0,T)$, and let $\Xi_\varepsilon:=\Xi*j_\varepsilon$ be a standard mollification in $(x,t)$.
For $\varepsilon>0$ sufficiently small, $\Xi_\varepsilon\in C_c^\infty(\R\times(0,T))$ and $\supp\Xi_\varepsilon\subset \widehat K_x\times \widehat I$ for fixed compact sets $\widehat K_x\supset K_x$ and $\widehat I\subset(0,T)$ containing $I$.
Since $\Xi\in C^1$, $\Xi_\varepsilon\to\Xi$ and $\nabla\Xi_\varepsilon\to\nabla\Xi$ uniformly on $\R\times(0,T)$ as $\varepsilon\to0^+$. Put
$$
\bar\rho:=\sup_{t\in \widehat I}\rho_t(\widehat K_x),\qquad
\bar u:=\sup_{t\in \widehat I}\|u(\cdot,t)\|_{L^\infty(d\rho_t)},\qquad
\|\alpha\|_{\widehat I}:=\max_{\widehat I}\alpha,\qquad
\|\sigma\|_{\widehat I}:=\max_{\widehat I}\sigma,
$$
which are finite by Lemma~\ref{uniform-bound}, \eqref{local-uniform-u}, and Assumption~\ref{coefficients}. Then,
$$
\Big|\iint_{\R\times(0,T)}\big(\partial_t\Xi_\varepsilon-\partial_t\Xi\big)\,d\rho_t\,dt\Big|
\le|\widehat I|\,\|\partial_t\Xi_\varepsilon-\partial_t\Xi\|_\infty\,\bar\rho\longrightarrow0
$$
and
$$
\Big|\iint_{\R\times(0,T)}\alpha(t)\big(\partial_x\Xi_\varepsilon-\partial_x\Xi\big)u\,d\rho_t\,dt\Big|
\le|\widehat I|\,\|\alpha\|_{\widehat I}\,\|\partial_x\Xi_\varepsilon-\partial_x\Xi\|_\infty\,\bar u\,\bar\rho\longrightarrow0
$$
as $\varepsilon\to0^+$, where the integrands in $t$ are continuous by Lemma~\ref{param-continuity}.
Since \eqref{weak1-rho} holds for each $\Xi_\varepsilon$, it holds for $\Xi$.

Similarly, for $\Psi\in C_c^1(\R\times(0,T))$ and its mollifications $\Psi_\varepsilon$, with the same notation,
$$
\Big|\iint_{\R\times(0,T)}\big(\partial_t\Psi_\varepsilon-\partial_t\Psi\big)u\,d\rho_t\,dt\Big|
\le|\widehat I|\,\|\partial_t\Psi_\varepsilon-\partial_t\Psi\|_\infty\,\bar u\,\bar\rho,
$$
$$
\Big|\iint_{\R\times(0,T)}\alpha(t)\big(\partial_x\Psi_\varepsilon-\partial_x\Psi\big)u^2\,d\rho_t\,dt\Big|
\le|\widehat I|\,\|\alpha\|_{\widehat I}\,\|\partial_x\Psi_\varepsilon-\partial_x\Psi\|_\infty\,\bar u^2\,\bar\rho,
$$
$$
\Big|\iint_{\R\times(0,T)}\sigma(t)\big(\Psi_\varepsilon-\Psi\big)u\,d\rho_t\,dt\Big|
\le|\widehat I|\,\|\sigma\|_{\widehat I}\,\|\Psi_\varepsilon-\Psi\|_\infty\,\bar u\,\bar\rho,
$$
where we used $|u(\cdot,t)|\le\bar u$ $\rho_t$-almost everywhere for $t\in \widehat I$ and the measurability in $t$ given by Lemma~\ref{borel-velocity}.
All three bounds tend to zero as $\varepsilon\to0^+$, and \eqref{weak2} passes to the limit.
\end{proof}

With Lemma~\ref{extension} available, the correspondence of test functions induced by the time change can be used in both directions.

\begin{lemma}\label{testfunctions}
Let Assumption~\ref{coefficients} hold and let $T\in(0,\infty]$, with $A(\infty):=A_\infty$.
The correspondence $\Phi(x,\tau)\longmapsto\varphi(x,t):=\Phi(x,A(t))$ is a bijection from $C_c^1(\R\times(0,A(T)))$ onto $C_c^1(\R\times(0,T))$, with inverse $\varphi\mapsto\Phi(x,\tau):=\varphi(x,A^{-1}(\tau))$.
For every such pair,
$$
\partial_x\varphi(x,t)=\partial_x\Phi(x,A(t))
\qquad \mbox{ and } \qquad
\partial_t\varphi(x,t)=A'(t)\partial_\tau\Phi(x,A(t)).
$$
\end{lemma}

\begin{proof}
By Lemma~\ref{properties-A}, $A$ is a $C^1$ diffeomorphism of $(0,T)$ onto $(0,A(T))$, with $C^1$ inverse. The composition of a $C^1$ function with the $C^1$ map $(x,t)\mapsto(x,A(t))$ is $C^1$, and the chain rule gives the stated identities.
Since $(x,t)\mapsto(x,A(t))$ is a homeomorphism of $\R\times(0,T)$ onto $\R\times(0,A(T))$, it maps compact sets onto compact sets in both directions, so compactness of supports is preserved. The same argument applied to $A^{-1}$ shows that the inverse correspondence is well defined and inverse to the first.
\end{proof}

\subsection{Equivalence of weak solutions}

With Lemmas~\ref{extension} and~\ref{testfunctions} at hand, the computation of Theorem~\ref{transformation} can be repeated at the level of the weak formulation.

\begin{theorem} 
\label{weak-equivalence}
Let Assumption~\ref{coefficients} hold and let $T\in(0,\infty]$.
Let $(\rho,u)$ be a weak solution of \eqref{NA} on $[0,T)$, and let $(\tilde\rho,w):=\mathcal T_T(\rho,u)$.
Then $(\tilde\rho,w)$ is a weak solution of \eqref{H} on $[0,A(T))$ with the same initial data.
Conversely, if $(\tilde\rho,w)$ is a weak solution of \eqref{H} on $[0,A(T))$, then $\mathcal T_T^{-1}(\tilde\rho,w)$ is a weak solution of \eqref{NA} on $[0,T)$ with the same initial data.
Consequently, $\mathcal T_T$ maps the set of weak solutions of \eqref{NA} on $[0,T)$ bijectively onto the set of weak solutions of \eqref{H} on $[0,A(T))$ with the same initial data.
\end{theorem}

\begin{proof}
We prove the direct implication and then the converse.
The converse requires a separate computation, because the inverse correspondence of test functions involves $A^{-1}$ and the derivative of $S\circ A^{-1}$.

We begin with conditions (a) and (b) for $(\tilde\rho,w)$.
For each $\tau$, $\tilde\rho_\tau=\rho_{A^{-1}(\tau)}$ is a nonnegative Radon measure, and $w(\cdot,\tau)$ is a bounded measurable function with respect to $\tilde\rho_\tau$, being a positive constant multiple of $u(\cdot,A^{-1}(\tau))$.
Weak-$*$ continuity of $\tau\mapsto\tilde\rho_\tau$ and of $\tau\mapsto w(\cdot,\tau)\tilde\rho_\tau =e^{S(A^{-1}(\tau))}u(\cdot,A^{-1}(\tau))\rho_{A^{-1}(\tau)}$ follows from condition (a) for $(\rho,u)$, the continuity of $A^{-1}$, and the continuity of the scalar factor $e^{S\circ A^{-1}}$.
If $I\subset[0,A(T))$ is compact, then $A^{-1}(I)\subset[0,T)$ is compact, and
$$
\sup_{\tau\in I}\|w(\cdot,\tau)\|_{L^\infty(d\tilde\rho_\tau)}
\le\Big(\max_{A^{-1}(I)}e^{S}\Big)\sup_{t\in A^{-1}(I)}\|u(\cdot,t)\|_{L^\infty(d\rho_t)}<\infty.
$$

Throughout the proof we use the substitution rule
\begin{equation}\label{substitution}
\int_0^TA'(t)\,G(A(t))\,dt=\int_0^{A(T)}G(\tau)\,d\tau,
\end{equation}
valid for every bounded Lebesgue measurable function $G$ vanishing outside a compact subset of $(0,A(T))$, since $A$ is a $C^1$ diffeomorphism of $(0,T)$ onto $(0,A(T))$ whose inverse is locally Lipschitz.

We now verify the mass equation. Let $\Phi\in C_c^\infty(\R\times(0,A(T)))$ and set $\Xi(x,t):=\Phi(x,A(t))$, which belongs to $C_c^1(\R\times(0,T))$ by Lemma~\ref{testfunctions}.
By Lemma~\ref{extension}, \eqref{weak1-rho} holds for $\Xi$. Using $\partial_t\Xi=A'(t)\partial_\tau\Phi(x,A(t))$, $\partial_x\Xi=\partial_x\Phi(x,A(t))$, and $\alpha(t)u=\alpha(t)e^{-S(t)}w(\cdot,A(t))=A'(t)w(\cdot,A(t))$ $\rho_t$-almost everywhere, we obtain $\int_0^T A'(t)\,G(A(t))\,dt=0$ where $G(\tau):=\int_\R\big(\partial_\tau\Phi(x,\tau)+w(x,\tau)\partial_x\Phi(x,\tau)\big)\,d\tilde\rho_\tau(x)$.
By Lemma~\ref{param-continuity}, applied to $\tilde\rho_\tau$ and to $w(\cdot,\tau)\tilde\rho_\tau$, whose local total variations are bounded on compact time intervals by Lemma~\ref{uniform-bound} and the bound on $w$ obtained above, $G$ is continuous on $(0,A(T))$ with compact support.
The substitution rule \eqref{substitution} therefore gives
$$
\int_0^{A(T)}\!\!\int_\R \big(\partial_\tau\Phi+w\partial_x\Phi\big)\,d\tilde\rho_\tau\,d\tau
=\int_0^TA'(t)\,G(A(t))\,dt=0,
$$
which is \eqref{weak1-H} for $(\tilde\rho,w)$.

The momentum equation is verified in the same way.
Let $Z\in C_c^\infty(\R\times(0,A(T)))$ and set $\Psi(x,t):=e^{S(t)}Z(x,A(t))$. Since $S\in C^1$ by Lemma~\ref{properties-A}(i), $\Psi\in C_c^1(\R\times(0,T))$ by the argument of Lemma~\ref{testfunctions}, and Lemma~\ref{extension} applies. We have
$$
\partial_t\Psi=\sigma(t)e^{S(t)}Z(x,A(t))+e^{S(t)}A'(t)\partial_\tau Z(x,A(t)),\qquad
\partial_x\Psi=e^{S(t)}\partial_xZ(x,A(t)).
$$
Write $w$ for $w(\cdot,A(t))$, so that $e^{S(t)}u=w$ $\rho_t$-almost everywhere. Then, $\partial_t\Psi u=\sigma(t)wZ+A'(t)w\,\partial_\tau Z$, $\sigma(t)\Psi u=\sigma(t)wZ$, and $\alpha(t)\partial_x\Psi u^2=\alpha(t)e^{-S(t)}w^2\partial_xZ=A'(t)w^2\partial_xZ$, with $Z$, $\partial_\tau Z$, $\partial_xZ$ evaluated at $(x,A(t))$. Substituting into \eqref{weak2}, the two terms $\sigma(t)wZ$ cancel exactly, and
$$
0=\int_0^T A'(t)\,G(A(t))\,dt,
\qquad
G(\tau):=\int_\R\big(w\,\partial_\tau Z+w^2\partial_xZ\big)\,d\tilde\rho_\tau.
$$
By Lemma~\ref{borel-velocity}, applied to $(\tilde\rho,w)$, and the bounds obtained at the beginning of the proof, $G$ is a bounded Borel function that vanishes outside a compact subset of $(0,A(T))$, so \eqref{substitution} gives
$$
0=\int_0^{A(T)}\!\!\int_\R\big(w\,\partial_\tau Z+w^2\partial_xZ\big)\,d\tilde\rho_\tau\,d\tau,
$$
which is \eqref{weak2-H} for $(\tilde\rho,w)$.

It remains to verify condition (d). Since $A^{-1}(\tau)\to0$ as $\tau\to0^+$, $\tilde\rho_\tau=\rho_{A^{-1}(\tau)}\stackrel{*}{\rightharpoonup}\rho_0$.
Moreover, $w(\cdot,\tau)\tilde\rho_\tau=e^{S(A^{-1}(\tau))}u(\cdot,A^{-1}(\tau))\rho_{A^{-1}(\tau)}$, where the scalar factor tends to $1$ and the measures converge weakly-$*$ to $u_0\rho_0$ with locally bounded total variation, and hence $w(\cdot,\tau)\tilde\rho_\tau\stackrel{*}{\rightharpoonup}u_0\rho_0$.
This completes the proof that $(\tilde\rho,w)$ is a weak solution of \eqref{H} on $[0,A(T))$.

We turn to the converse.
Let $(\tilde\rho,w)$ be a weak solution of \eqref{H} on $[0,A(T))$ and let $(\rho,u):=\mathcal T_T^{-1}(\tilde\rho,w)$, that is, $\rho_t=\tilde\rho_{A(t)}$ and $u(\cdot,t)=e^{-S(t)}w(\cdot,A(t))$.
Conditions (a) and (b) for $(\rho,u)$ follow as in the first part of the proof, with the roles of $A$ and $A^{-1}$ exchanged. Since $(\tilde\rho,w)$ satisfies (a) and (b), Lemma~\ref{extension}, applied to \eqref{H}, shows that its weak formulation holds for all test functions in $C_c^1(\R\times(0,A(T)))$.

For the mass equation, let $\Xi\in C_c^\infty(\R\times(0,T))$ and set $\Phi(x,\tau):=\Xi(x,A^{-1}(\tau))$, which belongs to $C_c^1(\R\times(0,A(T)))$ by Lemma~\ref{testfunctions}, so that \eqref{weak1-H} holds for $\Phi$, that is,
$$
\int_0^{A(T)}\!\!\int_\R\big(\partial_\tau\Phi+w\partial_x\Phi\big)\,d\tilde\rho_\tau\,d\tau=0.
$$
By the chain rule and $(A^{-1})'(A(t))=1/A'(t)$, we have $\partial_\tau\Phi(x,A(t))=\partial_t\Xi(x,t)/A'(t)$ and $\partial_x\Phi(x,A(t))=\partial_x\Xi(x,t)$.
Moreover, $w(\cdot,A(t))=e^{S(t)}u(\cdot,t)$ $\rho_t$-almost everywhere and $A'(t)e^{S(t)}=\alpha(t)$.
The substitution rule \eqref{substitution} therefore gives
$$
0=\int_0^T\!\!\int_\R\Big(\frac{\partial_t\Xi}{A'(t)}+e^{S(t)}u\,\partial_x\Xi\Big)A'(t)\,d\rho_t\,dt
=\int_0^T\!\!\int_\R\partial_t\Xi\,d\rho_t\,dt
+\int_0^T\!\!\int_\R\alpha(t)\partial_x\Xi u\,d\rho_t\,dt,
$$
which is \eqref{weak1-rho} for $(\rho,u)$.

For the momentum equation, let $\Psi\in C_c^\infty(\R\times(0,T))$ and set
$$
Z(x,\tau):=e^{-S(A^{-1}(\tau))}\Psi(x,A^{-1}(\tau)),
$$
so that $Z(x,A(t))=e^{-S(t)}\Psi(x,t)$. Since $S\circ A^{-1}$ and $\Psi(\cdot,A^{-1}(\cdot))$ are of class $C^1$, $Z\in C_c^1(\R\times(0,A(T)))$, and \eqref{weak2-H} holds for $Z$, that is,
$$
\int_0^{A(T)}\!\!\int_\R\big(w\,\partial_\tau Z+w^2\partial_xZ\big)\,d\tilde\rho_\tau\,d\tau=0.
$$
Using $S'=\sigma$ and $(A^{-1})'(A(t))=1/A'(t)$,
$$
\partial_xZ(x,A(t))=e^{-S(t)}\partial_x\Psi(x,t),\qquad
\partial_\tau Z(x,A(t))=\frac{e^{-S(t)}}{A'(t)}\big[\partial_t\Psi(x,t)-\sigma(t)\Psi(x,t)\big].
$$
Hence, $\rho_t$-almost everywhere,
\begin{align*}
w(\cdot,A(t))\,\partial_\tau Z(\cdot,A(t))\,A'(t)&=u\big(\partial_t\Psi-\sigma(t)\Psi\big),\\
w(\cdot,A(t))^2\,\partial_xZ(\cdot,A(t))\,A'(t)&=e^{S(t)}A'(t)\,u^2\partial_x\Psi=\alpha(t)u^2\partial_x\Psi.
\end{align*}
The substitution rule \eqref{substitution} gives
$$
\int_0^T\!\!\int_\R\Big(u\partial_t\Psi+\alpha(t)u^2\partial_x\Psi-\sigma(t)u\Psi\Big)\,d\rho_t\,dt=0,
$$
which is \eqref{weak2} for $(\rho,u)$. Finally, as $t\to0^+$, $\rho_t=\tilde\rho_{A(t)}\stackrel{*}{\rightharpoonup}\rho_0$ because $A(t)\to0$, and $u(\cdot,t)\rho_t=e^{-S(t)}w(\cdot,A(t))\tilde\rho_{A(t)}\stackrel{*}{\rightharpoonup}u_0\rho_0$ because $e^{-S(t)}\to1$.
\end{proof}

Weak solutions of the pressureless system are, however, not unique without an admissibility criterion. Selecting the physically relevant one is the purpose of the next section.

\section{Entropy solutions and equivalence of the Cauchy problems}
\label{entropy-cauchy-section}

We introduce admissibility conditions for \eqref{NA} and prove that they correspond, under $\mathcal T$, to those of the homogeneous system. The two conditions are an Oleinik-type one-sided estimate and the continuity of the kinetic energy at the initial time.
Both are formulated intrinsically, that is, in terms of $(\rho,u)$, $\alpha$, and $\sigma$ only, and their correspondence with the homogeneous conditions is a matter of tracking the effective time $A(t)$ and the factor $e^{S(t)}$.
Combined with Theorem~\ref{weak-equivalence}, this correspondence shows that the two Cauchy problems are equivalent.

\subsection{Admissibility conditions}
\label{admissibility-conditions}

The Oleinik condition compares the velocity at pairs of points, whereas $u(\cdot,t)$ is only defined $\rho_t$-almost everywhere. We must therefore specify which representative of $u(\cdot,t)$ the condition refers to.
This is done intrinsically, since defining the representative by transporting the velocity of the homogeneous problem would make the notion of entropy solution depend on the transformation itself.
The bound in \eqref{intrinsic-oleinik} is dictated by free particles. Two particles with initial positions $x_1<x_2$ and velocities $u_1,u_2$ are at positions $x_i+u_iA(t)$ with velocities $u_ie^{-S(t)}$, so that, before they collide, the ratio of the velocity difference to the distance is at most $e^{-S(t)}/A(t)$.

\begin{definition}[Admissible representative]\label{admissible-representative}
Let $(\rho,u)$ be a weak solution of \eqref{NA} on $[0,T)$ and let $t\in(0,T)$. A bounded Borel function $\hat u(\cdot,t):\R\to\R$ is an \emph{admissible representative} of $u(\cdot,t)$ if $\hat u(\cdot,t)=u(\cdot,t)$ $\rho_t$-almost everywhere and
\begin{equation}\label{intrinsic-oleinik}
\frac{\hat u(x_2,t)-\hat u(x_1,t)}{x_2-x_1}\le\frac{e^{-S(t)}}{A(t)}
\qquad\text{for all }x_1<x_2.
\end{equation}
We say that $u$ satisfies the \emph{intrinsic Oleinik condition} if $u(\cdot,t)$ has an admissible representative for almost every $t\in(0,T)$.
\end{definition}

Only the existence of an admissible representative is required. On an open set where $\rho_t$ vanishes, the weak formulation does not see the values of $\hat u(\cdot,t)$, which are then constrained only by \eqref{intrinsic-oleinik}.
The definition imposes the inequality at every pair of points, including points of vacuum. The next lemma shows that this requirement is equivalent to imposing the inequality only for pairs of points outside a $\rho_t$-null set.

\begin{lemma}[Equivalent formulations of the one-sided bound]\label{oleinik-forms}
Let $\mu$ be a nonnegative Radon measure on $\R$, let $v\in L^\infty(\R,d\mu)$, and let $k>0$. The following statements are equivalent.
\begin{enumerate}[label=(\roman*)]
\item there is a Borel function $\hat v:\R\to\R$ with $\hat v=v$ $\mu$-almost everywhere, $|\hat v|\le\|v\|_{L^\infty(d\mu)}$ on $\R$, and $\hat v(x_2)-\hat v(x_1)\le k(x_2-x_1)$ for all $x_1<x_2$.
\item there is a $\mu$-null set $N$ such that $v(x_2)-v(x_1)\le k(x_2-x_1)$ for all $x_1<x_2$ in $\R\setminus N$.
\end{enumerate}
\end{lemma}

\begin{proof}
If (i) holds, (ii) holds with $N:=\{\hat v\ne v\}$.
Conversely, assume (ii), put $M:=\|v\|_{L^\infty(d\mu)}$, and enlarge $N$ by a $\mu$-null set so that $|v|\le M$ on $E:=\R\setminus N$. If $\mu=0$, take $\hat v\equiv0$.
Otherwise $E\neq\emptyset$, and $g(x):=v(x)-kx$ is nonincreasing on $E$. Define, for $x\in\R$,
$$
\hat g(x):=
\begin{cases}
\sup\{g(y):y\in E,\ y\ge x\}, & \text{if }E\cap[x,\infty)\ne\emptyset,\\
\inf\{g(y):y\in E\}, & \text{otherwise}.
\end{cases}
$$
In the first case $\hat g(x)\le M-kx$, and in the second case $E$ is bounded above by $x$ and $\inf\limits_Eg\ge-M-kx$. Thus $\hat g$ is finite. It is nonincreasing on $\R$ and coincides with $g$ on $E$, because $g$ is nonincreasing on $E$. Hence $\tilde v(x):=\hat g(x)+kx$ satisfies $\tilde v=v$ on $E$ and $\tilde v(x_2)-\tilde v(x_1)\le k(x_2-x_1)$ for $x_1<x_2$.
Finally, let $\hat v:=\max\{-M,\min\{M,\tilde v\}\}$. The truncation does not change $\tilde v$ on $E$, and, being nondecreasing and $1$-Lipschitz, it satisfies $\hat v(x_2)-\hat v(x_1)\le\max\{\tilde v(x_2)-\tilde v(x_1),0\}\le k(x_2-x_1)$. The function $\hat v$ is Borel, since $\hat g$ is monotone.
\end{proof}

Admissible representatives correspond to each other under $\mathcal T$, as the following proposition shows.

\begin{proposition}\label{representative-correspondence}
Let $(\rho,u)$ be a weak solution of \eqref{NA} on $[0,T)$, let $(\tilde\rho,w)=\mathcal T(\rho,u)$ be as in Theorem~\ref{weak-equivalence}, and let $t\in(0,T)$, $\tau=A(t)$. A bounded Borel function $\hat u(\cdot,t)$ is an admissible representative of $u(\cdot,t)$ if and only if $\hat w(\cdot,\tau):=e^{S(t)}\hat u(\cdot,t)$ satisfies $\hat w(\cdot,\tau)=w(\cdot,\tau)$ $\tilde\rho_\tau$-almost everywhere and
\begin{equation}\label{homogeneous-oleinik}
\frac{\hat w(x_2,\tau)-\hat w(x_1,\tau)}{x_2-x_1}\le\frac1\tau
\qquad\text{for all }x_1<x_2.
\end{equation}
Consequently, $u$ satisfies the intrinsic Oleinik condition if and only if $w$ satisfies the Oleinik condition of Definition~\ref{admissible-representative} for \eqref{H}, namely \eqref{homogeneous-oleinik} for almost every $\tau\in(0,A(T))$.
\end{proposition}

\begin{proof}
Since $\tilde\rho_\tau=\rho_t$ and $w(\cdot,\tau)=e^{S(t)}u(\cdot,t)$ $\rho_t$-almost everywhere, the identities $\hat u(\cdot,t)=u(\cdot,t)$ and $\hat w(\cdot,\tau)=w(\cdot,\tau)$ hold almost everywhere with respect to the same measure simultaneously. Moreover, $e^{S(t)}$ is a positive constant at fixed $t$, so
$$
\hat w(x_2,\tau)-\hat w(x_1,\tau)=e^{S(t)}\big[\hat u(x_2,t)-\hat u(x_1,t)\big],
$$
and \eqref{intrinsic-oleinik} holds if and only if $\hat w(x_2,\tau)-\hat w(x_1,\tau)\le(x_2-x_1)/A(t)=(x_2-x_1)/\tau$. For the last statement, the set of times at which an admissible representative fails to exist is null in $(0,T)$ if and only if its image under $A$ is null in $(0,A(T))$, by Lemma~\ref{properties-A}(iii).
\end{proof}

The second admissibility condition concerns the kinetic energy at the initial time.
Since the transformation multiplies the velocity by $e^{S(t)}$ and $e^{S(t)}\to1$ as $t\to0^+$, the condition is insensitive to this factor.

\begin{lemma}[Energy condition]\label{energy}
Let $(\rho,u)$ be a weak solution of \eqref{NA} on $[0,T)$ and $(\tilde\rho,w)=\mathcal T(\rho,u)$. Then, $w(\cdot,\tau)^2\,\tilde\rho_\tau \stackrel{*}{\rightharpoonup} u_0^2\,\rho_0$ in $\M_{\mathrm{loc}}(\R)$ as $\tau\to0^+$ if and only if $u(\cdot,t)^2\,\rho_t \stackrel{*}{\rightharpoonup} u_0^2\,\rho_0$ in $\M_{\mathrm{loc}}(\R)$ as $t\to0^+$.
\end{lemma}

\begin{proof}
For every $t\in[0,T)$, $\tilde\rho_{A(t)}=\rho_t$ and $w(\cdot,A(t))^2=e^{2S(t)}u(\cdot,t)^2$ $\rho_t$-almost everywhere, so
$$
w(\cdot,A(t))^2\,\tilde\rho_{A(t)}=e^{2S(t)}u(\cdot,t)^2\,\rho_t
$$
as measures. Let $\zeta\in C_c(\R)$ with $\supp\zeta\subset K$.
By Lemma~\ref{uniform-bound} and \eqref{local-uniform-u}, applied on the compact interval $I=[0,t_0]$ for some $t_0\in(0,T)$, there is a constant $C$ such that $\int_K u(\cdot,t)^2\,d\rho_t\le C$ for $t\in I$.
Hence
$$
\Big|\int_\R\zeta\,w(\cdot,A(t))^2\,d\tilde\rho_{A(t)}-\int_\R\zeta\,u(\cdot,t)^2\,d\rho_t\Big|
\le\big(e^{2S(t)}-1\big)\|\zeta\|_\infty\,C\longrightarrow0
\qquad\text{as }t\to0^+.
$$
Since $A$ is a homeomorphism of $[0,T)$ onto $[0,A(T))$ with $A(0)=0$, the limit as $\tau\to0^+$ of the first integral, with $\tau=A(t)$, exists if and only if the limit as $t\to0^+$ of the second exists, and the two limits coincide.
\end{proof}

With both admissibility conditions formulated intrinsically, we can define entropy solutions of \eqref{NA}.

\begin{definition}[Entropy solution]\label{entropy-solution}
A weak solution $(\rho,u)$ of \eqref{NA} on $[0,T)$, in the sense of Definition~\ref{weak-solution}, is called an {\em entropy solution} if it satisfies the following two conditions, which together form the Oleinik--energy entropy condition.
\begin{enumerate}[label=(\roman*)]
\item The velocity $u$ satisfies the intrinsic Oleinik condition of Definition~\ref{admissible-representative}.
\item The energy measure has the prescribed weak-$*$ initial trace, that is,
$$
u(\cdot,t)^2\,\rho_t
\stackrel{*}{\rightharpoonup}
u_0^2\,\rho_0
\quad\text{in }\M_{\mathrm{loc}}(\R)
\qquad\text{as }t\to0^+.
$$
\end{enumerate}
Entropy solutions of \eqref{H} are defined by the same two conditions, with $(\rho,w)$ in place of $(\rho,u)$ and with the bound $1/\tau$ in place of $e^{-S(t)}/A(t)$ in \eqref{intrinsic-oleinik}, that is, with \eqref{homogeneous-oleinik}.
\end{definition}

For \eqref{H}, this is the notion of \cite[Definition~2]{HuangWang2001}, up to the choice of a representative of the velocity, as shown in Proposition~\ref{HW-equivalence}. We use $u^2\rho$ as the energy measure, without the conventional factor $1/2$. Physical kinetic energy is one half of this measure. The energy condition in Definition~\ref{entropy-solution} cannot be dropped.
Weak solutions that satisfy the one-sided bound but split a concentrated mass were exhibited for the homogeneous system by Bouchut and James \cite[Example~2.1]{BouchutJames1999} and by Huang and Wang \cite[(4.37)--(4.39)]{HuangWang2001}.
Moreover, Bouchut \cite[discussion following Proposition~3.3]{Bouchut1994} observed that, among the solutions issued from a collision of two particles, the sticky one has the smallest kinetic energy.
The simplest way to see the phenomenon in the present setting is to follow a single Dirac mass with zero initial velocity, which a weak solution satisfying the one-sided bound may keep at the origin, split into two halves, or spread into a fan.

\begin{example}[Splitting of a Dirac mass with zero initial velocity]\label{splitting}
Let Assumption~\ref{coefficients} hold and consider the data $\rho_0=\delta_0$ and $u_0=0$. For $t>0$, let
$$
\bar u(x,t)=
\begin{cases}
\dfrac{x\,e^{-S(t)}}{A(t)}, & |x|\le A(t),\\[2mm]
e^{-S(t)}\sgn x, & |x|>A(t).
\end{cases}
$$
Each of the following pairs, completed by $(\rho_0,u_0)$ at $t=0$, is a weak solution of \eqref{NA} on $[0,\infty)$.
$$
\begin{aligned}
&\text{(i)}\quad && \rho_t=\delta_0, && u(x,t)=0,\\[1mm]
&\text{(ii)}\quad && \rho_t=\tfrac12\,\delta_{-A(t)}+\tfrac12\,\delta_{A(t)}, && u=\bar u,\\[1mm]
&\text{(iii)}\quad && \rho_t=\frac{1}{2A(t)}\,\mathbf 1_{(-A(t),A(t))}\,dx, && u=\bar u.
\end{aligned}
$$
In (ii) the Dirac mass splits into two halves that move apart along $x=\pm A(t)$ with velocities $\pm e^{-S(t)}$, and in (iii) it spreads uniformly over the interval $(-A(t),A(t))$.
In each case, $\rho_t$ is the image of a probability measure $\bar\gamma$ on $[-1,1]$ with zero mean under the map $\xi\mapsto A(t)\xi$, namely $\bar\gamma=\delta_0$, $\bar\gamma=\tfrac12(\delta_{-1}+\delta_1)$, and $\bar\gamma=\tfrac12\mathbf 1_{[-1,1]}\,d\xi$, and the particle at $A(t)\xi$ has velocity $\bar u(A(t)\xi,t)=e^{-S(t)}\xi$.
Since $A'=\alpha e^{-S}$ and $S'=\sigma$, for test functions $\Xi,\Psi\in C_c^\infty(\R\times(0,\infty))$ we have
$$
\begin{gathered}
\frac{d}{dt}\Xi\big(A(t)\xi,t\big)=\partial_t\Xi+\alpha e^{-S}\xi\,\partial_x\Xi,\\
\frac{d}{dt}\Big[e^{-S}\Psi\big(A(t)\xi,t\big)\Big]=e^{-S}\partial_t\Psi+\alpha e^{-2S}\xi\,\partial_x\Psi-\sigma e^{-S}\Psi,
\end{gathered}
$$
with the derivatives of $\Xi$ and $\Psi$ evaluated at $(A(t)\xi,t)$. Integrating the first identity in $t$ and with respect to $\bar\gamma$ gives \eqref{weak1-rho}, and integrating the second one multiplied by $\xi$ gives \eqref{weak2}. Conditions (a) and (b) follow from the continuity of $A$ and $S$, and the initial conditions hold because $A(t)\to0$ and $\bar\gamma$ has zero mean. In all three cases the total mass equals $1$ and the total momentum vanishes for every $t$, and $u(\cdot,t)$ is continuous, nondecreasing, and piecewise affine with slopes $0$ and $e^{-S(t)}/A(t)$, so that it satisfies the intrinsic Oleinik condition.
The kinetic energies, however, differ. For every $\varphi\in C_c(\R)$, as $t\to0^+$,
$$
\int_\R\varphi\,u^2\,d\rho_t\longrightarrow
\begin{cases}
0 & \text{in case (i)},\\[1mm]
\varphi(0) & \text{in case (ii)},\\[1mm]
\tfrac13\,\varphi(0) & \text{in case (iii)},
\end{cases}
$$
since $\int_\R u^2\,d\rho_t$ equals $0$, $e^{-2S(t)}$, and $\tfrac13e^{-2S(t)}$, respectively. Hence only (i) satisfies the energy condition, and it is an entropy solution.
For $\alpha\equiv1$ and $\sigma\equiv0$, so that $A(t)=t$ and $S\equiv0$, these are, up to the normalization of the mass, the three solutions of Bouchut and James \cite[Example~2.1]{BouchutJames1999}.
\end{example}

The solutions (ii) and (iii) are obtained by letting the velocity of the Dirac mass split into a distribution of velocities with the same mean.
The next proposition shows that every such distribution produces a weak solution satisfying the intrinsic Oleinik condition, also in the presence of surrounding mass, and that the energy condition rejects all of them except the one without splitting, because the defect in the initial trace of $u^2\rho$ equals $W_0$ times the variance of $\gamma/W_0$.

\begin{proposition}[Nonuniqueness without the energy condition]
\label{nonuniqueness}
Let Assumption~\ref{coefficients} hold, let $\rho_\pm\ge0$, $W_0>0$, and $u_-<c<u_+$, and consider
$$
\rho_0=\rho_-\mathbf 1_{(-\infty,0)}\,dx+W_0\,\delta_0+\rho_+\mathbf 1_{(0,\infty)}\,dx,
\qquad
u_0=u_-\mathbf 1_{(-\infty,0)}+c\,\mathbf 1_{\{0\}}+u_+\mathbf 1_{(0,\infty)}.
$$
Let $\gamma$ be a nonnegative Borel measure on $[u_-,u_+]$ such that $\gamma([u_-,u_+])=W_0$ and $\int_{[u_-,u_+]}\xi\,d\gamma(\xi)=W_0c$.
For $t>0$, let $\gamma_t$ be the image of $\gamma$ under the map $[u_-,u_+]\ni\xi\mapsto A(t)\xi\in\R$, that is, the measure on $\R$ given by $\gamma_t(E):=\gamma(\{\xi\in[u_-,u_+]:A(t)\xi\in E\})$ for Borel sets $E\subset\R$, and define
$$
\rho^\gamma_t=\rho_-\mathbf 1_{(-\infty,u_-A(t))}\,dx+\gamma_t+\rho_+\mathbf 1_{(u_+A(t),\infty)}\,dx
$$
and
$$
u^\gamma(x,t)=
\begin{cases}
u_-e^{-S(t)}, & x<u_-A(t),\\[1mm]
\dfrac{x\,e^{-S(t)}}{A(t)}, & u_-A(t)\le x\le u_+A(t),\\[2mm]
u_+e^{-S(t)}, & x>u_+A(t),
\end{cases}
$$
together with $(\rho^\gamma_0,u^\gamma(\cdot,0))=(\rho_0,u_0)$. Then $(\rho^\gamma,u^\gamma)$ is a weak solution of \eqref{NA} on $[0,\infty)$, and $u^\gamma$ satisfies the intrinsic Oleinik condition.
Moreover, $(\rho^\gamma,u^\gamma)$ satisfies the energy condition of Definition~\ref{entropy-solution} if and only if $\gamma=W_0\,\delta_c$.
In particular, \eqref{NA} with these data has infinitely many weak solutions satisfying the intrinsic Oleinik condition, and exactly one member of this family is an entropy solution.
\end{proposition}

\begin{proof}
All the solutions in the statement describe free particles that never collide, and it is convenient to record their initial positions and velocities in one measure.
Let $\pi$ be the nonnegative Radon measure on $\R^2$, with coordinates $(x_0,\xi)$, given by
$$
\pi=\rho_-\mathbf 1_{(-\infty,0)}(x_0)\,dx_0\otimes\delta_{u_-}
+\delta_0\otimes \gamma
+\rho_+\mathbf 1_{(0,\infty)}(x_0)\,dx_0\otimes\delta_{u_+},
$$
where $\otimes$ denotes the product of a measure in the variable $x_0$ and a measure in the variable $\xi$.
The first and third terms describe the particles initially located at $x_0<0$ and at $x_0>0$, which have velocities $u_-$ and $u_+$, respectively, and the second term describes the Dirac mass at the origin, whose mass $W_0$ is distributed among the velocities $\xi\in[u_-,u_+]$ according to $\gamma$.
A particle issued from $x_0$ with velocity $\xi$ is located at $x_0+A(t)\xi$ and has velocity $e^{-S(t)}\xi$ at time $t$. Particles issued from $x_0<0$ remain in the region $x<u_-A(t)$, those issued from $x_0>0$ remain in $x>u_+A(t)$, and those issued from the origin with velocity $\xi\in[u_-,u_+]$ are located at $x=A(t)\xi$, where $u^\gamma(x,t)=x\,e^{-S(t)}/A(t)=e^{-S(t)}\xi$.
Consequently, for every $t>0$ and every bounded Borel function $\varphi$ with compact support,
$$
\int_\R\varphi\,u^\gamma(\cdot,t)^k\,d\rho^\gamma_t
=\int_{\R^2}\varphi\big(x_0+A(t)\xi\big)\,\big(e^{-S(t)}\xi\big)^k\,d\pi(x_0,\xi)
$$
for $k=0,1,2$. For $k=0,1$ this identity also holds at $t=0$, because $\int_{[u_-,u_+]}\xi\,d\gamma(\xi)=W_0c$.
If $\supp\varphi\subset[-R,R]$ and $t\le t_1$, the integrand vanishes unless $|x_0|\le R+A(t_1)\max\{|u_-|,|u_+|\}$, and $\pi$ is finite on the corresponding strip.
Dominated convergence therefore shows that the integrals with $k=0,1$ are continuous on $[0,\infty)$ and that the integral with $k=2$ is continuous on $(0,\infty)$. This gives conditions (a) and (b) of Definition~\ref{weak-solution}, since $u^\gamma$ is Borel and $|u^\gamma|\le e^{-S(t)}\max\{|u_-|,|u_+|\}$.

For $\Xi\in C_c^\infty(\R\times(0,\infty))$, the function $t\mapsto\Xi(x_0+A(t)\xi,t)$ is of class $C^1$ with compact support in $(0,\infty)$, and, since $A'=\alpha e^{-S}$,
$$
\frac{d}{dt}\,\Xi\big(x_0+A(t)\xi,t\big)
=\partial_t\Xi\big(x_0+A(t)\xi,t\big)+\alpha(t)e^{-S(t)}\xi\,\partial_x\Xi\big(x_0+A(t)\xi,t\big).
$$
Integrating in $t$ and then with respect to $\pi$, which is allowed by the bound above and Fubini's theorem, gives \eqref{weak1-rho}. Similarly, for $\Psi\in C_c^\infty(\R\times(0,\infty))$, the identities $A'=\alpha e^{-S}$ and $S'=\sigma$ give
$$
\frac{d}{dt}\Big[e^{-S(t)}\,\Psi\big(x_0+A(t)\xi,t\big)\Big]
=e^{-S}\,\partial_t\Psi+\alpha e^{-2S}\xi\,\partial_x\Psi-\sigma e^{-S}\,\Psi,
$$
with $\Psi$ and its derivatives evaluated at $(x_0+A(t)\xi,t)$, and multiplying by $\xi$ and integrating yields \eqref{weak2}.

As $t\to0^+$, $A(t)\to0$ and $e^{-S(t)}\to1$, so dominated convergence gives $\int_\R\varphi\,d\rho^\gamma_t\to\int_{\R^2}\varphi(x_0)\,d\pi(x_0,\xi)=\int_\R\varphi\,d\rho_0$ and $\int_\R\varphi\,u^\gamma\,d\rho^\gamma_t\to\int_{\R^2}\varphi(x_0)\,\xi\,d\pi(x_0,\xi) =\int_\R\varphi\,u_0\,d\rho_0$, where the last identity uses $\int_{[u_-,u_+]}\xi\,d\gamma(\xi)=W_0c$.
This is condition (d), and $(\rho^\gamma,u^\gamma)$ is a weak solution.

The function $u^\gamma(\cdot,t)$ is continuous, nondecreasing, and piecewise affine with slopes $0$ and $e^{-S(t)}/A(t)$. Hence it is an admissible representative of itself in the sense of Definition~\ref{admissible-representative}, for every $t>0$.

Finally, by the same argument,
$$
\int_\R\varphi\,u^\gamma(\cdot,t)^2\,d\rho^\gamma_t
\longrightarrow\int_{\R^2}\varphi(x_0)\,\xi^2\,d\pi(x_0,\xi)
=\int_\R\varphi\,u_0^2\,d\rho_0+\varphi(0)\Big(\int_{[u_-,u_+]}\xi^2\,d\gamma(\xi)-W_0c^2\Big),
$$
and, since $\int_{[u_-,u_+]}\xi\,d\gamma(\xi)=W_0c$ and $\gamma([u_-,u_+])=W_0$,
$$
\int_{[u_-,u_+]}\xi^2\,d\gamma(\xi)-W_0c^2=\int_{[u_-,u_+]}(\xi-c)^2\,d\gamma(\xi).
$$
Writing $\bar\gamma:=\gamma/W_0$, this defect is $W_0\operatorname{Var}_{\bar\gamma}(\xi)$. With the physical energy convention it is $\tfrac12W_0\operatorname{Var}_{\bar\gamma}(\xi)$.
The right-hand side vanishes if and only if $\gamma$ is concentrated at $\xi=c$, that is, $\gamma=W_0\,\delta_c$. Thus the energy condition holds exactly for this choice of $\gamma$. Since $u_-<c<u_+$, the measures $\gamma=W_0\big((1-\vartheta)\delta_c+\tfrac\vartheta2(\delta_{c-r}+\delta_{c+r})\big)$, with $0<r\le\min\{c-u_-,u_+-c\}$ and $\vartheta\in[0,1]$, provide infinitely many admissible choices.
\end{proof}

Example~\ref{splitting} corresponds to $\rho_\pm=0$, $W_0=1$, $c=0$, $u_\pm=\pm1$, and $\gamma=\delta_0$, $\gamma=\tfrac12(\delta_{-1}+\delta_1)$, and $\gamma=\tfrac12\mathbf 1_{[-1,1]}\,d\xi$. For $\alpha\equiv1$, $\sigma\equiv0$, $\rho_\pm=1$, $W_0=1$, $u_-=0$, $c=1$, and $u_+=2$, the choices $\gamma=\delta_1$ and $\gamma=\tfrac12\mathbf 1_{[0,2]}\,d\xi$ give the two solutions \cite[(4.38) and (4.39)]{HuangWang2001}, and the energy of the second exceeds the initial energy of the Dirac mass by $\int_0^2(\xi-1)^2\,\tfrac{d\xi}{2}=\tfrac13$.
For $\alpha\equiv1$ and $\sigma\equiv1$, the same choices give the counterexample of Jin \cite[(3.18)--(3.20)]{Jin2016}.

\subsection{Equivalence of the Cauchy problems}
\label{cauchy-equivalence-section}

The preceding results combine into the correspondence of entropy solutions.

\begin{theorem}[Equivalence of entropy solutions]\label{entropy-correspondence}
Let Assumption~\ref{coefficients} hold, let $T\in(0,\infty]$, let $(\rho,u)$ be a weak solution of \eqref{NA} on $[0,T)$, and let $(\tilde\rho,w)=\mathcal T_T(\rho,u)$. Then $(\rho,u)$ is an entropy solution of \eqref{NA} on $[0,T)$ if and only if $(\tilde\rho,w)$ is an entropy solution of \eqref{H} on $[0,A(T))$. Consequently, $\mathcal T_T$ maps the set of entropy solutions of \eqref{NA} on $[0,T)$ bijectively onto the set of entropy solutions of \eqref{H} on $[0,A(T))$ with the same initial data.
\end{theorem}

\begin{proof}
By Theorem~\ref{weak-equivalence}, $(\tilde\rho,w)$ is a weak solution of \eqref{H} on $[0,A(T))$ with the same initial data, and $\mathcal T_T$ is a bijection between the sets of weak solutions. By Proposition~\ref{representative-correspondence}, $u$ satisfies the intrinsic Oleinik condition if and only if $w$ satisfies the homogeneous one, and by Lemma~\ref{energy}, the energy conditions are equivalent. Hence $\mathcal T_T$ maps entropy solutions onto entropy solutions.
\end{proof}

Theorem~\ref{entropy-correspondence} can be stated at the level of solution spaces. The velocity of an entropy solution is determined only $\rho_t$-almost everywhere, and the uniqueness theorem of the homogeneous theory asserts equality only in this sense, so the natural solution spaces are quotient spaces.

\begin{definition}[Equivalent entropy solutions]\label{equiv-class}
Two entropy solutions $(\rho^1,u^1)$ and $(\rho^2,u^2)$ of \eqref{NA} on $[0,T)$ with the same initial data are equivalent, written $(\rho^1,u^1)\sim(\rho^2,u^2)$, if, for every $t\in[0,T)$, $\rho^1_t=\rho^2_t$ and $u^1(\cdot,t)\rho^1_t=u^2(\cdot,t)\rho^2_t$ as measures, or, equivalently, $\rho^1_t=\rho^2_t$ and $u^1(\cdot,t)=u^2(\cdot,t)$ $\rho^1_t$-almost everywhere, for every $t$. The same relation is used for \eqref{H}.
\end{definition}

Admissibility does not depend on the chosen member of an equivalence class.

\begin{lemma}[Independence of the representative]\label{entropy-well-defined}
Let $(\rho^1,u^1)$ be an entropy solution of \eqref{NA} on $[0,T)$, and let $(\rho^2,u^2)\in\mathcal P_T$ satisfy $\rho^2_t=\rho^1_t$ and $u^2(\cdot,t)=u^1(\cdot,t)$ $\rho^1_t$-almost everywhere for every $t$. Then $(\rho^2,u^2)$ is an entropy solution of \eqref{NA} on $[0,T)$.
\end{lemma}

\begin{proof}
All conditions of Definitions~\ref{weak-solution} and~\ref{entropy-solution} involve $u$ only through the measures $u(\cdot,t)\rho_t$ and $u(\cdot,t)^2\rho_t$, and through the $L^\infty(d\rho_t)$ norm and the $\rho_t$-almost everywhere class of $u(\cdot,t)$. These are the same for $u^1$ and $u^2$. In particular, an admissible representative of $u^1(\cdot,t)$ is also an admissible representative of $u^2(\cdot,t)$.
\end{proof}

For $T\in(0,\infty]$, we denote by $\mathcal E_{\mathrm{NA}}(T;\rho_0,u_0)$ the set of entropy solutions of \eqref{NA} on $[0,T)$ with initial data $(\rho_0,u_0)$, modulo the relation of Definition~\ref{equiv-class}, and by $\mathcal E_{\mathrm H}(A(T);\rho_0,u_0)$ the analogous quotient for \eqref{H} on $[0,A(T))$. The equivalence class of an entropy solution $(\rho,u)$ is denoted by $[(\rho,u)]$.

\begin{theorem}[Equivalence of the Cauchy problems]
\label{cauchy-equivalence}
Let Assumption~\ref{coefficients} hold, let $\rho_0\in\M_{\mathrm{loc}}(\R)$ be nonnegative, let $u_0\in L^\infty(\R,d\rho_0)$, and let $T\in(0,\infty]$.
The map
$$
[\mathcal T_T]:
\mathcal E_{\mathrm{NA}}(T;\rho_0,u_0)
\longrightarrow
\mathcal E_{\mathrm H}(A(T);\rho_0,u_0),
\qquad
[\mathcal T_T]\big([(\rho,u)]\big):=\big[\mathcal T_T(\rho,u)\big],
$$
is well defined and bijective, and its inverse is $[(\tilde\rho,w)]\mapsto[\mathcal T_T^{-1}(\tilde\rho,w)]$.
\end{theorem}

\begin{proof}
By Theorem~\ref{entropy-correspondence}, $\mathcal T_T$ is a bijection between the sets of entropy solutions with the given data, before passing to the quotient. It remains to check that $\mathcal T_T$ and $\mathcal T_T^{-1}$ respect the equivalence relations. If $(\rho^1,u^1)\sim(\rho^2,u^2)$, then $\tilde\rho^1_\tau=\rho^1_{A^{-1}(\tau)}=\rho^2_{A^{-1}(\tau)}=\tilde\rho^2_\tau$, and, multiplying the equality $u^1(\cdot,t)\rho^1_t=u^2(\cdot,t)\rho^2_t$ at $t=A^{-1}(\tau)$ by the common scalar $e^{S(t)}$, we obtain $w^1(\cdot,\tau)\tilde\rho^1_\tau=w^2(\cdot,\tau)\tilde\rho^2_\tau$ for every $\tau$. The same argument applies to $\mathcal T_T^{-1}$. Hence $[\mathcal T_T]$ and the map induced by $\mathcal T_T^{-1}$ are well defined and mutually inverse.
\end{proof}

By Theorem~\ref{cauchy-equivalence}, existence of an entropy solution of \eqref{NA} on $[0,T)$ is equivalent to existence of an entropy solution of \eqref{H} on $[0,A(T))$ with the same data, and the same holds for uniqueness modulo Definition~\ref{equiv-class}. Both questions are therefore reduced to the homogeneous theory, provided that the solution class of Definition~\ref{entropy-solution} is the one in which that theory is formulated. This is the subject of the next section.

\section{Global existence and uniqueness}
\label{global-existence-uniqueness}

In this section we prove existence and uniqueness of entropy solutions of \eqref{NA}. By Theorem~\ref{cauchy-equivalence}, both questions reduce to the homogeneous system on $[0,A_\infty)$. Existence is obtained from the theorem of Huang and Wang \cite{HuangWang2001}, after we recall their formulation in the present notation. Uniqueness is proved directly on finite time strips. Finally, we show that the intrinsic solution class coincides with the class of the homogeneous theory.

\subsection{Existence}
\label{bridge}

We recall the formulation of \cite[Definitions~1 and~2]{HuangWang2001} in the notation of the present paper. There a pair $(m,w)$, with $m(\cdot,\tau)\in BV_{\mathrm{loc}}(\R)$ and $w$ bounded and measurable with respect to $\partial_xm$, is a weak solution of \eqref{H} if the measures $\partial_xm$ and $w\,\partial_xm$ are weakly continuous in time and
\begin{equation}\label{HW-weak}
\begin{gathered}
\iint_{\R\times(0,\infty)} \partial_\tau\varphi\, m\,dx\,d\tau-
\iint_{\R\times(0,\infty)} \varphi\, w\,dm\,d\tau=0,\\
\iint_{\R\times(0,\infty)}\big(\partial_\tau\psi\,w+\partial_x\psi\,w^2\big)\,dm\,d\tau=0
\end{gathered}
\end{equation}
for all $\varphi,\psi\in C_c^\infty(\R\times(0,\infty))$. Here $dm$ stands for the measure $\partial_xm$, and the integrals with respect to $dm\,d\tau$ are Lebesgue--Stieltjes integrals on $\R\times(0,\infty)$.
The initial data are attained in the sense that $\partial_xm$ and $w\,\partial_xm$ converge weakly to $\rho_0$ and $u_0\rho_0$ as $\tau\to0^+$. Such a weak solution is an entropy solution if, for all $x_1<x_2$ and almost every $\tau>0$,
\begin{equation}\label{HW-Oleinik}
\frac{w(x_2,\tau)-w(x_1,\tau)}{x_2-x_1}\le \frac1\tau,
\end{equation}
and if $w(\cdot,\tau)^2\partial_xm(\cdot,\tau)$ converges weakly to $u_0^2\rho_0$ as $\tau\to0^+$. Since $w(\cdot,\tau)$ is determined only $\partial_xm$-almost everywhere, \eqref{HW-Oleinik} can be read either for all $x_1<x_2$, for a suitable representative, or for $x_1<x_2$ outside an $\partial_xm$-null set. By Lemma~\ref{oleinik-forms}, the two readings are equivalent, and both are equivalent to \eqref{homogeneous-oleinik}. These definitions are stated on the whole half-plane $\tau>0$, and on a strip $\R\times[0,T)$ we use them with test functions in $C_c^\infty(\R\times(0,T))$ and with all conditions restricted to $\tau<T$. Throughout this comparison, boundedness on $[0,T)$ means a uniform bound on each $\R\times[0,T_0]$ with $T_0<T$. Thus the comparison is with the strip-local version of their class.

We now make precise the integrals with respect to $dm\,d\tau$. Let $\rho_\tau:=\partial_xm(\cdot,\tau)$, and let $\mu$ be the measure associated with $(\rho_\tau)$ by Lemma~\ref{kernel}(iii). The Lebesgue--Stieltjes measure $dm\,d\tau$ on $\R\times(0,T)$ is the kernel measure associated with $\tau\mapsto \partial_xm(\cdot,\tau)$, see Ding and Wang \cite[remark following Theorem~2]{DingWang1996}, that is, the Borel measure that assigns to each rectangle $[a,b)\times J$, with $a<b$ and $J$ an interval whose closure is contained in $(0,T)$, the value $\int_J\rho_\tau([a,b))\,d\tau$. This value equals $\mu([a,b)\times J)$. These rectangles, together with the empty set, form a $\pi$-system that generates the Borel $\sigma$-algebra of $\R\times(0,T)$, and $\R\times(0,T)$ is a countable union of such rectangles, each of finite $\mu$-measure by Lemma~\ref{uniform-bound}. By the uniqueness theorem for measures that coincide on a generating $\pi$-system and are $\sigma$-finite on it, $dm\,d\tau$ coincides with $\mu$ on the Borel subsets of $\R\times(0,T)$. Consequently, for a Borel velocity bounded on every strip $\R\times(0,\tau_1]$ with $\tau_1<T$, the integrals in \eqref{HW-weak} coincide with the iterated integrals of Definition~\ref{weak-solution}, by Lemma~\ref{kernel}(iii) applied to the positive and negative parts of the integrands. Such a velocity is jointly Borel measurable and is therefore measurable for the completion of the kernel measure, as required in the definition of weak solution recalled above.

We now recall the existence theorem of the homogeneous theory in the present notation.

\begin{theorem}[Huang and Wang {\cite[Theorem~1]{HuangWang2001}}]
\label{HW}
Let $\rho_0\in\M_{\mathrm{loc}}(\R)$ be nonnegative and let $u_0\in L^\infty(\R,d\rho_0)$. Then system \eqref{H} admits a global entropy solution $(\tilde m,w)$ in the sense recalled above. The solution constructed in the proof of this theorem, with $\tilde m$ and $w$ given by \cite[(2.35) and (2.29)]{HuangWang2001}, satisfies $|w|\le\|u_0\|_{L^\infty(d\rho_0)}$.
\end{theorem}

The bound is part of the definition of the velocity $w$, which takes values in $[-M,M]$ with $M=\|u_0\|_{L^\infty(d\rho_0)}$.

The entropy solution of \eqref{NA} is obtained by transporting the solution of Theorem~\ref{HW} through $\mathcal T^{-1}$.

\begin{theorem}[Existence]
\label{existence}
Let Assumption~\ref{coefficients} hold, let $\rho_0\in\M_{\mathrm{loc}}(\R)$ be nonnegative, and let $u_0\in L^\infty(\R,d\rho_0)$. Then \eqref{NA} admits an entropy solution on $[0,\infty)$, whose velocity satisfies $|u(\cdot,t)|\le e^{-S(t)}\|u_0\|_{L^\infty(d\rho_0)}$ $\rho_t$-almost everywhere.
\end{theorem}

\begin{proof}
Let $(\tilde m,w)$ be the global entropy solution of Theorem~\ref{HW}, with $\tilde m$ and $w$ as in Theorem~\ref{HW}. Put $\tilde\rho_\tau:=\partial_x\tilde m(\cdot,\tau)$ for $\tau>0$, $\tilde\rho_0:=\rho_0$, and $w(\cdot,0):=u_0$, where $u_0$ denotes a representative with $|u_0|\le M:=\|u_0\|_{L^\infty(d\rho_0)}$ everywhere.
We show that $(\tilde\rho,w)$ is an entropy solution of \eqref{H} on $[0,\infty)$ in the sense of Definitions~\ref{weak-solution} and~\ref{entropy-solution}.
Conditions (a), (b), and (d) are part of the definition of weak solution recalled above, which requires $w$ to be bounded and measurable with respect to $\partial_x\tilde m$ and the measures $\partial_x\tilde m$ and $w\,\partial_x\tilde m$ to be weakly continuous in time with initial values $\rho_0$ and $u_0\rho_0$, together with the bound $|w|\le M$ of Theorem~\ref{HW}.
For the identities of Definition~\ref{weak-solution}, we use the form in which \eqref{HW-weak} is verified in the proof of existence \cite[Lemma~2.8, (2.53), and (2.57)--(2.59)]{HuangWang2001}. There, the integrals with respect to $dm\,d\tau$ are computed at fixed $\tau$ through the spatial primitives of $w(\cdot,\tau)\tilde\rho_\tau$ and $w(\cdot,\tau)^2\tilde\rho_\tau$, followed by an integration in $\tau$. Hence \eqref{HW-weak} holds in the form
$$
\int_0^\infty\!\!\int_\R\partial_\tau\varphi\,\tilde m\,dx\,d\tau-\int_0^\infty\!\!\int_\R\varphi\,w\,d\tilde\rho_\tau\,d\tau=0,
\qquad
\int_0^\infty\!\!\int_\R\big(\partial_\tau\psi w+\partial_x\psi w^2\big)\,d\tilde\rho_\tau\,d\tau=0,
$$
for all $\varphi,\psi\in C_c^\infty(\R\times(0,\infty))$, where the inner integrals are Borel functions of $\tau$ by Lemma~\ref{borel-velocity}.
The second identity is \eqref{weak2-H}. In the first identity, take $\varphi=\partial_x\Phi$ with $\Phi\in C_c^\infty(\R\times(0,\infty))$. Since $\tilde m(\cdot,\tau)$ is a nondecreasing function with $\partial_x\tilde m(\cdot,\tau)=\tilde\rho_\tau$, an integration by parts at fixed $\tau$ gives $\int_\R\partial_x\partial_\tau\Phi\tilde m(\cdot,\tau)\,dx=-\int_\R\partial_\tau\Phi\,d\tilde\rho_\tau$, and the first identity becomes \eqref{weak1-H}.
Thus $(\tilde\rho,w)$ is a weak solution of \eqref{H} on $[0,\infty)$.
By the definition of entropy solution recalled above, $w(\cdot,\tau)$ satisfies \eqref{HW-Oleinik} for all $x_1<x_2$ and almost every $\tau>0$, so Lemma~\ref{oleinik-forms}, applied with the empty exceptional set, provides an admissible representative of $w(\cdot,\tau)$ for almost every $\tau$. The energy condition is also part of that definition.
Therefore $(\tilde\rho,w)$ is an entropy solution of \eqref{H} on $[0,\infty)$ in the sense of Definition~\ref{entropy-solution}, and so is its restriction to $[0,A_\infty)$.
By Theorem~\ref{entropy-correspondence} with $T=\infty$, the pair
$$
\rho_t:=\tilde\rho_{A(t)},
\qquad
u(x,t):=e^{-S(t)}w(x,A(t))
$$
is an entropy solution of \eqref{NA} on $[0,\infty)$, and the bound on $u$ follows from $|w|\le\|u_0\|_{L^\infty(d\rho_0)}$.
If $A_\infty<\infty$, only the restriction of the homogeneous solution to $[0,A_\infty)$ is used.
\end{proof}

\subsection{Uniqueness on finite time strips}
\label{strips}

Huang and Wang also prove that entropy solutions defined for all $\tau>0$, with a velocity bounded on the whole half-plane, are unique \cite[Theorem~2]{HuangWang2001}. We do not use this theorem, for two reasons. First, when $A_\infty<\infty$ we need uniqueness on a strip $[0,T)$ with $T$ finite, for velocities that are bounded only on each strip $\R\times[0,T_0]$ with $T_0<T$. Uniqueness on the half-line does not imply uniqueness on strips, because an entropy solution on a strip need not be the restriction of an entropy solution on the half-line. A continuation from a positive time $\tau_0$ satisfies the one-sided bound only with the constant $1/(\tau-\tau_0)$, which is larger than $1/\tau$. Second, in their proof the passage to the initial time uses that an exceptional set of initial points is Lebesgue null, and this set may have positive measure, as we explain in Section~\ref{immediate-concentration}. We therefore give a direct proof of uniqueness on strips.

Throughout this subsection and Section~\ref{lagrangian}, $(\rho,w)$ denotes a solution of \eqref{H}, without the tilde used for the transformed solutions. The main result of this subsection is that two entropy solutions of \eqref{H} on a strip $[0,T)$ with the same initial data have the same density and the same momentum at every time.

The proof follows the three steps of the argument of the homogeneous theory, namely a description of the solution at positive times, its representation through a generalized potential built from its state at a positive time $s$, and the passage $s\to0^+$. We carry out the first two steps through a Lagrangian description in mass coordinates, in which every concentration is an interval of mass coordinates and its velocity is read off directly, and the third step through a passage to the limit in the generalized potentials themselves. The Oleinik condition enters only through the forward uniqueness of characteristics, and the energy condition only through the vanishing of the velocity variance at the initial time.

The proof uses the forward uniqueness of characteristics, which requires the one-sided bound at every pair of points, including points of vacuum, for a single velocity function defined jointly in space and time. The Oleinik condition of Definition~\ref{entropy-solution} requires the bound only for some representative at almost every time. A Borel velocity is provided by Lemma~\ref{borel-velocity}, and the next lemma constructs a Borel velocity that satisfies the one-sided bound at every pair of points and still defines the given momentum at every time.

\begin{lemma}\label{joint-velocity}
Let $T\in(0,\infty]$, and let $(\rho,w)$ satisfy conditions (a) and (b) of Definition~\ref{weak-solution} for \eqref{H} on $[0,T)$. Assume that, for almost every $\tau\in(0,T)$, there is a bounded Borel function that coincides with $w(\cdot,\tau)$ $\rho_\tau$-almost everywhere and satisfies \eqref{homogeneous-oleinik}.
Then there is a Borel function $\check w:\R\times[0,T)\to\R$, bounded on $\R\times[0,\tau_1]$ for every $\tau_1\in(0,T)$, such that $\check w(\cdot,\tau)=w(\cdot,\tau)$ $\rho_\tau$-almost everywhere for every $\tau\in[0,T)$ and, for almost every $\tau\in(0,T)$,
$$
\check w(x_2,\tau)-\check w(x_1,\tau)\le\frac{x_2-x_1}{\tau}
\qquad\text{for all }x_1<x_2.
$$
If $\|w(\cdot,\tau)\|_{L^\infty(d\rho_\tau)}\le M$ for every $\tau\in[0,T)$, then $\check w$ can be chosen with $|\check w|\le M$ everywhere.
\end{lemma}

\begin{proof}
The idea is to replace, at each time, the velocity minus $x/\tau$ by its essential upper envelope on the right, which is monotone by construction and depends measurably on $(x,\tau)$.
Let $h$ be the Borel function of Lemma~\ref{borel-velocity} for $(\rho,w)$.
Choose $0=s_0<s_1<s_2<\dots$ with $s_l\to T$, put $M_l:=\sup\limits_{\tau\in[0,s_l]}\|w(\cdot,\tau)\|_{L^\infty(d\rho_\tau)}$, which is finite by \eqref{local-uniform-u}, and define $B(\tau):=M_l$ for $\tau\in(s_{l-1},s_l]$. If a bound $M$ as in the statement is given, put $B\equiv M$ instead. Then $B$ is Borel, bounded on every interval $(0,\tau_1]$ with $\tau_1<T$, and $|h(x,\tau)|\le B(\tau)$ for all $x\in\R$ and $\tau\in(0,T)$.
Put $g(y,\tau):=h(y,\tau)-y/\tau$ for $\tau\in(0,T)$, with $g(y,0):=0$, and let $\mathcal N$ be the set of $\tau\in(0,T)$ with $\rho_\tau\ne0$. Since $\rho_\tau\ne0$ if and only if $\rho_\tau((-n,n))>0$ for some $n\in\N$, $\mathcal N$ is Borel by Lemma~\ref{kernel}(i).
For $\tau\in\mathcal N$ and $x\in\R$, define
$$
H(x,\tau):=\mathop{\mathrm{ess\,sup}}_{y\in[x,\infty)}g(y,\tau),
\qquad
L(\tau):=\mathop{\mathrm{ess\,inf}}_{y\in\R}g(y,\tau),
$$
where the essential bounds are taken with respect to $\rho_\tau$, and the essential supremum over a $\rho_\tau$-null set is $-\infty$.

We first show that $H$ and $L$ are Borel. For rational $q$, $H(x,\tau)>q$ if and only if $\rho_\tau(\{y\ge x:g(y,\tau)>q\})>0$, and $L(\tau)<q$ if and only if $\rho_\tau(\{y:g(y,\tau)<q\})>0$.
By Lemma~\ref{kernel}(ii), applied with $Z=\R$ to the indicator function of the Borel set $\{(x,\tau,y):y\ge x,\ |y|<n,\ g(y,\tau)>q\}$, the function $(x,\tau)\mapsto\rho_\tau(\{y\in(-n,n):y\ge x,\ g(y,\tau)>q\})$ is Borel, and letting $n\to\infty$ shows that $(x,\tau)\mapsto\rho_\tau(\{y\ge x:g(y,\tau)>q\})$ is Borel. The same argument applies to $L$.
Since $\{H>a\}$ is the union of the sets $\{H>q\}$ over rational $q>a$, and $\{L<a\}$ is the union of the sets $\{L<q\}$ over rational $q<a$, the functions $H$ and $L$ are Borel on $\R\times\mathcal N$ and on $\mathcal N$, respectively.

Next we show that $\Gamma(x,\tau):=\max\{H(x,\tau),L(\tau)\}$ is finite. Fix $\tau\in\mathcal N$ and $x\in\R$.
If $\rho_\tau([x,\infty))>0$, choose $r>x$ with $\rho_\tau([x,r])>0$. Since $g(y,\tau)\le B(\tau)-x/\tau$ for $y\ge x$ and $g(y,\tau)\ge-B(\tau)-r/\tau$ for $y\in[x,r]$, we obtain $-B(\tau)-r/\tau\le H(x,\tau)\le B(\tau)-x/\tau$. Moreover, the essential infimum over $\R$ does not exceed the essential infimum over $[x,\infty)$, which does not exceed $H(x,\tau)$, so $L(\tau)\le H(x,\tau)$ and $\Gamma(x,\tau)=H(x,\tau)$.
If $\rho_\tau([x,\infty))=0$, then $H(x,\tau)=-\infty$, $\rho_\tau$ is concentrated on $(-\infty,x)$, and there is $r<x$ with $\rho_\tau([r,x))>0$. Since $g(y,\tau)\ge-B(\tau)-x/\tau$ for $y<x$ and $g(y,\tau)\le B(\tau)-r/\tau$ for $y\ge r$, we obtain $-B(\tau)-x/\tau\le L(\tau)\le B(\tau)-r/\tau$, and $\Gamma(x,\tau)=L(\tau)$.
In both cases $\Gamma(x,\tau)\in\R$. Moreover, $\Gamma(\cdot,\tau)$ is nonincreasing, because $H(\cdot,\tau)$ is.

Define $\Upsilon(x,\tau):=\max\{-B(\tau),\min\{B(\tau),x/\tau+\Gamma(x,\tau)\}\}$ for $\tau\in\mathcal N$, and $\Upsilon(x,\tau):=0$ for $\tau\in(0,T)\setminus\mathcal N$. Then $\Upsilon$ is Borel on $\R\times(0,T)$ and $|\Upsilon|\le B$.
For $\tau\in\mathcal N$ and $x_1<x_2$, the monotonicity of $\Gamma(\cdot,\tau)$ gives $x_2/\tau+\Gamma(x_2,\tau)-x_1/\tau-\Gamma(x_1,\tau)\le(x_2-x_1)/\tau$, and the truncation, which is nondecreasing and $1$-Lipschitz, preserves this inequality, as in the proof of Lemma~\ref{oleinik-forms}. Hence $\Upsilon(\cdot,\tau)$ satisfies \eqref{homogeneous-oleinik} for every $\tau\in(0,T)$.

We now show that $\Upsilon(\cdot,\tau)=h(\cdot,\tau)$ $\rho_\tau$-almost everywhere whenever $\tau\in\mathcal N$ and $w(\cdot,\tau)$ has a representative as in the statement.
For such $\tau$, the function $h(\cdot,\tau)$ satisfies statement (ii) of Lemma~\ref{oleinik-forms} with $\mu=\rho_\tau$ and $k=1/\tau$, so statement (i) gives a Borel function $v$ with $v=h(\cdot,\tau)$ $\rho_\tau$-almost everywhere, $|v|\le\|w(\cdot,\tau)\|_{L^\infty(d\rho_\tau)}\le B(\tau)$, and $v(x_2)-v(x_1)\le(x_2-x_1)/\tau$ for all $x_1<x_2$.
Then $g_v(y):=v(y)-y/\tau$ is nonincreasing on $\R$ and coincides with $g(\cdot,\tau)$ $\rho_\tau$-almost everywhere, so $H(\cdot,\tau)$ and $L(\tau)$ are the corresponding essential bounds of $g_v$.
Let $x\in\supp\rho_\tau$. Since $g_v\le g_v(x)$ on $[x,\infty)$, we have $H(x,\tau)\le g_v(x)$.
If $\rho_\tau(\{x\})>0$, then $H(x,\tau)\ge g_v(x)$, since $\{x\}$ has positive measure.
If $\rho_\tau((x,x+\delta))>0$ for every $\delta>0$ and $g_v$ is continuous at $x$, then $H(x,\tau)\ge g_v(x+\delta)$ for every $\delta>0$, because $g_v\ge g_v(x+\delta)$ on the set $(x,x+\delta)$ of positive measure, and letting $\delta\to0^+$ gives $H(x,\tau)\ge g_v(x)$.
In both cases $\rho_\tau([x,\infty))>0$, so $\Gamma(x,\tau)=H(x,\tau)=g_v(x)$ and $\Upsilon(x,\tau)=\max\{-B(\tau),\min\{B(\tau),v(x)\}\}=v(x)$.
The remaining points form a $\rho_\tau$-null set. Indeed, $\R\setminus\supp\rho_\tau$ is $\rho_\tau$-null, and the remaining points of $\supp\rho_\tau$ carry no atom and belong either to the countable set of discontinuity points of the monotone function $g_v$ or to the set $D$ of points $x\in\supp\rho_\tau$ such that $\rho_\tau((x,x+\delta_x))=0$ for some $\delta_x>0$. The set $D$ is countable, because, if $x<x'$ both belong to $D$, then $x'\in\supp\rho_\tau$ does not lie in the $\rho_\tau$-null open interval $(x,x+\delta_x)$, so the intervals $(x,x+\delta_x)$, $x\in D$, are pairwise disjoint. A countable set of points without atoms is $\rho_\tau$-null.
Consequently, $\Upsilon(\cdot,\tau)=v=h(\cdot,\tau)$ $\rho_\tau$-almost everywhere.

Finally, the set of times at which the conclusion holds must be selected measurably. Define
$$
\mathcal J:=\bigcap_{n\ge1}\Big\{\tau\in(0,T):\int_{(-n,n)}\min\{1,|\Upsilon(y,\tau)-h(y,\tau)|\}\,d\rho_\tau(y)=0\Big\},
$$
which is a Borel set by Lemma~\ref{kernel}(ii), applied with $\Upsilon$ extended by $0$ at $\tau=0$. It contains $(0,T)\setminus\mathcal N$, where $\rho_\tau=0$, and every $\tau\in\mathcal N$ at which $w(\cdot,\tau)$ has a representative as in the statement, so its complement in $(0,T)$ is a Lebesgue null set.
Put $\check w:=\Upsilon$ on $\R\times\mathcal J$ and $\check w:=h$ on $\R\times([0,T)\setminus\mathcal J)$. Then $\check w$ is Borel, $|\check w(x,\tau)|\le B(\tau)$ for $\tau>0$ and $|\check w(x,0)|\le\|w(\cdot,0)\|_{L^\infty(d\rho_0)}$, $\check w(\cdot,\tau)=w(\cdot,\tau)$ $\rho_\tau$-almost everywhere for every $\tau\in[0,T)$, and $\check w(\cdot,\tau)$ satisfies \eqref{homogeneous-oleinik} for every $\tau\in\mathcal J$.
\end{proof}

We first fix the setting. Let $T\in(0,\infty]$, let $(\rho,w)$ be an entropy solution of \eqref{H} on $[0,T)$ with initial data $(\rho_0,u_0)$, and fix $T_0\in(0,T)$. By Lemma~\ref{joint-velocity}, we may replace $w$ by a jointly Borel function, still denoted by $w$, such that $w(\cdot,\tau)\rho_\tau$ and $w(\cdot,\tau)^2\rho_\tau$ are unchanged for every $\tau\in[0,T)$, and \eqref{homogeneous-oleinik} holds for all $x_1<x_2$ and almost every $\tau\in(0,T)$. The weak identities and the conditions of Definitions~\ref{weak-solution} and~\ref{entropy-solution} depend on the velocity only through these measures, so $(\rho,w)$ is still an entropy solution. The same lemma shows that
$$
M:=\sup_{(x,\tau)\in\R\times[0,T_0]}|w(x,\tau)|
$$
is finite. We write $P_\tau:=w(\cdot,\tau)\rho_\tau$ for the momentum. Since $|P_\tau|\le M\rho_\tau$, condition (d) of Definition~\ref{weak-solution} and Definition~\ref{entropy-solution}(ii) give
\begin{equation}\label{initial-traces}
\rho_\tau\stackrel{*}{\rightharpoonup}\rho_0,\qquad
P_\tau\stackrel{*}{\rightharpoonup}u_0\rho_0,\qquad
w(\cdot,\tau)^2\rho_\tau\stackrel{*}{\rightharpoonup}u_0^2\rho_0
\qquad\text{as }\tau\to0^+,
\end{equation}
and $|u_0|\le M$ $\rho_0$-almost everywhere.

The generalized potentials are built with the oriented integrals $\int_a^yf\,d\mu$ of Theorem~\ref{main-introduction}. For $a\in\R$ and $y_1<y_2$, $\int_a^{y_2}f\,d\mu-\int_a^{y_1}f\,d\mu=\int_{[y_1,y_2)}f\,d\mu$, and $y\mapsto\int_a^yf\,d\mu$ is left-continuous. The same notation is used for integrals with respect to the Lebesgue measure. Fix $\bar a\in\R$ with $\rho_0(\{\bar a\})=0$. For $0\le s\le\tau\le T_0$ and $x,y\in\R$, define
$$
G_s(y;x,\tau):=(\tau-s)\int_{\bar a}^ydP_s+\int_{\bar a}^y(\eta-x)\,d\rho_s(\eta),
\qquad
\mathcal V_s(x,\tau):=\inf_{y\in\R}G_s(y;x,\tau),
$$
with $P_0:=u_0\rho_0$, and put $G:=G_0$ and $\mathcal V:=\mathcal V_0$. Thus $G(y;x,\tau)=\int_{\bar a}^y(\tau u_0(\eta)+\eta-x)\,d\rho_0(\eta)$ is the potential of Theorem~\ref{main-introduction}(iii) with $A(t)$ replaced by $\tau$, and $G_s$ is the analogous potential built from the state of the solution at time $s$. A change of base point adds to these potentials a quantity of the form $(\tau-s)c_1+c_2-xc_3$, with constants $c_1,c_2,c_3$, which is annihilated by $\partial_x^2$ and by $\partial_x\partial_\tau$. The function $\mathcal V$ depends only on the data and on $\bar a$.

The first two steps lead to the following representation of the solution at positive times.

\begin{proposition}\label{positive-representation}
For every $s\in(0,T_0)$,
$$
\rho=-\partial_x^2\mathcal V_s,\qquad P=\partial_x\partial_\tau\mathcal V_s
\qquad\text{in }\mathcal D'(\R\times(s,T_0)),
$$
where $\rho$ and $P$ denote the measures $\rho_\tau\,d\tau$ and $P_\tau\,d\tau$.
\end{proposition}

The proof is given in Section~\ref{lagrangian}. It involves the solution only on $[s,T_0]$, it uses the Oleinik condition only through the forward uniqueness of characteristics, and it does not use the energy condition. The third step is the convergence of $\mathcal V_s$ to $\mathcal V$ as $s\to0^+$, which we obtain from three lemmas. The first localizes the minimization, the second gives the upper bound, and the third gives the lower bound, which is the only place where the energy condition enters.

\begin{lemma}[Localization]\label{localization}
For $0\le s<\tau\le T_0$ and $x\in\R$, put $K:=[x-M\tau,x+M\tau]$. Then $\mathcal V_s(x,\tau)=\inf\limits_{y\in K}G_s(y;x,\tau)$. Moreover, $\mathcal V_s$ is concave in $(x,\tau)$, and there is $s_0>0$ such that $\sup\limits_{s\le s_0}\sup\limits_{(x,\tau)\in\mathcal Q}|\mathcal V_s(x,\tau)|<\infty$ for every compact $\mathcal Q\subset\R\times(0,T_0)$.
\end{lemma}

\begin{proof}
For $y_1<y_2\le x-M\tau$,
$$
G_s(y_2)-G_s(y_1)=\int_{[y_1,y_2)}\big((\tau-s)w(\eta,s)+\eta-x\big)\,d\rho_s(\eta)\le0,
$$
because $(\tau-s)w+\eta-x<M\tau-M\tau=0$ for $\eta<x-M\tau$, $\rho_s$-almost everywhere, where $w(\cdot,0):=u_0$ and $|u_0|\le M$ $\rho_0$-almost everywhere. Hence $G_s(\cdot;x,\tau)$ is nonincreasing on $(-\infty,x-M\tau]$. In the same way it is nondecreasing on $[x+M\tau,\infty)$, and the infimum over $\R$ is the infimum over $K$.
For fixed $y$, $G_s(y;x,\tau)=(\tau-s)\int_{\bar a}^ydP_s+\int_{\bar a}^y\eta\,d\rho_s-x\int_{\bar a}^yd\rho_s$ is affine in $(x,\tau)$, so $\mathcal V_s$ is an infimum of affine functions and is concave.
For the bound, let $\mathcal Q\subset[-R,R]\times[\tau_-,\tau_+]$ and $J:=[-R-MT_0-|\bar a|,R+MT_0+|\bar a|]$. By Lemma~\ref{uniform-bound}, $\sup\limits_{s\le s_0}\rho_s(J)<\infty$ for every $s_0<T_0$, and $|G_s(y;x,\tau)|\le(MT_0+|J|)\sup\limits_{s\le s_0}\rho_s(J)$ for $y\in K$ and $(x,\tau)\in\mathcal Q$. The infimum over $K$ obeys the same bound.
\end{proof}

\begin{lemma}[Upper bound]\label{upper}
For every $x$ and $\tau\in(0,T_0)$, $\limsup\limits_{s\to0^+}\mathcal V_s(x,\tau)\le\mathcal V(x,\tau)$.
\end{lemma}

\begin{proof}
Let $y$ satisfy $\rho_0(\{y\})=0$. Then the boundary of $[\bar a,y)$, or of $[y,\bar a)$, is $\rho_0$-null. The measures $\rho_s$ and $P_s$ converge weakly-$*$ by \eqref{initial-traces}, and $|P_s|\le M\rho_s$. Approximating the indicator function of the interval by continuous functions and estimating the error by $M\rho_s$ of a small neighborhood of the two endpoints, whose limit superior is small, we obtain $G_s(y;x,\tau)\to G(y;x,\tau)$. Hence $\limsup\mathcal V_s\le G(y;x,\tau)$ for every such $y$.
The points $y$ with $\rho_0(\{y\})=0$ are dense, and $G(\cdot;x,\tau)$ is left-continuous. Hence the infimum of $G$ over these points equals $\mathcal V(x,\tau)$.
\end{proof}

The lower bound requires the joint distribution of position and velocity. For a Borel map $f$ and a measure $\mu$, we denote by $f_\#\mu$ the image measure, $f_\#\mu(B):=\mu(f^{-1}(B))$, and by $\mu\llcorner E$ the restriction $\mu(\cdot\cap E)$ of $\mu$ to a Borel set $E$. For $s\in(0,T_0]$ let
$$
\pi_s:=(\mathrm{id},w(\cdot,s))_\#\rho_s,
$$
which is a nonnegative Radon measure on $\R\times[-M,M]$, and let $\pi_0:=(\mathrm{id},u_0)_\#\rho_0$, where $u_0$ is a Borel representative with $|u_0|\le M$.

\begin{lemma}[Zero variance in the limit]\label{zero-variance}
$\pi_s\stackrel{*}{\rightharpoonup}\pi_0$ as $s\to0^+$, in duality with $C_c(\R\times[-M,M])$.
\end{lemma}

\begin{proof}
The family $(\pi_s)_{s\le s_0}$ is bounded on every set $J\times[-M,M]$ with $J$ compact, so every sequence $s_k\to0^+$ has a subsequence with $\pi_{s_k}\stackrel{*}{\rightharpoonup}\pi$. For $\varphi\in C_c(\R)$ and $k=0,1,2$, the function $(\eta,v)\mapsto\varphi(\eta)v^k$ is continuous with compact support in $\R\times[-M,M]$, and \eqref{initial-traces} gives
$$
\int\varphi(\eta)v^k\,d\pi(\eta,v)=\int\varphi\,u_0^k\,d\rho_0,\qquad k=0,1,2.
$$
The case $k=0$ shows that the first marginal of $\pi$ is $\rho_0$. Disintegrating $\pi=\int\pi_\eta\,d\rho_0(\eta)$ with probability measures $\pi_\eta$ on $[-M,M]$, the cases $k=1,2$ give $\int v\,d\pi_\eta=u_0(\eta)$ and $\int v^2\,d\pi_\eta=u_0(\eta)^2$ for $\rho_0$-almost every $\eta$. Hence $\int(v-u_0(\eta))^2\,d\pi_\eta(v)=0$, so $\pi_\eta=\delta_{u_0(\eta)}$ and $\pi=\pi_0$. Since the limit does not depend on the subsequence, the whole family converges.
\end{proof}

\begin{lemma}[Lower bound]\label{lower}
For every $x$ and $\tau\in(0,T_0)$, $\liminf\limits_{s\to0^+}\mathcal V_s(x,\tau)\ge\mathcal V(x,\tau)$.
\end{lemma}

\begin{proof}
Choose $s_k\to0^+$ realizing the limit inferior and, by Lemma~\ref{localization}, points $y_k\in K$ with $G_{s_k}(y_k;x,\tau)\le\mathcal V_{s_k}(x,\tau)+s_k$. Passing to a subsequence, $y_k\to y^*\in K$. Put $g_s(\eta,v):=(\tau-s)v+\eta-x$ and $g(\eta,v):=\tau v+\eta-x$, so that $g_s\to g$ uniformly on compact sets and
$$
G_s(y;x,\tau)=\int_{\bar a}^y g_s\,d\pi_s,
$$
where the oriented integral is taken over $\{\eta\in[\bar a,y)\}$ or $\{\eta\in[y,\bar a)\}$.
Choose $b<\min\{\bar a,\min K\}$ with $\rho_0(\{b\})=0$. Then $G_s(y)=H_s(y)-H_s(\bar a)$ with $H_s(y):=\int_{\{b\le\eta<y\}}g_s\,d\pi_s$ for $y>b$, and the same identity holds for $G$ with $H(y):=\int_{\{b\le\eta<y\}}g\,d\pi_0$. Since $\rho_0(\{\bar a\})=\rho_0(\{b\})=0$, Lemma~\ref{zero-variance} gives $H_{s_k}(\bar a)\to H(\bar a)$, as in the proof of Lemma~\ref{upper}.

It remains to identify the limits of $H_{s_k}(y_k)$. Fix $\delta>0$ with $\rho_0(\{y^*-\delta\})=0$ and $y^*-\delta>b$. For $k$ large, $y_k>y^*-\delta$ and
$$
H_{s_k}(y_k)=H_{s_k}(y^*-\delta)+\int g_{s_k}\,d\sigma_k,
\qquad
\sigma_k:=\mathbf 1_{\{y^*-\delta\le\eta<y_k\}}\,\pi_{s_k}.
$$
The first term converges to $H(y^*-\delta)$. The measures $\sigma_k$ are supported in a fixed compact set and have bounded mass, so, after passing to a subsequence, $\sigma_k\stackrel{*}{\rightharpoonup}\sigma$, and $\int g_{s_k}\,d\sigma_k\to\int g\,d\sigma$.
We identify $\sigma$. For nonnegative $\phi\in C_c(\R\times[-M,M])$, $\int\phi\,d\sigma\le\limsup\int\phi\,\mathbf 1_{\{y^*-\delta\le\eta\le y^*+\delta\}}\,d\pi_{s_k}\le\int\phi\,\mathbf 1_{\{y^*-\delta\le\eta\le y^*+\delta\}}\,d\pi_0$, because the set is closed. Hence $\sigma\le\pi_0\llcorner\{y^*-\delta\le\eta\le y^*+\delta\}$. If $\phi$ is supported in $\{\eta>y^*\}$, then $\int\phi\,d\sigma_k=0$ for $k$ large, and if $\phi$ is supported in $\{y^*-\delta<\eta<y^*\}$, then $\int\phi\,d\sigma_k=\int\phi\,d\pi_{s_k}$ for $k$ large. Since $\rho_0(\{y^*-\delta\})=0$, it follows that
$$
\sigma=\pi_0\llcorner\{y^*-\delta<\eta<y^*\}+\sigma\llcorner\{\eta=y^*\}.
$$
The measure $\pi_0$ is concentrated on the graph of $u_0$, and $\pi_0\llcorner\{\eta=y^*\}=\rho_0(\{y^*\})\,\delta_{(y^*,u_0(y^*))}$. Therefore $\sigma\llcorner\{\eta=y^*\}=\theta\,\rho_0(\{y^*\})\,\delta_{(y^*,u_0(y^*))}$ for some $\theta\in[0,1]$, and
$$
\lim_{k\to\infty}H_{s_k}(y_k)=H(y^*)+\theta\,\rho_0(\{y^*\})\big(\tau u_0(y^*)+y^*-x\big)
=(1-\theta)H(y^*)+\theta H(y^*+0).
$$
Consequently,
$$
\liminf_{s\to0^+}\mathcal V_s(x,\tau)\ge(1-\theta)G(y^*;x,\tau)+\theta\,G(y^*+0;x,\tau)\ge\mathcal V(x,\tau),
$$
since $G(y^*+0;x,\tau)$ is a limit of values $G(y;x,\tau)$ with $y\downarrow y^*$.
\end{proof}

\begin{theorem}\label{convergence}
For every $x\in\R$ and $\tau\in(0,T_0)$, $\mathcal V_s(x,\tau)\to\mathcal V(x,\tau)$ as $s\to0^+$. Moreover, the functions $\mathcal V_s$ are concave and uniformly bounded on every compact subset of $\R\times(0,T_0)$, for $s$ small.
\end{theorem}

\begin{proof}
The pointwise convergence follows from Lemmas~\ref{upper} and~\ref{lower}, and the concavity and the uniform bounds from Lemma~\ref{localization}.
\end{proof}

We can now prove the main result of this subsection.

\begin{theorem}[Uniqueness on strips]\label{strip-uniqueness}
Let $T\in(0,\infty]$, let $\rho_0\in\M_{\mathrm{loc}}(\R)$ be nonnegative, and let $u_0\in L^\infty(\R,d\rho_0)$. If $(\rho^1,w^1)$ and $(\rho^2,w^2)$ are entropy solutions of \eqref{H} on $[0,T)$ with initial data $(\rho_0,u_0)$, then, for every $\tau\in[0,T)$,
$$
\rho^1_\tau=\rho^2_\tau,
\qquad
w^1(\cdot,\tau)\rho^1_\tau=w^2(\cdot,\tau)\rho^2_\tau
$$
as Radon measures.
\end{theorem}

\begin{proof}
Fix $T_0\in(0,T)$, replace $w^1$ and $w^2$ by the representatives described above, and let $M$ be the larger of the two corresponding bounds. Let $O\Subset\R\times(0,T_0)$ be open. For $s$ small, $O\subset\R\times(s,T_0)$, and Proposition~\ref{positive-representation} gives $\rho^i=-\partial_x^2\mathcal V^i_s$ and $P^i=\partial_x\partial_\tau\mathcal V^i_s$ in $\mathcal D'(O)$ for $i=1,2$, with the obvious notation. By Theorem~\ref{convergence}, both $\mathcal V^1_s$ and $\mathcal V^2_s$ converge pointwise to the same function $\mathcal V$, which depends only on the data, and they are uniformly bounded on $O$. By dominated convergence, $\mathcal V^i_s\to\mathcal V$ in $L^1(O)$, hence in $\mathcal D'(O)$ together with all derivatives. Since the left-hand sides do not depend on $s$, $\rho^i=-\partial_x^2\mathcal V$ and $P^i=\partial_x\partial_\tau\mathcal V$ in $\mathcal D'(O)$. Since every compact subset of $\R\times(0,T_0)$ is contained in such a set $O$, these identities hold in $\mathcal D'(\R\times(0,T_0))$. Hence $\rho^1_\tau=\rho^2_\tau$ and $P^1_\tau=P^2_\tau$ for almost every $\tau\in(0,T_0)$. For each $\varphi\in C_c(\R)$, the differences of the corresponding pairings are continuous functions of $\tau$ by condition (a) of Definition~\ref{weak-solution}, so the identities hold for every $\tau\in[0,T_0]$. Since $T_0<T$ is arbitrary, the theorem follows.
\end{proof}

Without the energy condition, the limit $\pi$ in Lemma~\ref{zero-variance} may have velocity dispersion over an atom of $\rho_0$. A part of $\pi$ over that atom may then select its lower velocities, and the corresponding limit in Lemma~\ref{lower} may lie strictly below both $G(y^*)$ and $G(y^*+0)$. This is the mechanism of the splitting solutions of Proposition~\ref{nonuniqueness}, in which the defect in the initial energy trace is the mass of the atom times the variance of its velocity distribution.

The argument does not follow individual particles from the initial time. In the passage to the initial time, Huang and Wang \cite[(4.21)--(4.36)]{HuangWang2001} use instead an exceptional set of initial points, which we now describe. For $\tau>0$, let $\gamma_\tau:\R\to\R$ be a nondecreasing map, where $\gamma_\tau(\eta)$ is the position at time $\tau$ of the particle issued from $\eta$. For $x\in\R$ and $\tau>0$, the smallest and the largest initial points whose particles reach $x$ at time $\tau$ are
$$
\xi^-(x,\tau):=\sup\{\eta\in\R:\gamma_\tau(\eta)<x\},
\qquad
\xi^+(x,\tau):=\inf\{\eta\in\R:\gamma_\tau(\eta)>x\},
$$
with $\sup\emptyset:=-\infty$ and $\inf\emptyset:=+\infty$, and the exceptional set is
$$
\Omega:=\bigcap_{\tau>0}\,\bigcup_{x\in\R}\big\{\eta\in\R:\xi^-(x,\tau)<\xi^+(x,\tau)\text{ and }\xi^-(x,\tau)\le\eta\le\xi^+(x,\tau)\big\}.
$$
Thus $\Omega$ consists of the initial points whose particles belong to a concentration at every positive time. The set used in \cite[Section~4]{HuangWang2001}, denoted there by $\Omega_s$, is the set $\Omega$ obtained with $\gamma_\tau:=\gamma(\cdot,\tau)$, where $\gamma$ is the limit of the approximate characteristics issued at time $0$, denoted by $X$ in \cite{HuangWang2001} \cite[(3.4), Lemma~3.1(1), and Lemma~3.2(4)]{HuangWang2001}. Huang and Wang deduce that $\Omega_s$ has Lebesgue measure zero from the countability of the curves of discontinuity, and the Lebesgue measure is the one used in the subsequent limit \cite[(4.23)--(4.24)]{HuangWang2001}. This countability does not imply that $\Omega$ is Lebesgue null. In Section~\ref{immediate-concentration} we give bounded initial data for which the solution of Theorem~\ref{HW} consists, at every positive time, of finitely many atoms, while $\Omega$ contains every point of $(0,1)$ that is not a dyadic rational, and therefore has Lebesgue measure at least $1$, both for the forward generalized characteristics of the homogeneous theory and for every family of nondecreasing maps with $(\gamma_\tau)_\#\rho_0=\rho_\tau$ (Proposition~\ref{dyadic-exceptional}). For the curves $\gamma$ themselves, $\Omega_s$ contains every point of $(0,1)$ outside a countable set (Proposition~\ref{HW-exceptional}). The example concerns an auxiliary set in the proof and is consistent with Theorem~\ref{strip-uniqueness}.

We can now prove uniqueness for \eqref{NA}.

\begin{theorem}[Uniqueness]
\label{uniqueness}
Let Assumption~\ref{coefficients} hold and let $T\in(0,\infty]$. Let $(\rho^1,u^1)$ and $(\rho^2,u^2)$ be entropy solutions of \eqref{NA} on $[0,T)$ with the same initial data. Then, for every $t\in[0,T)$,
$$
\rho^1_t=\rho^2_t,
\qquad
u^1(\cdot,t)\rho^1_t=u^2(\cdot,t)\rho^2_t
$$
as Radon measures, that is, $(\rho^1,u^1)\sim(\rho^2,u^2)$.
\end{theorem}

\begin{proof}
By Theorem~\ref{entropy-correspondence}, $(\tilde\rho^i,w^i):=\mathcal T(\rho^i,u^i)$, $i=1,2$, are entropy solutions of \eqref{H} on $[0,A(T))$ with the same initial data. By Theorem~\ref{strip-uniqueness}, $\tilde\rho^1_\tau=\tilde\rho^2_\tau$ and $w^1(\cdot,\tau)\tilde\rho^1_\tau=w^2(\cdot,\tau)\tilde\rho^2_\tau$ for every $\tau\in[0,A(T))$. Evaluating at $\tau=A(t)$ and multiplying the second identity by $e^{-S(t)}$ gives the assertion.
\end{proof}

\begin{theorem}[Global existence and uniqueness]
\label{existence-uniqueness}
Let Assumption~\ref{coefficients} hold. For every nonnegative $\rho_0\in\M_{\mathrm{loc}}(\R)$ and every $u_0\in L^\infty(\R,d\rho_0)$, problem \eqref{NA} admits a global entropy solution, which is unique modulo the relation of Definition~\ref{equiv-class}.
\end{theorem}

\begin{proof}
Combine Theorems~\ref{existence} and~\ref{uniqueness}.
\end{proof}

\subsection{Equivalence with the Huang--Wang class}
\label{class-equivalence}

The proofs above use the formulation of \cite{HuangWang2001} only for existence. We now show that this solution class coincides with the intrinsic one, up to the choice of a representative of the velocity, so that the homogeneous theory itself can be formulated without a potential.

\begin{proposition}[Equivalence with the Huang--Wang class]
\label{HW-equivalence}
Let $T\in(0,\infty]$ and let $(\rho,w)$ satisfy conditions (a) and (b) of Definition~\ref{weak-solution} for \eqref{H} on $[0,T)$.
Then $(\rho,w)$ is a weak solution of \eqref{H} on $[0,T)$ in the sense of Definition~\ref{weak-solution} if and only if there are a potential $m\in L^\infty_{\mathrm{loc}}([0,T);BV_{\mathrm{loc}}(\R))$ with $\partial_xm=\rho$ and a jointly Borel velocity $\check w$ with $\check w(\cdot,\tau)=w(\cdot,\tau)$ $\rho_\tau$-almost everywhere for every $\tau\in[0,T)$ such that $(m,\check w)$ is a weak solution of \eqref{H} on $[0,T)$ in the strip-local sense of \cite[Definition~1]{HuangWang2001} specified above.
In that case, $m$ may be chosen with the initial normalization of Proposition~\ref{intrinsic-potential}, $\check w$ may be chosen as the representative $h$ of Lemma~\ref{borel-velocity}, and \eqref{HW-weak} holds for $(m,\check w)$ and all $\varphi,\psi\in C_c^1(\R\times(0,T))$.
The same equivalence holds for entropy solutions in the sense of Definition~\ref{entropy-solution} and of \cite[Definition~2]{HuangWang2001}, with $\check w$ chosen as the jointly Borel admissible representative of Lemma~\ref{joint-velocity}. The representative $h$ need not satisfy \eqref{HW-Oleinik} in vacuum.
\end{proposition}

\begin{proof}
The weak equations, measure traces, and velocity bounds depend only on the $\rho_\tau$-almost everywhere velocity class at each time. The pointwise Oleinik inequality will be handled by the admissible representative in Lemma~\ref{joint-velocity}.
For a jointly Borel representative bounded on finite strips, all the compactly supported integrands in \eqref{HW-weak} are Borel and absolutely integrable for the locally finite kernel measure $\mu(dx\,d\tau)=\rho_\tau(dx)\,d\tau$. Indeed, their absolute integrals are bounded by the spatial local mass bound times a strip velocity bound and the length of the time support. Lemma~\ref{kernel}(iii), applied to positive and negative parts, therefore identifies their Lebesgue--Stieltjes integrals with the iterated integrals in Definition~\ref{weak-solution}.

Let $(\rho,w)$ be a weak solution in the sense of Definition~\ref{weak-solution}, and let $\check w:=h$ be given by Lemma~\ref{borel-velocity}.
Then $(\rho,\check w)$ is also a weak solution in the sense of Definition~\ref{weak-solution}, and $\check w$ has the measurability and boundedness required in the definition recalled in Section~\ref{bridge}, on every strip $\R\times(0,\tau_1]$ with $\tau_1<T$.
By Corollary~\ref{intrinsic-cumulative}, there is a potential $m$ with $\partial_xm=\rho$ and $\partial_\tau m+\check w\rho=0$ in $\mathcal D'(\R\times(0,T))$.
Testing the latter identity with $\varphi\in C_c^\infty(\R\times(0,T))$ and using the coincidence of the integrals just established gives the first identity in \eqref{HW-weak}.
By the same coincidence, the second identity in \eqref{HW-weak} is \eqref{weak2-H} with $Z=\psi$.
The continuity in time and the initial conditions are contained in Definition~\ref{weak-solution}.

Conversely, let $m$ and $\check w$ be as in the statement.
By the coincidence of the integrals established above, the first identity in \eqref{HW-weak} states that $\partial_\tau m+\check w\rho=0$ in $\mathcal D'(\R\times(0,T))$, and Corollary~\ref{intrinsic-cumulative} gives \eqref{weak1-H} for $(\rho,\check w)$.
The momentum identity, the continuity in time, and the initial conditions coincide in the two formulations, so $(\rho,\check w)$, and therefore $(\rho,w)$, is a weak solution in the sense of Definition~\ref{weak-solution}.

For the extension to $C^1$ test functions, the second identity in \eqref{HW-weak} is \eqref{weak2-H}, so it extends by Lemma~\ref{extension}.
In the first identity, $m$ is locally integrable on $\R\times[0,T)$, since $m\in L^\infty_{\mathrm{loc}}([0,T);BV_{\mathrm{loc}}(\R))$, and $\check w\rho$ is a locally finite measure on $\R\times(0,T)$ by Lemma~\ref{uniform-bound}.
If $\varphi_\varepsilon$ is a standard mollification of $\varphi\in C_c^1(\R\times(0,T))$, then $\partial_\tau\varphi_\varepsilon\to\partial_\tau\varphi$ and $\varphi_\varepsilon\to\varphi$ uniformly as $\varepsilon\to0^+$, with supports in a fixed compact set, and both integrals in the first identity pass to the limit.

For entropy solutions, choose instead the jointly Borel admissible representative $\check w$ of Lemma~\ref{joint-velocity}, which is bounded on every strip $\R\times(0,\tau_1]$ with $\tau_1<T$, so that the preceding arguments apply to it. It satisfies the literal pointwise inequality \eqref{HW-Oleinik} for almost every time, and its energy measure equals that of $w$ at every time. Conversely, a pointwise Oleinik representative in the homogeneous class is an admissible representative for the intrinsic definition. This proves the entropy equivalence without asserting invariance of the pointwise inequality under arbitrary modifications in vacuum.
\end{proof}

Combined with Theorem~\ref{strip-uniqueness}, Proposition~\ref{HW-equivalence} shows that two entropy solutions of \eqref{H} on $[0,T)$ in the strip-local sense of Section~\ref{bridge}, with the same initial data and with jointly Borel velocities bounded on every strip $\R\times[0,T_0]$, $T_0<T$, have the same density and the same momentum at every time.

\section{Lagrangian description of entropy solutions}
\label{lagrangian}

In this section we prove Proposition~\ref{positive-representation}. Huang and Wang derive the corresponding representation \cite[(4.1)--(4.20)]{HuangWang2001} from the uniqueness theorem for the linear equation $\partial_\tau\omega+w\,\partial_x\omega=0$ \cite[Theorem~3.1]{HuangWang2001}, whose proof uses an approximation of the velocity due to Wang and Ding \cite{WangDing1997}, and from an evaluation of a potential along generalized characteristics. We follow the same plan with two changes. The linear uniqueness theorem is replaced by a description of the solution in mass coordinates, which is elementary in one space dimension and identifies the motion of every concentration directly. The evaluation along characteristics is replaced by an inequality along straight lines, which is first obtained along almost every line parallel to a given segment and then on the segment itself by continuity.
The description is Lagrangian in the sense of continuum mechanics. Instead of the density and the velocity at a fixed point $x$, we follow each portion of mass, identified by its mass coordinate $\ell$, which measures the mass lying to its left, and we describe the solution by the position $Y(\ell,\tau)$ of that portion of mass at time $\tau$. Lagrangian coordinates for sticky particles appear in E, Rykov, and Sinai \cite{ERykovSinai}, who take the initial position as Lagrangian coordinate, in Natile and Savar\'e \cite{NatileSavare2009}, through the monotone rearrangement of the density, and in Brenier, Gangbo, Savar\'e, and Westdickenberg \cite{BGSW2013}, who write the system in mass coordinates.
Throughout this section, $(\rho,w)$, $T_0$, $M$, $P_\tau$, $\bar a$, $G_s$, and $\mathcal V_s$ are as in Section~\ref{strips}, and we fix $s\in(0,T_0)$. In particular, $w$ is the jointly Borel function provided by Lemma~\ref{joint-velocity}. It satisfies $|w(x,\tau)|\le M$ at every point $(x,\tau)\in\R\times[0,T_0]$, and there is a Lebesgue null set of times in $(0,T)$ outside which $w(\cdot,\tau)$ satisfies \eqref{homogeneous-oleinik} for all $x_1<x_2$, including points of vacuum. These pointwise properties are not shared by an arbitrary representative of the velocity, and they are used in Lemmas~\ref{forward-uniqueness}, \ref{label-velocity}, and~\ref{label-momentum}.

\subsection{Mass coordinates}
\label{mass-coordinates}

We first describe an entropy solution through the positions of its mass coordinates. The weak formulation is first recorded on closed time intervals.

\begin{lemma}\label{closed-intervals}
Let $\mu_\tau$ denote either $\rho_\tau$ or $P_\tau$. For $0<\tau_1<\tau_2\le T_0$ and every $\Phi\in C^1(\R\times[\tau_1,\tau_2])$ that vanishes outside $K\times[\tau_1,\tau_2]$ for a compact set $K$,
\begin{equation}\label{closed-identity}
\int_{\tau_1}^{\tau_2}\!\!\int_\R(\partial_\tau\Phi+w\partial_x\Phi)\,d\mu_\tau\,d\tau=\int_\R\Phi(\cdot,\tau_2)\,d\mu_{\tau_2}-\int_\R\Phi(\cdot,\tau_1)\,d\mu_{\tau_1}.
\end{equation}
\end{lemma}

\begin{proof}
Let $\vartheta\in C^1(\R)$ be nondecreasing with $\vartheta=0$ on $(-\infty,1]$ and $\vartheta=1$ on $[2,\infty)$, and put $\chi_k(\tau):=\vartheta(k(\tau-\tau_1))\,\vartheta(k(\tau_2-\tau))$. For $k$ large, $\chi_k\Phi\in C_c^1(\R\times(\tau_1,\tau_2))$. By Lemma~\ref{extension}, applied with $\alpha\equiv1$ and $\sigma\equiv0$, the identities \eqref{weak1-H} and \eqref{weak2-H} hold for test functions in $C_c^1(\R\times(0,T))$, and \eqref{weak2-H} is the identity for $P$, since $w\,dP_\tau=w^2\,d\rho_\tau$. Hence
$$
\int_{\tau_1}^{\tau_2}\chi_k\int_\R(\partial_\tau\Phi+w\partial_x\Phi)\,d\mu_\tau\,d\tau+\int_{\tau_1}^{\tau_2}\chi_k'(\tau)\int_\R\Phi(\cdot,\tau)\,d\mu_\tau\,d\tau=0.
$$
The inner integral in the first term is bounded by a constant times $\sup\limits_\tau\rho_\tau(K)$, which is finite by Lemma~\ref{uniform-bound}, so the first term converges to the left-hand side of \eqref{closed-identity}. The function $\tau\mapsto\int\Phi(\cdot,\tau)\,d\mu_\tau$ is continuous on $[\tau_1,\tau_2]$ by Lemma~\ref{param-continuity}, applied with $J=[\tau_1,\tau_2]$ and the local bounds of Lemma~\ref{uniform-bound}. Since $\chi_k'$ is a nonnegative bump of integral $1$ concentrating at $\tau_1$ minus a nonnegative bump of integral $1$ concentrating at $\tau_2$, the second term converges to $\int\Phi(\cdot,\tau_1)\,d\mu_{\tau_1}-\int\Phi(\cdot,\tau_2)\,d\mu_{\tau_2}$.
\end{proof}

We call a family $(\mu_\tau)_{\tau\in[s,T_0]}$ of nonnegative measures in $\M_{\mathrm{loc}}(\R)$ a \emph{transport family} if $\tau\mapsto\mu_\tau$ is weakly-$*$ continuous and \eqref{closed-identity} holds for $s\le\tau_1<\tau_2\le T_0$ and all $\Phi$ as in Lemma~\ref{closed-intervals}. By Lemma~\ref{closed-intervals}, $(\rho_\tau)$ is a transport family. Every transport family satisfies $\sup\limits_\tau\mu_\tau(K)<\infty$ for compact $K$, by the argument of Lemma~\ref{uniform-bound}, and the integrals in \eqref{closed-identity} are well defined by Lemma~\ref{kernel}(ii), whose proof applies without change to families indexed by $[s,T_0]$.

The Oleinik condition enters the proof only through the following lemma.

\begin{lemma}[Forward uniqueness]\label{forward-uniqueness}
Let $s\le t_0<T_0$, and let $\gamma_1,\gamma_2:[t_0,T_0]\to\R$ be absolutely continuous with $\gamma_i'(\tau)=w(\gamma_i(\tau),\tau)$ for almost every $\tau\in(t_0,T_0)$. If $\gamma_1(t_1)=\gamma_2(t_1)$ for some $t_1\in[t_0,T_0]$, then $\gamma_1=\gamma_2$ on $[t_1,T_0]$. If $\gamma_1(t_0)\le\gamma_2(t_0)$, then $\gamma_1\le\gamma_2$ on $[t_0,T_0]$.
\end{lemma}

\begin{proof}
Put $\delta:=\gamma_2-\gamma_1$, which is absolutely continuous. By the choice of $w$ recalled at the beginning of this section, the null set of times outside which \eqref{homogeneous-oleinik} holds for $w(\cdot,\tau)$ and all pairs of points does not depend on $\gamma_1$ and $\gamma_2$. Hence, for almost every $\tau$, both equations hold at $\tau$ and \eqref{homogeneous-oleinik} holds for $w(\cdot,\tau)$ and all pairs of points. At such a $\tau$, if $\delta(\tau)>0$, then \eqref{homogeneous-oleinik} with $x_1=\gamma_1(\tau)$ and $x_2=\gamma_2(\tau)$ gives $\delta'(\tau)\le\delta(\tau)/\tau\le\delta(\tau)/s$, and if $\delta(\tau)<0$, then in the same way $-\delta'(\tau)\le-\delta(\tau)/s$. Moreover, $\delta'=0$ almost everywhere on $\{\delta=0\}$. Hence $|\delta|'\le|\delta|/s$ almost everywhere, and Gronwall's lemma gives $|\delta(\tau)|\le|\delta(t_1)|e^{(\tau-t_1)/s}$ for $\tau\ge t_1$. This proves the first assertion.
For the second, suppose that $\delta(t_2)<0$ for some $t_2>t_0$, and let $t_*$ be the largest $\tau\in[t_0,t_2]$ with $\delta(\tau)\ge0$. By continuity $\delta(t_*)=0$, and the first assertion gives $\delta=0$ on $[t_*,T_0]$, which contradicts $\delta(t_2)<0$.
\end{proof}

The next lemma constructs the cumulative mass function of a transport family.

\begin{lemma}\label{cumulative-function}
Let $(\mu_\tau)$ be a transport family. Fix $\theta\in C_c^\infty(\R)$ with $\theta\ge0$ and $\int_\R\theta\,dx=1$, as in the proof of Proposition~\ref{intrinsic-potential}, and define, for $x\in\R$ and $\tau\in[s,T_0]$,
$$
m_\mu(x,\tau):=\int_\R\theta(y)\int_y^xd\mu_\tau\,dy-\int_s^\tau\!\!\int_\R\theta\,w\,d\mu_r\,dr,
$$
where $\int_y^xd\mu_\tau$ is the oriented integral of Theorem~\ref{main-introduction}. Then $m_\mu$ is locally bounded and Borel on $\R\times[s,T_0]$, the function $m_\mu(\cdot,\tau)$ is nondecreasing and left-continuous with $m_\mu(x_2,\tau)-m_\mu(x_1,\tau)=\mu_\tau([x_1,x_2))$ for $x_1<x_2$, and
\begin{equation}\label{cumulative-identity}
\int_\R\Phi(\cdot,\tau_2)\,m_\mu(\cdot,\tau_2)\,dx-\int_\R\Phi(\cdot,\tau_1)\,m_\mu(\cdot,\tau_1)\,dx=\int_{\tau_1}^{\tau_2}\Big(\int_\R\partial_\tau\Phi\,m_\mu\,dx-\int_\R\Phi\,w\,d\mu_\tau\Big)\,d\tau
\end{equation}
for $s\le\tau_1<\tau_2\le T_0$ and all $\Phi$ as in Lemma~\ref{closed-intervals}.
\end{lemma}

\begin{proof}
For every $y$, $\int_y^{x_2}d\mu_\tau-\int_y^{x_1}d\mu_\tau=\mu_\tau([x_1,x_2))$, which gives the formula for the increments. If $x_n$ increases to $x$, then $\int_y^{x_n}d\mu_\tau\to\int_y^xd\mu_\tau$ for every $y$, and these integrals are bounded by $\mu_\tau(K)$ for $y\in\supp\theta$, where $K$ is a compact interval containing $\supp\theta$ and all the points $x_n$. The dominated convergence theorem therefore shows that $m_\mu(\cdot,\tau)$ is left-continuous. Local boundedness follows from $\sup\limits_\tau\mu_\tau(K)<\infty$ and $|w|\le M$.
For $\varphi\in C_c(\R)$ put
$$
g_\varphi(z):=\int_\R\!\int_\R\varphi(x)\theta(y)\big(\mathbf 1_{\{y\le z<x\}}-\mathbf 1_{\{x\le z<y\}}\big)\,dx\,dy.
$$
Since $\int_y^xd\mu_\tau=\int(\mathbf 1_{\{y\le z<x\}}-\mathbf 1_{\{x\le z<y\}})\,d\mu_\tau(z)$, Fubini's theorem gives
$$
\int_\R\varphi(x)\int_\R\theta(y)\int_y^xd\mu_\tau\,dy\,dx=\int_\R g_\varphi\,d\mu_\tau.
$$ Writing $\theta_-(z):=\int_{(-\infty,z]}\theta$ and $\Phi_+(z):=\int_{(z,\infty)}\varphi$, we have $g_\varphi=\theta_-\Phi_+-(1-\theta_-)\big(\int\varphi\,dx-\Phi_+\big)$, so $g_\varphi\in C_c^1(\R)$ and $g_\varphi'=(\int\varphi)\theta-\varphi$. In the same way, $\int\theta(y)\int_y^xd\mu_\tau\,dy=\int\kappa_x\,d\mu_\tau$ with $\kappa_x(z):=\mathbf 1_{\{z<x\}}\theta_-(z)-\mathbf 1_{\{z\ge x\}}(1-\theta_-(z))$, which is bounded, Borel, and vanishes for $z$ outside a compact set. By Lemma~\ref{kernel}(ii), applied to $\kappa_x$ and to the bounded Borel function $(y,r)\mapsto\theta(y)w(y,r)$, the function $m_\mu(x,\cdot)$ is Borel for every $x$. Since $m_\mu$ is, in addition, left-continuous in $x$, it is Borel on $\R\times[s,T_0]$. Indeed, $m_\mu(x,\tau)$ is the limit as $n\to\infty$ of $m_\mu\big((\lceil2^nx\rceil-1)2^{-n},\tau\big)$, and for each $n$ and each integer $j$ this function coincides on the Borel set $\{x:\lceil2^nx\rceil=j\}$ with the Borel function $\tau\mapsto m_\mu((j-1)2^{-n},\tau)$.
Now let $\Phi$ be as in Lemma~\ref{closed-intervals}, and put $g(z,\tau):=g_{\Phi(\cdot,\tau)}(z)$. Then $g$ is $C^1$, vanishes outside a set $K'\times[\tau_1,\tau_2]$ with $K'$ compact, $\partial_\tau g=g_{\partial_\tau\Phi}$, and $\partial_zg=(\int\Phi(\cdot,\tau)\,dx)\theta-\Phi(\cdot,\tau)$. By \eqref{closed-identity},
$$
\Big[\int g\,d\mu_\tau\Big]_{\tau_1}^{\tau_2}=\int_{\tau_1}^{\tau_2}\Big(\int g_{\partial_\tau\Phi}\,d\mu_\tau+\int\Phi(\cdot,\tau)\,dx\int\theta\,w\,d\mu_\tau-\int\Phi\,w\,d\mu_\tau\Big)\,d\tau.
$$
The function $C_\mu(\tau):=\int_s^\tau\int\theta\,w\,d\mu_r\,dr$ is Lipschitz, so
$$
\Big[C_\mu(\tau)\int\Phi(\cdot,\tau)\,dx\Big]_{\tau_1}^{\tau_2}=\int_{\tau_1}^{\tau_2}\Big(C_\mu(\tau)\int\partial_\tau\Phi\,dx+\int\Phi(\cdot,\tau)\,dx\int\theta\,w\,d\mu_\tau\Big)\,d\tau.
$$
Subtracting the second identity from the first and using $\int\partial_\tau\Phi m_\mu\,dx=\int g_{\partial_\tau\Phi}\,d\mu_\tau-H\int\partial_\tau\Phi\,dx$ gives \eqref{cumulative-identity}.
\end{proof}

Since $m_\mu(\cdot,\tau)$ is left-continuous and nondecreasing with Stieltjes measure $\mu_\tau$, we have $\int\partial_x\varphi\,m_\mu(\cdot,\tau)\,dx=-\int\varphi\,d\mu_\tau$ for $\varphi\in C_c^1(\R)$. We now introduce mass coordinates. We write $\Leb$ for the Lebesgue measure on $\R$.

\begin{lemma}\label{mass-labels}
Let $(\mu_\tau)$ be a transport family with $\mu_s\ne0$, and write $m:=m_\mu$.
\begin{enumerate}[label=(\roman*)]
\item For $s\le\tau_1\le\tau_2\le T_0$ and $x\in\R$, with $h:=\tau_2-\tau_1$,
$$
m(x-Mh,\tau_2)\le m(x,\tau_1)\le m(x+Mh,\tau_2).
$$
In particular, the open interval $I:=\big(\inf\limits_xm(x,\tau),\sup\limits_xm(x,\tau)\big)$ does not depend on $\tau$.
\item For $\ell\in I$ and $\tau\in[s,T_0]$, the number $Y(\ell,\tau):=\inf\{x\in\R:m(x,\tau)>\ell\}$ is finite, and for every $x\in\R$
$$
Y(\ell,\tau)<x\quad\text{if and only if}\quad\ell<m(x,\tau).
$$
The function $Y(\cdot,\tau)$ is nondecreasing and right-continuous, and $Y(\cdot,\tau)_\#(\Leb\llcorner I)=\mu_\tau$.
\item $|Y(\ell,\tau_2)-Y(\ell,\tau_1)|\le M|\tau_2-\tau_1|$, and $Y$ is Borel on $I\times[s,T_0]$.
\item Let $J$ be an open interval and $Z:J\to\R$ a nondecreasing right-continuous function with $Z_\#(\Leb\llcorner J)=\mu_\tau$. Then there is $c\in\R$ such that $J=I+c$ and $Z(\kappa)=Y(\kappa-c,\tau)$ for every $\kappa\in J$.
\end{enumerate}
\end{lemma}

\begin{proof}
(i) Let $\varphi\in C_c^1(\R)$ be nonnegative and put $\Phi(x,\tau):=\varphi(x-M(\tau-\tau_1))$. Then $\partial_\tau\Phi=-M\partial_x\Phi$, so $\int\partial_\tau\Phi m\,dx=M\int\Phi\,d\mu_\tau$, and \eqref{cumulative-identity} gives
$$
\int\varphi(x)\,m(x+Mh,\tau_2)\,dx-\int\varphi\,m(\cdot,\tau_1)\,dx=\int_{\tau_1}^{\tau_2}\!\!\int\Phi\,(M-w)\,d\mu_\tau\,d\tau\ge0.
$$
Hence $m(x+Mh,\tau_2)\ge m(x,\tau_1)$ for almost every $x$, and for every $x$ because both sides are left-continuous in $x$. The choice $\Phi(x,\tau):=\varphi(x+M(\tau-\tau_1))$ gives the other inequality in the same way, with $-(M+w)\le0$ in place of $M-w$. Letting $x\to\pm\infty$ in both inequalities shows that the infimum and the supremum of $m(\cdot,\tau)$ do not depend on $\tau$. In particular, the total mass $\mu_\tau(\R)=\sup m(\cdot,\tau)-\inf m(\cdot,\tau)$ does not depend on $\tau$.

(ii) Since $\ell<\sup m(\cdot,\tau)$, the set $\{m(\cdot,\tau)>\ell\}$ is nonempty, and since $\ell>\inf m(\cdot,\tau)$, there is $x'$ with $m\le m(x',\tau)<\ell$ on $(-\infty,x']$. Hence $Y(\ell,\tau)$ is finite. If $Y(\ell,\tau)<x$, there is $x'<x$ with $m(x',\tau)>\ell$, so $m(x,\tau)>\ell$. Conversely, if $\ell<m(x,\tau)$, left continuity gives $x'<x$ with $m(x',\tau)>\ell$, so $Y(\ell,\tau)\le x'<x$. Monotonicity is clear. If $\ell_n$ decreases to $\ell$ and $x>Y(\ell,\tau)$, then $\ell<m(x',\tau)$ for some $x'<x$, hence $\ell_n<m(x',\tau)$ for $n$ large and $Y(\ell_n,\tau)\le x'$, which proves right continuity. Finally, for $x_1<x_2$,
$$
\Leb\big(\{\ell\in I:x_1\le Y(\ell,\tau)<x_2\}\big)=\Leb\big(\{\ell\in I:m(x_1,\tau)\le\ell<m(x_2,\tau)\}\big)=\mu_\tau([x_1,x_2)),
$$
since the values of $m(\cdot,\tau)$ lie in the closure of $I$. The intervals $[x_1,x_2)$ form a $\pi$-system that generates the Borel sets, both measures are finite on each of them, and $\R$ is a countable union of such intervals, so $Y(\cdot,\tau)_\#(\Leb\llcorner I)=\mu_\tau$.

(iii) If $m(x,\tau_1)>\ell$, then $m(x+Mh,\tau_2)>\ell$ by (i), so $Y(\ell,\tau_2)\le x+Mh$, and taking the infimum over $x$ gives $Y(\ell,\tau_2)\le Y(\ell,\tau_1)+Mh$. The inequality $Y(\ell,\tau_1)\le Y(\ell,\tau_2)+Mh$ follows in the same way from the first inequality in (i). By (ii), $\{Y<x\}=\bigcup_{q\in\mathbb Q,\,q<x}\{(\ell,\tau):\ell<m(q,\tau)\}$, which is Borel because $m(q,\cdot)$ is Borel.

(iv) For $x\in\R$, the set $\{\kappa\in J:Z(\kappa)<x\}$ is an initial segment of $J$, and it is open on the right because $Z$ is right-continuous. Hence it equals $J\cap(-\infty,b(x))$ for some $b(x)$ in the closure of $J$ in $[-\infty,\infty]$, with $b(x)=\inf J$ if the set is empty. For $x_1<x_2$ the set $\{x_1\le Z<x_2\}=J\cap[b(x_1),b(x_2))$ has measure $\mu_\tau([x_1,x_2))<\infty$. If $b(x_1)=-\infty$ for some $x_1$, then $\inf J=-\infty$, and choosing $x_2$ with $b(x_2)>-\infty$, which exists because $Z$ is real-valued, gives a set of infinite measure. In the same way $b(x)<\infty$ for every $x$. Thus $b$ is finite and $b(x_2)-b(x_1)=\mu_\tau([x_1,x_2))=m(x_2,\tau)-m(x_1,\tau)$, so $b-m(\cdot,\tau)$ is a constant $c$. Consequently, for $\kappa\in J$, $Z(\kappa)<x$ if and only if $\kappa-c<m(x,\tau)$, which by (ii) means $Y(\kappa-c,\tau)<x$ whenever $\kappa-c\in I$. Since $Z$ is real-valued, every $\kappa\in J$ satisfies $b(x)\le\kappa<b(x')$ for suitable $x,x'$, so $J\subset[\inf b,\sup b]$, and $J\subset(\inf b,\sup b)$ because $J$ is open. Conversely, every point of $(\inf b,\sup b)$ lies strictly between two values of $b$, which belong to the closure of the interval $J$, so it belongs to $J$. Hence $J=(\inf b,\sup b)=I+c$. Hence $Z(\kappa)=Y(\kappa-c,\tau)$ for $\kappa\in J$.
\end{proof}

We call the points $\ell\in I$ \emph{mass coordinates}, and $Y(\ell,\tau)$ the position at time $\tau$ of the portion of mass with coordinate $\ell$. The number $\ell$ is the mass lying to the left of the corresponding portion of mass, up to the additive constant fixed by $\theta$, and it does not change in time. The next lemma identifies the velocity of every mass coordinate. It contains the fact that a concentration moves with the velocity assigned to it at the point of concentration. Its proof does not use the Oleinik condition. The velocity of a concentration is identified through the monotonicity of $Y(\cdot,\tau)$.

\begin{lemma}\label{label-velocity}
In the setting of Lemma~\ref{mass-labels}, for almost every $\ell\in I$ the function $Y(\ell,\cdot)$ satisfies $\partial_\tau Y(\ell,\tau)=w(Y(\ell,\tau),\tau)$ for almost every $\tau\in(s,T_0)$.
\end{lemma}

\begin{proof}
Fix $\varphi\in C_c^1(\R)$ with support in $[c_1,c_2]$, fix $\ell_0\in I$, and put, for $\ell\in I$ and $z\in\R$,
$$
k(\ell,z):=\mathbf 1_{\{\ell\ge\ell_0\}}\int_{(z,\infty)}\varphi\,dx-\mathbf 1_{\{\ell<\ell_0\}}\int_{(-\infty,z]}\varphi\,dx.
$$
By Lemma~\ref{mass-labels}(ii), $m(x,\tau)-\ell_0=\int_I\big(\mathbf 1_{\{\ell\ge\ell_0,\ Y(\ell,\tau)<x\}}-\mathbf 1_{\{\ell<\ell_0,\ Y(\ell,\tau)\ge x\}}\big)\,d\ell$, and Fubini's theorem gives
$$
\int_\R\varphi\,\big(m(\cdot,\tau)-\ell_0\big)\,dx=\int_Ik(\ell,Y(\ell,\tau))\,d\ell.
$$
Let $J_0:=\big[m(c_1-M(T_0-s),s),\,m(c_2+M(T_0-s),s)\big)$, which may be empty. By Lemma~\ref{mass-labels}(ii) and (iii), for $\ell\in I\setminus J_0$ the curve $Y(\ell,\cdot)$ stays on one side of $[c_1,c_2]$, so $\tau\mapsto k(\ell,Y(\ell,\tau))$ is constant. Moreover, $k(\ell,Y(\ell,\tau))=0$ for every $\tau$ if $\ell\ge\max\{\ell_0,m(c_2+M(T_0-s),s)\}$ or $\ell<\min\{\ell_0,m(c_1-M(T_0-s),s)\}$, so all the integrals over $I$ below are integrals over a bounded interval. For $\ell\in I$ and $\tau\in[s,T_0]$, we denote by $\partial_\tau Y(\ell,\tau)$ the upper limit of the difference quotients of $Y(\ell,\cdot)$ at $\tau$ over rational increments $h$ with $\tau+h\in[s,T_0]$. Since $Y$ is Borel and continuous in $\tau$, $\partial_\tau Y$ is Borel and bounded by $M$, the set of $(\ell,\tau)$ such that $Y(\ell,\cdot)$ is differentiable at $\tau$ is Borel, and $\partial_\tau Y(\ell,\tau)$ is the derivative at every such point. For every $\ell$ the function $\tau\mapsto k(\ell,Y(\ell,\tau))$ is Lipschitz, and since $\partial_zk(\ell,z)=-\varphi(z)$, its derivative is $-\varphi(Y(\ell,\tau))\,\partial_\tau Y(\ell,\tau)$ wherever $Y(\ell,\cdot)$ is differentiable at $\tau$. Therefore
$$
\Big[\int_Ik(\ell,Y(\ell,\tau))\,d\ell\Big]_{\tau_1}^{\tau_2}=-\int_{\tau_1}^{\tau_2}\!\!\int_I\varphi(Y(\ell,\tau))\,\partial_\tau Y(\ell,\tau)\,d\ell\,d\tau.
$$
On the other hand, \eqref{cumulative-identity} with $\Phi(x,\tau)=\varphi(x)$ and Lemma~\ref{mass-labels}(ii) give
$$
\Big[\int_\R\varphi\,\big(m(\cdot,\tau)-\ell_0\big)\,dx\Big]_{\tau_1}^{\tau_2}=-\int_{\tau_1}^{\tau_2}\!\!\int_\R\varphi\,w\,d\mu_\tau\,d\tau=-\int_{\tau_1}^{\tau_2}\!\!\int_I\varphi(Y(\ell,\tau))\,w(Y(\ell,\tau),\tau)\,d\ell\,d\tau.
$$
Let $e(\ell,\tau):=\partial_\tau Y(\ell,\tau)-w(Y(\ell,\tau),\tau)$ where $Y(\ell,\cdot)$ is differentiable at $\tau$ and $e:=0$ elsewhere, so that $e$ is Borel and $|e|\le2M$. Comparing the two identities, $\int_{\tau_1}^{\tau_2}\int_I\varphi(Y(\ell,\tau))\,e(\ell,\tau)\,d\ell\,d\tau=0$ for all $\tau_1<\tau_2$, so $\int_I\varphi(Y(\ell,\tau))\,e(\ell,\tau)\,d\ell=0$ for almost every $\tau$. We apply this to a countable set $\mathcal D\subset C_c^1(\R)$ such that every $\varphi\in C_c(\R)$ is the uniform limit of a sequence in $\mathcal D$ with supports in a fixed compact set. By dominated convergence, there is a null set $N_0\subset(s,T_0)$ such that, for $\tau\notin N_0$, $\int_I\varphi(Y(\ell,\tau))\,e(\ell,\tau)\,d\ell=0$ for every $\varphi\in C_c(\R)$. Enlarging $N_0$ by a null set, we may also assume, by Fubini's theorem, that for $\tau\notin N_0$ $Y(\ell,\cdot)$ is differentiable at $\tau$ for almost every $\ell\in I$.

Fix $\tau\in(s,T_0)\setminus N_0$. For an atom $x$ of $\mu_\tau$, the level set $L_x:=\{\ell\in I:Y(\ell,\tau)=x\}$ is an interval of length $\mu_\tau(\{x\})$ by Lemma~\ref{mass-labels}(ii). Let $\ell<\ell'$ be points of $L_x$ such that $Y(\ell,\cdot)$ and $Y(\ell',\cdot)$ are differentiable at $\tau$. For small $h>0$ we have $Y(\ell,\tau\pm h)\le Y(\ell',\tau\pm h)$, and since $Y(\ell,\tau)=Y(\ell',\tau)$, dividing by $h$ and by $-h$ and letting $h\to0^+$ gives $\partial_\tau Y(\ell,\tau)\le\partial_\tau Y(\ell',\tau)$ and $\partial_\tau Y(\ell,\tau)\ge\partial_\tau Y(\ell',\tau)$. Hence $\partial_\tau Y(\cdot,\tau)$ equals a constant $c_x$ almost everywhere on $L_x$. If $\ell$ lies in no such interval, then $x=Y(\ell,\tau)$ is not an atom, and Lemma~\ref{mass-labels}(ii) gives $\ell=m(x,\tau)$, because $\ell\ge m(x,\tau)$, $\ell<m(x',\tau)$ for every $x'>x$, and $m(\cdot,\tau)$ is continuous at $x$.
Define $E(x):=c_x-w(x,\tau)$ if $x$ is an atom of $\mu_\tau$, $E(x):=e(m(x,\tau),\tau)$ if $x$ is not an atom and $m(x,\tau)\in I$, and $E(x):=0$ otherwise. Then $E$ is bounded and Borel, since the atoms are countable, and $e(\ell,\tau)=E(Y(\ell,\tau))$ for almost every $\ell\in I$. Consequently $\int\varphi E\,d\mu_\tau=\int_I\varphi(Y(\ell,\tau))\,e(\ell,\tau)\,d\ell=0$ for every $\varphi\in C_c(\R)$, so $E=0$ $\mu_\tau$-almost everywhere and $e(\cdot,\tau)=0$ almost everywhere on $I$. By Fubini's theorem, $e(\ell,\cdot)=0$ almost everywhere on $(s,T_0)$ for almost every $\ell\in I$.
\end{proof}

We call $\ell\in I$ a \emph{good mass coordinate} of $(\mu_\tau)$ if $Y(\ell,\cdot)$ has the property of Lemma~\ref{label-velocity}. By Lemma~\ref{forward-uniqueness}, if $\ell\le\ell'$ are good mass coordinates and $Y(\ell,\tau)<Y(\ell',\tau)$ for some $\tau$, then $Y(\ell,\tau')<Y(\ell',\tau')$ for every $\tau'\in[s,\tau]$. In words, concentrations do not split forward in time.

From now on, $m:=m_\rho$, and $I$ and $Y$ are the interval of mass coordinates and the position function of $(\rho_\tau)$. We assume $\rho_s\ne0$. Otherwise, the proof of Lemma~\ref{mass-labels}(i), which does not use $\mu_s\ne0$, shows that $\rho_\tau=0$ for every $\tau\in[s,T_0]$, and Proposition~\ref{positive-representation} is trivial. The momentum is transported in the same way.

\begin{lemma}\label{label-momentum}
Put $V(\ell):=w(Y(\ell,s),s)$ for $\ell\in I$. Then $P_\tau=Y(\cdot,\tau)_\#(V\,\Leb\llcorner I)$ for every $\tau\in[s,T_0]$, where, for a bounded Borel function $f$ on $I$, $Y(\cdot,\tau)_\#(f\,\Leb\llcorner I)$ denotes the locally finite signed measure that assigns to each bounded Borel set $B\subset\R$ the value $\int_{\{\ell\in I:\,Y(\ell,\tau)\in B\}}f\,d\ell$.
\end{lemma}

\begin{proof}
Put $\nu_\tau:=P_\tau+(M+1)\rho_\tau=(w(\cdot,\tau)+M+1)\rho_\tau$. Then $\nu_\tau\ge\rho_\tau\ge0$, $\tau\mapsto\nu_\tau$ is weakly-$*$ continuous, and \eqref{closed-identity} holds for $\nu$ by Lemma~\ref{closed-intervals} and linearity. Thus $(\nu_\tau)$ is a transport family, and we denote by $I'$ and $\widetilde Y$ its interval of mass coordinates and its position function.
Fix $\ell_0\in I$ and put $\Lambda(\ell):=\int_{\ell_0}^\ell(V+M+1)\,d\ell'$ for $\ell\in I$. Since $1\le V+M+1\le2M+1$, the map $\Lambda$ is increasing and bi-Lipschitz from $I$ onto an open interval $J$, and $\Lambda_\#\big((V+M+1)\Leb\llcorner I\big)=\Leb\llcorner J$. Hence $Z:=Y(\Lambda^{-1}(\cdot),s)$ is nondecreasing and right-continuous on $J$, and
$$
\begin{aligned}
Z_\#(\Leb\llcorner J)&=Y(\cdot,s)_\#\big((V+M+1)\Leb\llcorner I\big)\\
&=(w(\cdot,s)+M+1)\,Y(\cdot,s)_\#(\Leb\llcorner I)=(w(\cdot,s)+M+1)\rho_s=\nu_s,
\end{aligned}
$$
where we used $V=w(Y(\cdot,s),s)$. By Lemma~\ref{mass-labels}(iv) applied to $(\nu_\tau)$ at $\tau=s$, after adding a constant to $\Lambda$ we may assume $J=I'$ and $\widetilde Y(\Lambda(\ell),s)=Y(\ell,s)$ for every $\ell\in I$.
Let $N_\rho\subset I$ and $N_\nu\subset I'$ be Lebesgue null sets outside which the mass coordinates are good for $(\rho_\tau)$ and for $(\nu_\tau)$, respectively. Since $\Lambda^{-1}$ is Lipschitz, the set $N:=N_\rho\cup\Lambda^{-1}(N_\nu)$ is Lebesgue null. For every $\ell\in I\setminus N$, the curves $Y(\ell,\cdot)$ and $\widetilde Y(\Lambda(\ell),\cdot)$ are absolutely continuous, solve $\gamma'=w(\gamma,\tau)$ for almost every $\tau\in(s,T_0)$, and coincide at $\tau=s$. By Lemma~\ref{forward-uniqueness}, $\widetilde Y(\Lambda(\ell),\tau)=Y(\ell,\tau)$ for every $\tau\in[s,T_0]$ and every $\ell\in I\setminus N$. Hence, for every $\tau$,
$$
\nu_\tau=\widetilde Y(\cdot,\tau)_\#(\Leb\llcorner I')=\widetilde Y(\Lambda(\cdot),\tau)_\#\big((V+M+1)\Leb\llcorner I\big)=Y(\cdot,\tau)_\#\big((V+M+1)\Leb\llcorner I\big),
$$
and subtracting $(M+1)\rho_\tau=Y(\cdot,\tau)_\#\big((M+1)\Leb\llcorner I\big)$ gives the claim.
\end{proof}

\subsection{Representation at positive times}
\label{representation-proof}

We now use the mass coordinates to prove Proposition~\ref{positive-representation}. A concentration at time $\tau$ is a level interval of $Y(\cdot,\tau)$. Since concentrations do not split forward in time, a mass coordinate that separates concentrations at time $\tau$ also separates them at every earlier time. We formalize this through the following notion. For $\tau\in[s,T_0]$, a real number $\lambda$ with $\inf I\le\lambda\le\sup I$ is a \emph{$\tau$-cut} if there is a set $B\subset\R$ with $x'\in B$ whenever $x'<x\in B$, such that $\{\ell\in I:\ell<\lambda\}$ and $\{\ell\in I:Y(\ell,\tau)\in B\}$ differ by a Lebesgue null set. By Lemma~\ref{mass-labels}(ii), $m(x,\tau)$ is a $\tau$-cut for every $x\in\R$, with $B=(-\infty,x)$. We put
$$
\Theta_\tau(\ell):=Y(\ell,s)+(\tau-s)V(\ell),\qquad\ell\in I,
$$
which is the position at time $\tau$ of the portion of mass with coordinate $\ell$, moving freely with its velocity at time $s$.

\begin{lemma}\label{cuts}
Let $\tau\in[s,T_0]$.
\begin{enumerate}[label=(\roman*)]
\item Every $\tau$-cut is a $\tau'$-cut for every $\tau'\in[s,\tau]$.
\item If $\lambda_1$ and $\lambda_2$ are $\tau$-cuts, then $\int_{\lambda_1}^{\lambda_2}\big(Y(\ell,\tau)-\Theta_\tau(\ell)\big)\,d\ell=0$.
\end{enumerate}
\end{lemma}

\begin{proof}
(i) Let $\lambda$ be a $\tau$-cut with set $B$, and let $N\subset I$ be a null set that contains all mass coordinates that are not good and all those at which the two sets in the definition differ. For $\ell<\lambda\le\ell'$ in $I\setminus N$ we have $Y(\ell,\tau)\in B$ and $Y(\ell',\tau)\notin B$, hence $Y(\ell,\tau)<Y(\ell',\tau)$, and therefore $Y(\ell,\tau')<Y(\ell',\tau')$. Let $B'$ be the set of $x'$ such that $x'\le Y(\ell,\tau')$ for some $\ell\in I\setminus N$ with $\ell<\lambda$. Then $B'$ has the required monotonicity, it contains $Y(\ell,\tau')$ for every $\ell\in I\setminus N$ with $\ell<\lambda$, and it does not contain $Y(\ell',\tau')$ for $\ell'\in I\setminus N$ with $\ell'\ge\lambda$. Thus $\lambda$ is a $\tau'$-cut with set $B'$.

(ii) We may assume $\lambda_1<\lambda_2$. For $\tau'\in[s,\tau]$, let $B_1(\tau')$ and $B_2(\tau')$ be sets for which $\lambda_1$ and $\lambda_2$ are $\tau'$-cuts, which exist by (i). Then $[\lambda_1,\lambda_2)\cap I$ and $\{\ell\in I:Y(\ell,\tau')\in B_2(\tau')\setminus B_1(\tau')\}$ differ by a null set. The integral $\int_{\lambda_1}^{\lambda_2}Y(\ell,\tau')\,d\ell$ need not be finite when $\lambda_1$ or $\lambda_2$ is an endpoint of $I$, since a measure of finite mass may have an infinite first moment. We therefore work with displacements. By Lemma~\ref{mass-labels}(iii), $|Y(\ell,\tau')-Y(\ell,s)|\le M(\tau'-s)$, so
$$
f(\tau'):=\int_{\lambda_1}^{\lambda_2}\big(Y(\ell,\tau')-Y(\ell,s)\big)\,d\ell
$$
is finite and Lipschitz on $[s,\tau]$. Since $Y(\ell,\tau')-Y(\ell,s)=\int_s^{\tau'}\partial_\tau Y(\ell,r)\,dr$ for every $\ell$ and $|\partial_\tau Y|\le M$, Fubini's theorem gives $f(\tau')=\int_s^{\tau'}\int_{\lambda_1}^{\lambda_2}\partial_\tau Y(\ell,r)\,d\ell\,dr$. Hence, for almost every $\tau'$, by Lemma~\ref{label-velocity}, by $P_{\tau'}=w(\cdot,\tau')\rho_{\tau'}$, and by Lemma~\ref{label-momentum},
$$
f'(\tau')=\int_{\lambda_1}^{\lambda_2}w(Y(\ell,\tau'),\tau')\,d\ell=P_{\tau'}\big(B_2(\tau')\setminus B_1(\tau')\big)=\int_{\lambda_1}^{\lambda_2}V(\ell)\,d\ell.
$$
Here $\rho_{\tau'}(B_2(\tau')\setminus B_1(\tau'))=\lambda_2-\lambda_1$ is finite, and $P_{\tau'}(B_2(\tau')\setminus B_1(\tau'))$ stands for $\int_{B_2(\tau')\setminus B_1(\tau')}w(\cdot,\tau')\,d\rho_{\tau'}$. The second equality follows from $Y(\cdot,\tau')_\#(\Leb\llcorner I)=\rho_{\tau'}$, and the last one from Lemma~\ref{label-momentum}, both applied to the bounded sets $(B_2(\tau')\setminus B_1(\tau'))\cap[-n,n]$, together with the dominated convergence theorem as $n\to\infty$, since $|w|\le M$ and $|V|\le M$. Integrating from $s$ to $\tau$ and using $f(s)=0$ gives $f(\tau)=(\tau-s)\int_{\lambda_1}^{\lambda_2}V\,d\ell$, which is the claim. Note that the integrand in (ii) is bounded by $2M(\tau-s)$.
\end{proof}

Statement (ii) says that the barycenter of every group of mass coordinates bounded by cuts moves as if the corresponding particles moved freely. This is one half of the variational characterization of the sticky dynamics. The other half, an inequality, is obtained below from the momentum equation along straight lines.

We recall that $m=m_\rho$ is the cumulative mass function of Lemma~\ref{cumulative-function}, and that $I$ and $Y$ are the interval of mass coordinates and the position function of $(\rho_\tau)$ given by Lemma~\ref{mass-labels}. Put $m_*:=m(\bar a,s)$ and define, for $(x,\tau)\in\R\times[s,T_0]$,
$$
Q(x,\tau):=\int_{m_*}^{m(x,\tau)}V\,d\ell,
\qquad
\widetilde{\mathcal V}_s(x,\tau):=\int_{m_*}^{m(x,\tau)}\big(\Theta_\tau(\ell)-x\big)\,d\ell.
$$
These integrals are finite. Indeed, if $m_*\le\ell<m(x,\tau)$, then $Y(\ell,s)\ge \bar a$ and $Y(\ell,\tau)<x$ by Lemma~\ref{mass-labels}(ii), so $\bar a\le Y(\ell,s)<x+M(\tau-s)$ by Lemma~\ref{mass-labels}(iii). In the same way, $x-M(\tau-s)\le Y(\ell,s)<\bar a$ if $m(x,\tau)\le\ell<m_*$. Hence $\Theta_\tau$ is bounded on the range of integration, even when $\rho_s$ has an infinite first moment.
By Lemma~\ref{label-momentum} and Lemma~\ref{mass-labels}(ii), $Q(x_2,\tau)-Q(x_1,\tau)=P_\tau([x_1,x_2))$ for $x_1<x_2$, the function $Q(\cdot,\tau)$ is left-continuous, and $|Q|\le M|m-m_*|$. Lemma~\ref{mass-labels}(i), applied with $\tau_1=s$, gives $m(x-M(T_0-s),s)\le m(x,\tau)\le m(x+M(T_0-s),s)$ for $\tau\in[s,T_0]$, so $m$ and $Q$ are locally bounded on $\R\times[s,T_0]$. Moreover, $Q$ is Borel on $\R\times[s,T_0]$, since it is the composition of the Borel function $m$ with the continuous function $\lambda\mapsto\int_{m_*}^{\lambda}V\,d\ell$, and so is $Q-v\,m$ for every $v\in\R$.

\begin{lemma}\label{candidate-potential}
\begin{enumerate}[label=(\roman*)]
\item For $x,y\in\R$ and $\tau\in[s,T_0]$, $G_s(y;x,\tau)=\int_{m_*}^{m(y,s)}\big(\Theta_\tau(\ell)-x\big)\,d\ell$. In particular, $\widetilde{\mathcal V}_s(y,s)=G_s(y;y,s)$.
\item $\mathcal V_s\le\widetilde{\mathcal V}_s$ on $\R\times(s,T_0]$.
\item $\widetilde{\mathcal V}_s$ is locally Lipschitz on $\R\times[s,T_0]$. At every $(x,\tau)\in\R\times(s,T_0)$ with $\rho_\tau(\{x\})=0$, its partial derivatives exist and
$$
\partial_x\widetilde{\mathcal V}_s(x,\tau)=m_*-m(x,\tau),\qquad\partial_\tau\widetilde{\mathcal V}_s(x,\tau)=Q(x,\tau).
$$
\end{enumerate}
\end{lemma}

\begin{proof}
(i) For $y\ge \bar a$, Lemma~\ref{mass-labels}(ii) shows that $\{\ell\in I:Y(\ell,s)\in[\bar a,y)\}=[m_*,m(y,s))\cap I$. Since $Y(\cdot,s)_\#(\Leb\llcorner I)=\rho_s$, $P_s=w(\cdot,s)\rho_s$, and $(\tau-s)w(Y(\ell,s),s)+Y(\ell,s)-x=\Theta_\tau(\ell)-x$, the formula follows. The case $y<\bar a$ is the same. At $\tau=s$ we have $\Theta_s(\ell)=Y(\ell,s)$.

(ii) Let $\tau\in(s,T_0]$ and $x\in\R$. By Lemma~\ref{cuts}(i), $\lambda:=m(x,\tau)$ is an $s$-cut, with a set $B$ that is empty, equal to $\R$, or a half-line $(-\infty,y)$ or $(-\infty,y]$. If $B=(-\infty,y)$, then $\{\ell<\lambda\}$ and $\{\ell<m(y,s)\}$ differ by a null set by Lemma~\ref{mass-labels}(ii), so $\lambda=m(y,s)$ and $\widetilde{\mathcal V}_s(x,\tau)=G_s(y;x,\tau)\ge\mathcal V_s(x,\tau)$ by (i). If $B=(-\infty,y]$, then in the same way $\lambda=\lim\limits_{y'\to y^+}m(y',s)$, and $\widetilde{\mathcal V}_s(x,\tau)=\lim\limits_{y'\to y^+}G_s(y';x,\tau)\ge\mathcal V_s(x,\tau)$. If $B$ is empty, then $\lambda=m(x,\tau)=\inf I$, and Lemma~\ref{mass-labels}(i), applied with $\tau_1=s$ and $\tau_2=\tau$, gives $m(y,s)\le m(x,\tau)=\inf I$ for $y:=x-M(\tau-s)$. Hence $m(y,s)=\lambda$, and $\widetilde{\mathcal V}_s(x,\tau)=G_s(y;x,\tau)\ge\mathcal V_s(x,\tau)$ by (i). If $B=\R$, then $\lambda=m(x,\tau)=\sup I$, and the other inequality of Lemma~\ref{mass-labels}(i) gives $m(y,s)\ge m(x,\tau)=\sup I$ for $y:=x+M(\tau-s)$. Hence $m(y,s)=\lambda$, and the conclusion follows as before.

(iii) Let $x_1<x_2$ and $\lambda_i:=m(x_i,\tau)$. Then
$$
\widetilde{\mathcal V}_s(x_2,\tau)-\widetilde{\mathcal V}_s(x_1,\tau)=-(x_2-x_1)(\lambda_1-m_*)+\int_{\lambda_1}^{\lambda_2}\big(\Theta_\tau-x_2\big)\,d\ell.
$$
For $\ell\in[\lambda_1,\lambda_2)$ we have $x_1\le Y(\ell,\tau)<x_2$ by Lemma~\ref{mass-labels}(ii), and Lemma~\ref{cuts}(ii) gives $\int_{\lambda_1}^{\lambda_2}(\Theta_\tau-x_2)\,d\ell=\int_{\lambda_1}^{\lambda_2}(Y(\ell,\tau)-x_2)\,d\ell$, where both integrands are bounded. Hence
\begin{equation}\label{x-increment}
\big|\widetilde{\mathcal V}_s(x_2,\tau)-\widetilde{\mathcal V}_s(x_1,\tau)+(x_2-x_1)(\lambda_1-m_*)\big|\le(x_2-x_1)(\lambda_2-\lambda_1).
\end{equation}
Next let $x\in\R$, $s\le\tau_1<\tau_2\le T_0$, and $\lambda_i:=m(x,\tau_i)$. Then
$$
\widetilde{\mathcal V}_s(x,\tau_2)-\widetilde{\mathcal V}_s(x,\tau_1)=(\tau_2-\tau_1)Q(x,\tau_2)+\int_{\lambda_1}^{\lambda_2}\big(\Theta_{\tau_1}-x\big)\,d\ell.
$$
The number $\lambda_2$ is a $\tau_2$-cut, hence a $\tau_1$-cut by Lemma~\ref{cuts}(i), and $\lambda_1$ is a $\tau_1$-cut. By Lemma~\ref{cuts}(ii), the last integral equals $\int_{\lambda_1}^{\lambda_2}(Y(\ell,\tau_1)-x)\,d\ell$. If $\lambda_1<\lambda_2$, then for $\ell\in[\lambda_1,\lambda_2)\cap I$ we have $Y(\ell,\tau_1)\ge x>Y(\ell,\tau_2)$, so $0\le Y(\ell,\tau_1)-x\le M(\tau_2-\tau_1)$ by Lemma~\ref{mass-labels}(iii). If $\lambda_2<\lambda_1$, then for $\ell\in[\lambda_2,\lambda_1)\cap I$ we have $Y(\ell,\tau_1)<x\le Y(\ell,\tau_2)$, so $-M(\tau_2-\tau_1)\le Y(\ell,\tau_1)-x<0$. In both cases the last integral has absolute value at most $M(\tau_2-\tau_1)|\lambda_2-\lambda_1|$. Hence
\begin{equation}\label{t-increment}
\big|\widetilde{\mathcal V}_s(x,\tau_2)-\widetilde{\mathcal V}_s(x,\tau_1)-(\tau_2-\tau_1)Q(x,\tau_2)\big|\le M(\tau_2-\tau_1)\,|\lambda_2-\lambda_1|.
\end{equation}
Since $m$ and $Q$ are locally bounded on $\R\times[s,T_0]$, \eqref{x-increment} and \eqref{t-increment} show that $\widetilde{\mathcal V}_s$ is Lipschitz in each variable, with constants that are uniform on compact sets, and hence locally Lipschitz.
Now let $\rho_\tau(\{x\})=0$. Then $m(\cdot,\tau)$ is continuous at $x$, and \eqref{x-increment} gives $\partial_x\widetilde{\mathcal V}_s(x,\tau)=m_*-m(x,\tau)$. By Lemma~\ref{mass-labels}(i), applied on $[\tau,\tau+h]$ and on $[\tau-h,\tau]$,
$$
m(x-Mh,\tau)\le m(x,\tau\pm h)\le m(x+Mh,\tau),
$$
so $m(x,\cdot)$ and $Q(x,\cdot)$ are continuous at $\tau$. Then \eqref{t-increment}, applied with $\tau_1=\tau$ or $\tau_2=\tau$, gives $\partial_\tau\widetilde{\mathcal V}_s(x,\tau)=Q(x,\tau)$.
\end{proof}

The inequality along lines uses the following consequence of the momentum equation.

\begin{lemma}\label{line-monotonicity}
For $v\in\R$ put $Z_v:=Q-v\,m$. Then, for every nonnegative $\varphi\in C_c^1(\R)$ and $s\le\tau_1\le\tau_2\le T_0$,
$$
\int_\R\varphi(\xi)\,Z_v\big(\xi+v(\tau_2-s),\tau_2\big)\,d\xi\le\int_\R\varphi(\xi)\,Z_v\big(\xi+v(\tau_1-s),\tau_1\big)\,d\xi.
$$
\end{lemma}

\begin{proof}
We first derive from the description in mass coordinates the analogue of \eqref{cumulative-identity} for $Q$, in which the momentum equation \eqref{weak2-H} enters only through Lemma~\ref{label-momentum}. For $\Phi$ as in Lemma~\ref{closed-intervals},
\begin{equation}\label{Q-identity}
\Big[\int_\R\Phi(\cdot,\tau)\,Q(\cdot,\tau)\,dx\Big]_{\tau_1}^{\tau_2}=\int_{\tau_1}^{\tau_2}\Big(\int_\R\partial_\tau\Phi\,Q\,dx-\int_\R\Phi\,w\,dP_\tau\Big)\,d\tau.
\end{equation}
By Lemma~\ref{mass-labels}(ii), $Q(x,\tau)=\int_I\big(\mathbf 1_{\{\ell\ge m_*,\ Y(\ell,\tau)<x\}}-\mathbf 1_{\{\ell<m_*,\ Y(\ell,\tau)\ge x\}}\big)V(\ell)\,d\ell$, so that
$$
\int_\R\Phi(\cdot,\tau)\,Q(\cdot,\tau)\,dx=\int_I k_\Phi(\ell,Y(\ell,\tau),\tau)\,V(\ell)\,d\ell,
$$
where
$$
k_\Phi(\ell,z,\tau):=\mathbf 1_{\{\ell\ge m_*\}}\int_{(z,\infty)}\Phi(x,\tau)\,dx-\mathbf 1_{\{\ell<m_*\}}\int_{(-\infty,z]}\Phi(x,\tau)\,dx.
$$
If $\Phi$ vanishes outside $[c_1,c_2]\times[\tau_1,\tau_2]$, then, by Lemma~\ref{mass-labels}(ii) and (iii), $k_\Phi(\ell,Y(\ell,\tau),\tau)=0$ for all $\tau$ unless $\ell$ lies in the bounded interval between $\min\{m_*,m(c_1-M(T_0-s),s)\}$ and $\max\{m_*,m(c_2+M(T_0-s),s)\}$. For a good mass coordinate $\ell$, the function $\tau\mapsto k_\Phi(\ell,Y(\ell,\tau),\tau)$ is Lipschitz with derivative $k_{\partial_\tau\Phi}(\ell,Y(\ell,\tau),\tau)-\Phi(Y(\ell,\tau),\tau)\,w(Y(\ell,\tau),\tau)$ for almost every $\tau$. Integrating in $\tau$ and $\ell$, using the representation above with $\partial_\tau\Phi$ in place of $\Phi$ and Lemma~\ref{label-momentum} in the form $\int_I\Phi(Y(\ell,\tau),\tau)\,w(Y(\ell,\tau),\tau)\,V(\ell)\,d\ell=\int\Phi(\cdot,\tau)\,w(\cdot,\tau)\,dP_\tau$, we obtain \eqref{Q-identity}.

Now put $\Phi(x,\tau):=\varphi(x-v(\tau-s))$, so that $\partial_\tau\Phi=-v\partial_x\Phi$. Since $Z_v(\cdot,\tau)$ is left-continuous with Stieltjes measure $P_\tau-v\rho_\tau$, we have $\int\partial_\tau\Phi Z_v\,dx=v\int\Phi\,d(P_\tau-v\rho_\tau)$. By \eqref{cumulative-identity} for $m$ and \eqref{Q-identity},
$$
\Big[\int_\R\Phi\,Z_v\,dx\Big]_{\tau_1}^{\tau_2}=\int_{\tau_1}^{\tau_2}\!\!\int_\R\Phi\,\big(2vw-w^2-v^2\big)\,d\rho_\tau\,d\tau=-\int_{\tau_1}^{\tau_2}\!\!\int_\R\Phi\,(w-v)^2\,d\rho_\tau\,d\tau\le0,
$$
where we used $dP_\tau=w\,d\rho_\tau$. The change of variables $x=\xi+v(\tau-s)$ gives the claim.
\end{proof}

\begin{proposition}\label{line-inequality}
For all $(x_0,\tau_0)\in\R\times(s,T_0]$ and $y_0\in\R$, $\widetilde{\mathcal V}_s(x_0,\tau_0)\le G_s(y_0;x_0,\tau_0)$.
\end{proposition}

\begin{proof}
Put $v:=(x_0-y_0)/(\tau_0-s)$. The set of $(x,\tau)$ with $\rho_\tau(\{x\})>0$ is Borel, since $\rho_\tau(\{x\})=\lim\limits_{k\to\infty}\big(m(x+1/k,\tau)-m(x,\tau)\big)$, and each of its sections $\{x:\rho_\tau(\{x\})>0\}$ is countable, so it is a Lebesgue null set. By Rademacher's theorem and Lemma~\ref{candidate-potential}(iii), $\widetilde{\mathcal V}_s$ is differentiable at almost every point of $\R\times(s,T_0)$, with gradient $(m_*-m,Q)$. The linear map $(\xi,\tau)\mapsto(\xi+v(\tau-s),\tau)$ preserves Lebesgue null sets, so, by Fubini's theorem, for almost every $\xi$ the Lipschitz function $\tau\mapsto\widetilde{\mathcal V}_s(\xi+v(\tau-s),\tau)$ has derivative $Q-vm+vm_*=Z_v+vm_*$, evaluated at $(\xi+v(\tau-s),\tau)$, for almost every $\tau\in(s,\tau_0)$. Hence, for almost every $\xi$,
$$
\widetilde{\mathcal V}_s\big(\xi+v(\tau_0-s),\tau_0\big)-\widetilde{\mathcal V}_s(\xi,s)=\int_s^{\tau_0}\Big(Z_v\big(\xi+v(\tau-s),\tau\big)+vm_*\Big)\,d\tau.
$$
Let $\varphi\in C_c^1(\R)$ be nonnegative, with support in $(y_0-\delta,y_0)$ and $\int\varphi=1$. Multiplying by $\varphi(\xi)$, integrating, and using Fubini's theorem and Lemma~\ref{line-monotonicity} with $\tau_1=s$, we obtain
$$
\int_\R\varphi(\xi)\Big(\widetilde{\mathcal V}_s\big(\xi+v(\tau_0-s),\tau_0\big)-\widetilde{\mathcal V}_s(\xi,s)\Big)\,d\xi\le(\tau_0-s)\int_\R\varphi(\xi)\big(Z_v(\xi,s)+vm_*\big)\,d\xi.
$$
As $\delta\to0^+$, the left-hand side converges to $\widetilde{\mathcal V}_s(x_0,\tau_0)-\widetilde{\mathcal V}_s(y_0,s)$ by continuity of $\widetilde{\mathcal V}_s$. Since $m(\cdot,s)$ is left-continuous and $Q(\xi,s)=\int_{m_*}^{m(\xi,s)}V\,d\ell$, we have $Z_v(\xi,s)\to Z_v(y_0,s)$ as $\xi\to y_0^-$, and the right-hand side converges to $(\tau_0-s)\big(Z_v(y_0,s)+vm_*\big)$. Therefore
$$
\widetilde{\mathcal V}_s(x_0,\tau_0)\le\widetilde{\mathcal V}_s(y_0,s)+(\tau_0-s)\,Q(y_0,s)-(x_0-y_0)\big(m(y_0,s)-m_*\big).
$$
By Lemma~\ref{candidate-potential}(i), the right-hand side equals
$$
\int_{m_*}^{m(y_0,s)}\Big(Y(\ell,s)-y_0+(\tau_0-s)V(\ell)-(x_0-y_0)\Big)\,d\ell
=\int_{m_*}^{m(y_0,s)}\big(\Theta_{\tau_0}(\ell)-x_0\big)\,d\ell,
$$
which is $G_s(y_0;x_0,\tau_0)$.
\end{proof}

\begin{proof}[Proof of Proposition~\ref{positive-representation}]
By Proposition~\ref{line-inequality}, $\widetilde{\mathcal V}_s\le\inf\limits_yG_s=\mathcal V_s$ on $\R\times(s,T_0]$, and by Lemma~\ref{candidate-potential}(ii), $\mathcal V_s\le\widetilde{\mathcal V}_s$. Hence $\mathcal V_s=\widetilde{\mathcal V}_s$, and by Lemma~\ref{candidate-potential}(iii), for every $\tau\in(s,T_0)$ the Lipschitz function $\mathcal V_s(\cdot,\tau)$ has derivative $m_*-m(\cdot,\tau)$ almost everywhere. Let $\varphi\in C_c^\infty(\R\times(s,T_0))$. Integrating by parts in $x$ and using $\int\partial_x\varphi\,m(\cdot,\tau)\,dx=-\int\varphi\,d\rho_\tau$,
$$
\iint\mathcal V_s\,\partial_x^2\varphi\,dx\,d\tau=-\iint(m_*-m)\,\partial_x\varphi\,dx\,d\tau=-\int\!\!\int\varphi\,d\rho_\tau\,d\tau,
$$
so $\partial_x^2\mathcal V_s=-\rho$. In the same way,
$$
\iint\mathcal V_s\,\partial_x\partial_\tau\varphi\,dx\,d\tau=-\iint(m_*-m)\,\partial_\tau\varphi\,dx\,d\tau=\iint m\,\partial_\tau\varphi\,dx\,d\tau=\int\!\!\int\varphi\,dP_\tau\,d\tau,
$$
where the last equality is \eqref{cumulative-identity} on an interval $[\tau_1,\tau_2]$ containing the support of $\varphi$ in time, whose boundary terms vanish. Hence $\partial_x\partial_\tau\mathcal V_s=P$.
\end{proof}

The proof shows more than Proposition~\ref{positive-representation}. In mass coordinates, the infimum defining $\mathcal V_s(x,\tau)$ is realized by the cut $m(x,\tau)$, and in physical space it equals $G_s(y;x,\tau)$ or $\lim\limits_{y'\to y^+}G_s(y';x,\tau)$ for some $y\in\R$, by the proof of Lemma~\ref{candidate-potential}(ii). The position map $Y(\cdot,\tau)$ is determined by the data at time $s$. This is the variational characterization of the sticky dynamics, with the barycentric identity of Lemma~\ref{cuts} as the equality and Proposition~\ref{line-inequality} as the inequality.

\section{Generalized potential and characteristics}
\label{characteristics}

The theory of Huang and Wang \cite{HuangWang2001} uses two geometric constructions, namely a generalized potential, whose minimizers determine the entropy solution, and generalized characteristics, which describe where the mass initially located at a point $\eta$ is found at later times. The formulas used below transfer through $\mathcal T$ by the substitution $\tau=A(t)$. The transfer of the potential is a direct computation. The transfer of the characteristics requires more care, because at a concentration point an admissible velocity representative may have different one-sided limits, and the characteristic equation is then a differential inclusion between them.

\subsection{Generalized potential}
\label{generalized-potential}

If $\rho_0=0$, choose $m^\sharp=q^\sharp=0$ and omit the minimizer selections below. The zero density and zero momentum are the entropy solution by Theorem~\ref{uniqueness}. Throughout the rest of this subsection assume $\rho_0\ne0$.
Let $\bar a$ and $G$ be as in Section~\ref{strips}, so that $\rho_0(\{\bar a\})=0$, and define
$$
F(y;x,t):=G\big(y;x,A(t)\big)=\int_{\bar a}^y\big(A(t)u_0(\eta)+\eta-x\big)\,d\rho_0(\eta),
$$
which is the potential of Theorem~\ref{main-introduction}(iii). The potential of \cite[(2.1) and (2.2)]{HuangWang2001} is written with base point $0$ and an endpoint convention at $0$. It differs from $G$ by a quantity that does not depend on $y$, so both potentials have the same minimizing sequences, and all the objects defined below from their minimizers coincide. We set
$$
\nu(x,t):=\inf_{y\in\R}F(y;x,t),
\qquad
\Scal(x,t):=\big\{y\in\R:\ \exists\, y_n\to y\ \text{with}\ F(y_n;x,t)\to\nu(x,t)\big\}.
$$
Since $F(\cdot;x,t)$ is only left-continuous, the infimum need not be attained, and $\Scal(x,t)$ is the set of limits of minimizing sequences. From $\Scal(x,t)$ and $\rho_0$, two distinguished points $y_*(x,t)\le y^*(x,t)$ of $\Scal(x,t)\cap\supp\rho_0$ are selected in \cite[(2.5)--(2.7)]{HuangWang2001}. They determine the backward generalized characteristics. The integrand $A(t)u_0(\eta)+\eta-x$ is the signed distance at time $t$ between the point $x$ and the free particle issued from $\eta$ with velocity $u_0(\eta)$.

\begin{proposition}\label{potential-transfer}
Let $\Scal_H(x,\tau)$, $y_*^H(x,\tau)$, and $y^{*H}(x,\tau)$ denote the objects defined in the same way from $G(\cdot;x,\tau)$. Then, for every $(x,t)$,
$$
\begin{gathered}
\nu(x,t)=\mathcal V(x,A(t)),\qquad
\Scal(x,t)=\Scal_H(x,A(t)),\\
y_*(x,t)=y_*^H(x,A(t)),\qquad
y^*(x,t)=y^{*H}(x,A(t)),
\end{gathered}
$$
and $F(\cdot;x,t)$ attains its infimum at a point $y$ if and only if $G(\cdot;x,A(t))$ does.
\end{proposition}

\begin{proof}
By definition, $F(\cdot;x,t)$ and $G(\cdot;x,A(t))$ are the same function of $y$, and $\mathcal V$ is the infimum of $G$. All the other objects are defined from this function and from $\rho_0$ alone, without further reference to the time variable.
\end{proof}

The value of the infimum determines the solution, without any selection among the minimizers.

\begin{corollary}\label{potential-solution}
Let Assumption~\ref{coefficients} hold and let $(\rho,u)$ be the entropy solution of \eqref{NA}. Then, for every $t>0$, $\nu(\cdot,t)$ is concave and $-\partial_x^2\nu(\cdot,t)=\rho_t$ in $\mathcal D'(\R)$. Moreover, $\partial_x\partial_t\nu=\alpha(t)\,u\,\rho_t\,dt$ in $\mathcal D'(\R\times(0,\infty))$.
\end{corollary}

\begin{proof}
Let $(\tilde\rho,w):=\mathcal T(\rho,u)$ and $P_\tau:=w(\cdot,\tau)\tilde\rho_\tau$. By Proposition~\ref{potential-transfer}, $\nu(x,t)=\mathcal V(x,A(t))$. In the proof of Theorem~\ref{strip-uniqueness} we obtained $\tilde\rho=-\partial_x^2\mathcal V$ and $P=\partial_x\partial_\tau\mathcal V$ in $\mathcal D'(\R\times(0,T_0))$ for every $T_0<A_\infty$. By Lemma~\ref{localization}, $\mathcal V$ is concave on $\R\times(0,T_0)$, hence continuous and locally Lipschitz. For $\varphi\in C_c^\infty(\R)$, the functions $\tau\mapsto\int_\R\mathcal V(\cdot,\tau)\varphi''\,dx$ and $\tau\mapsto\int_\R\varphi\,d\tilde\rho_\tau$ are continuous on $(0,A_\infty)$, the second by condition (a) of Definition~\ref{weak-solution}, and the first identity gives $\int_\R\mathcal V(\cdot,\tau)\varphi''\,dx=-\int_\R\varphi\,d\tilde\rho_\tau$ for almost every $\tau$, hence for every $\tau\in(0,A_\infty)$. At $\tau=A(t)$ this is $-\partial_x^2\nu(\cdot,t)=\rho_t$, since $\tilde\rho_{A(t)}=\rho_t$, and $\nu(\cdot,t)=\mathcal V(\cdot,A(t))$ is concave.

For the second identity, let $\Xi\in C_c^\infty(\R\times(0,\infty))$ and $\Psi(x,\tau):=\Xi(x,A^{-1}(\tau))$. By Lemma~\ref{testfunctions}, $\Psi\in C_c^1(\R\times(0,A_\infty))$, the derivatives of $\Psi$ and $\partial_\tau\Psi$ in $x$ of every order are continuous, and $\partial_t\Xi(x,t)=A'(t)\partial_\tau\Psi(x,A(t))$. The identity $P=\partial_x\partial_\tau\mathcal V$ extends by mollification to test functions with these properties, because $\mathcal V$ is locally bounded and $P$ is a locally finite measure. The substitution $\tau=A(t)$ gives
$$
\iint\nu\,\partial_x\partial_t\Xi\,dx\,dt=\iint\mathcal V\,\partial_x\partial_\tau\Psi\,dx\,d\tau=\int_0^{A_\infty}\!\!\int_\R\Psi\,dP_\tau\,d\tau=\int_0^\infty\!\!\int_\R\Xi\,A'(t)e^{S(t)}u\,d\rho_t\,dt,
$$
where the last equality uses $\tilde\rho_{A(t)}=\rho_t$ and $w(\cdot,A(t))=e^{S(t)}u(\cdot,t)$ $\rho_t$-almost everywhere. Since $A'(t)e^{S(t)}=\alpha(t)$, the second identity follows.
\end{proof}

The representation of the solution through the minimizers transfers as well.
Following \cite[(2.35) and (2.36)]{HuangWang2001}, define
$$
m^\sharp(x,t):=
\begin{cases}
\displaystyle\int_{\bar a}^{y_*(x,t)}d\rho_0, &\text{if } F(\cdot;x,t)\text{ attains its infimum at }y_*(x,t),\\[3mm]
\displaystyle\int_{\bar a}^{y_*(x,t)}d\rho_0+\rho_0(\{y_*(x,t)\}), &\text{otherwise},
\end{cases}
$$
and let $q^\sharp(x,t)$ be defined by the same formula with $d\rho_0$ replaced by $e^{-S(t)}u_0\,d\rho_0$, and $\rho_0(\{y_*\})$ by $e^{-S(t)}u_0(y_*)\rho_0(\{y_*\})$.

\begin{proposition}[Representation formulas]\label{representation}
Let $(\rho,u)$ be the entropy solution of Theorem~\ref{existence}, with the zero-density convention above. Then $m^\sharp$ satisfies \eqref{potential-system}, so that it coincides almost everywhere, up to an additive constant, with the potential of Proposition~\ref{intrinsic-potential}, and
\begin{equation}\label{momentum-potential}
\partial_xq^\sharp(\cdot,t)=u(\cdot,t)\,\rho_t
\end{equation}
as Radon measures, for almost every $t>0$.
\end{proposition}

\begin{proof}
If $\rho_0=0$, both asserted identities hold for the zero primitives. Assume henceforth that $\rho_0\ne0$.
Let $\tilde m$ and $\tilde q$ be the mass and momentum primitives \cite[(2.35) and (2.36)]{HuangWang2001} of the solution of Theorem~\ref{HW}, so that $(\tilde m,w)$ is that solution. They are given by the formulas above with $G(\cdot;x,\tau)$, $y_*^H(x,\tau)$, and $u_0\,d\rho_0$ in place of $F(\cdot;x,t)$, $y_*(x,t)$, and $e^{-S(t)}u_0\,d\rho_0$, except that the integrals start at $0$ with the endpoint convention of \cite{HuangWang2001}. By Proposition~\ref{potential-transfer}, there are constants $c$ and $c'$ such that $m^\sharp(x,t)=\tilde m(x,A(t))+c$ and $q^\sharp(x,t)=e^{-S(t)}\big(\tilde q(x,A(t))+c'\big)$. By \cite[Lemma~2.8]{HuangWang2001}, $\partial_x\tilde q(\cdot,\tau)=w(\cdot,\tau)\,\partial_x\tilde m(\cdot,\tau)$ for almost every $\tau>0$.

First, $\partial_xm^\sharp(\cdot,t)=\partial_x\tilde m(\cdot,A(t))=\tilde\rho_{A(t)}=\rho_t$. Next, let $\Xi\in C_c^\infty(\R\times(0,\infty))$ and $\Phi(x,\tau):=\Xi(x,A^{-1}(\tau))$, which belongs to $C_c^1(\R\times(0,A_\infty))$ by Lemma~\ref{testfunctions}. By the substitution $\tau=A(t)$, $\partial_t\Xi(x,t)=A'(t)\partial_\tau\Phi(x,A(t))$, and $\iint c\,\partial_t\Xi\,dx\,dt=0$,
$$
\begin{aligned}
\iint_{\R\times(0,\infty)} m^\sharp\,\partial_t\Xi\,dx\,dt
&=\iint_{\R\times(0,A_\infty)} \tilde m\,\partial_\tau\Phi\,dx\,d\tau
=\iint_{\R\times(0,A_\infty)}\Phi\,w\,d\tilde\rho_\tau\,d\tau\\
&=\int_0^\infty\!\!\int_\R\Xi\,A'(t)e^{S(t)}u\,d\rho_t\,dt,
\end{aligned}
$$
where the second equality is the first identity in the proof of Theorem~\ref{existence}, which extends to $C_c^1$ test functions by mollification as at the end of the proof of Proposition~\ref{HW-equivalence}, and the third follows from the substitution $\tau=A(t)$ and $w(\cdot,A(t))=e^{S(t)}u(\cdot,t)$ $\rho_t$-almost everywhere. Since $A'e^S=\alpha$, this states that $\partial_tm^\sharp=-\alpha u\rho$, so $m^\sharp$ satisfies \eqref{potential-system}, and the uniqueness part of Proposition~\ref{intrinsic-potential} applies. Finally, for every $t$ such that this identity holds at $\tau=A(t)$, which is the case for almost every $t$ by Lemma~\ref{properties-A}(iii), $\partial_xq^\sharp(\cdot,t)=e^{-S(t)}\,w(\cdot,A(t))\,\tilde\rho_{A(t)}=u(\cdot,t)\,\rho_t$, which is \eqref{momentum-potential}.
\end{proof}

\subsection{Generalized characteristics}
\label{generalized-characteristics}

The minimizers of the potential trace out curves in the $(x,t)$ plane. At points where the admissible velocity representative is continuous, the differential inclusion below reduces to the corresponding characteristic equation. At a concentration point, an admissible velocity representative may have different one-sided limits, and the differential inclusion expresses the relation between these traces and the characteristic speed. The motion of an atom must also satisfy the weak equations, which impose additional compatibility conditions on its speed and weight. One-sided limits are available for admissible representatives. If $\hat w(\cdot,\tau)$ satisfies \eqref{homogeneous-oleinik}, then $x\mapsto\hat w(x,\tau)-x/\tau$ is nonincreasing, so the limits $\hat w(x\pm,\tau)$ exist at every $x$ and $\hat w(x+,\tau)\le\hat w(x-,\tau)$. Likewise, if $\hat u(\cdot,t)$ satisfies \eqref{intrinsic-oleinik}, then $\hat u(x\pm,t)$ exist and $\hat u(x+,t)\le\hat u(x-,t)$. Moreover, if $\hat w(\cdot,A(t))=e^{S(t)}\hat u(\cdot,t)$ as in Proposition~\ref{representative-correspondence}, then $\hat w(x\pm,A(t))=e^{S(t)}\hat u(x\pm,t)$.

Let $\hat w$ be a family of admissible representatives of $w(\cdot,\tau)$, defined for almost every $\tau\in(0,A_\infty)$, and let $X_H$ be a locally Lipschitz curve on $(0,A_\infty)$ such that, for almost every $\tau$,
\begin{equation}\label{char-inclusion-H}
X_H'(\tau)\in\big[\hat w(X_H(\tau)+,\tau),\,\hat w(X_H(\tau)-,\tau)\big].
\end{equation}
At a continuity point of $\hat w(\cdot,\tau)$, the two one-sided limits coincide and \eqref{char-inclusion-H} reduces to $X_H'(\tau)=\hat w(X_H(\tau),\tau)$.

\begin{theorem}\label{characteristics-transfer}
Let $X_H$ and $\hat w$ be as above, and let $\hat u(\cdot,t):=e^{-S(t)}\hat w(\cdot,A(t))$.
Then $X(t):=X_H(A(t))$ is locally Lipschitz on $(0,\infty)$ and, for almost every $t>0$,
\begin{equation}\label{char-inclusion-NA}
X'(t)\in\alpha(t)\big[\hat u(X(t)+,t),\,\hat u(X(t)-,t)\big].
\end{equation}
In particular, at every such $t$ at which $\hat u(\cdot,t)$ is continuous at $X(t)$, $X'(t)=\alpha(t)\,\hat u(X(t),t)$.
\end{theorem}

\begin{proof}
The composition of the locally Lipschitz function $X_H$ with the $C^1$ function $A$ is locally Lipschitz. Let $N\subset(0,A_\infty)$ be the null set of those $\tau$ at which $X_H$ is not differentiable, $\hat w(\cdot,\tau)$ is not defined, or \eqref{char-inclusion-H} fails. By Lemma~\ref{properties-A}(iii), $A^{-1}(N)$ is null. At every $t\notin A^{-1}(N)$, the chain rule gives $X'(t)=A'(t)X_H'(A(t))$, and \eqref{char-inclusion-H} at $\tau=A(t)$, multiplied by $A'(t)>0$, yields
$$
X'(t)\in A'(t)\big[\hat w(X(t)+,A(t)),\,\hat w(X(t)-,A(t))\big]
=A'(t)e^{S(t)}\big[\hat u(X(t)+,t),\,\hat u(X(t)-,t)\big].
$$
Since $A'(t)e^{S(t)}=\alpha(t)$, this is \eqref{char-inclusion-NA}. At a continuity point the interval reduces to the single point $\alpha(t)\hat u(X(t),t)$.
\end{proof}

Theorem~\ref{characteristics-transfer} transfers the differential inclusion for every homogeneous curve satisfying \eqref{char-inclusion-H}. Its endpoints involve $\alpha(t)\hat u$, so the spatial speed of a characteristic differs from the velocity by the factor $\alpha(t)$, a distinction that is absent when $\alpha\equiv1$, as in the homogeneous theory and in Jin \cite{Jin2016}.

The forward generalized characteristics constructed in \cite[(2.31)--(2.33)]{HuangWang2001} transfer in the same way. For $\eta\in\supp\rho_0$, let $x_H(\eta,\tau)$ denote the forward generalized characteristic issued from $(\eta,0)$ defined there. It is nondecreasing in $\eta$ and Lipschitz continuous in $\tau>0$.

\begin{proposition}[Forward characteristics]\label{minimizer-char}
For $\eta\in\supp\rho_0$ and $t>0$, put $X(\eta,t):=x_H(\eta,A(t))$. Then $X(\eta,\cdot)$ is locally Lipschitz on $(0,\infty)$, $X(\cdot,t)$ is nondecreasing, and, for every $x\in\R$ and $t>0$,
\begin{equation}\label{HW-duality}
X(\eta,t)=x\qquad\text{whenever}\quad
y_*(x-0,t)<\eta<y_*(x+0,t),\ \ \eta\in\supp\rho_0 .
\end{equation}
Moreover, let $w$ be the velocity of Theorem~\ref{HW}, let $u(x,t):=e^{-S(t)}w(x,A(t))$ be the velocity of the entropy solution constructed in the proof of Theorem~\ref{existence}, and put $M:=\|u_0\|_{L^\infty(d\rho_0)}$.
Then, for every $\eta\in\supp\rho_0$ and every $t>0$ at which $X(\eta,\cdot)$ is differentiable and $|X(\eta,t)-y_*(X(\eta,t),t)|\le MA(t)$,
\begin{equation}\label{forward-speed}
\partial_tX(\eta,t)=\alpha(t)\,u\big(X(\eta,t),t\big).
\end{equation}
\end{proposition}

\begin{proof}
Local Lipschitz continuity in $t$ follows from that of $x_H(\eta,\cdot)$ and Lemma~\ref{properties-A}, and monotonicity in $\eta$ is inherited from $x_H$. By \cite[(2.46)]{HuangWang2001}, $x_H(\eta,\tau)=x$ whenever $y_*^H(x-0,\tau)<\eta<y_*^H(x+0,\tau)$ and $\eta\in\supp\rho_0$. Taking $\tau=A(t)$ and using $y_*(\cdot,t)=y_*^H(\cdot,A(t))$ from Proposition~\ref{potential-transfer}, whose one-sided limits in $x$ therefore coincide as well, gives \eqref{HW-duality}.

Since $A$ is a $C^1$ bijection of $[0,\infty)$ onto $[0,A_\infty)$ with $A'>0$, the function $X(\eta,\cdot)$ is differentiable at $t$ if and only if $x_H(\eta,\cdot)$ is differentiable at $A(t)$.
By \cite[(2.34)]{HuangWang2001}, at every such $\tau=A(t)$ at which, in addition, $|x_H(\eta,\tau)-y_*^H(x_H(\eta,\tau),\tau)|\le M\tau$, we have $\partial_\tau x_H(\eta,\tau)=w(x_H(\eta,\tau),\tau)$.
Since $y_*(\cdot,t)=y_*^H(\cdot,A(t))$, this condition is the one in the statement, and the chain rule gives
$$
\partial_tX(\eta,t)=A'(t)\,w\big(X(\eta,t),A(t)\big)=A'(t)e^{S(t)}u\big(X(\eta,t),t\big)=\alpha(t)\,u\big(X(\eta,t),t\big).
$$
\end{proof}

Relation \eqref{HW-duality} expresses the formation of delta shocks in terms of the potential. When $y_*(\cdot,t)$ jumps at $x$, the forward characteristics issued from all points of $\supp\rho_0$ in the open interval $(y_*(x-0,t),y_*(x+0,t))$ reach the point $x$ at time $t$. By \eqref{forward-speed}, the curves $X(\eta,\cdot)$ move with speed $\alpha(t)u$ at the times at which this speed identity applies.

\section{Examples}
\label{examples}

By Theorem~\ref{cauchy-equivalence}, all the information on the coefficients is contained in the effective time $A$ and the factor $e^{-S}$.
We first compute $A$ for several representative choices of $\alpha$ and $\sigma$, recovering the autonomous theory and the constant-damping theory of Jin \cite{Jin2016} as special cases, and we then treat the Riemann problem, for which the formulas of De~la~cruz and Neves \cite{DelaCruzNeves2025} allow a direct comparison. Finally, we give bounded data for which the particles issued from almost every initial point belong to a concentration at every positive time.

\subsection{Explicit effective times}

The simplest case shows that the homogeneous system is the base point of the family \eqref{NA}. If $\alpha\equiv1$ and $\sigma\equiv0$, then $S\equiv0$ and $A(t)=t$, so the transformation is the identity, and \eqref{NA} coincides with the homogeneous system \eqref{H}.

The first nontrivial case is that of constant damping, $\alpha\equiv1$ and $\sigma\equiv\sigma_0>0$. Here $S(t)=\sigma_0 t$ and
$$
A(t)=\int_0^t e^{-\sigma_0 s}\,ds=\frac{1-e^{-\sigma_0 t}}{\sigma_0},
\qquad
A_\infty=\frac{1}{\sigma_0}.
$$
The whole evolution of \eqref{NA} on $[0,\infty)$ therefore corresponds to the evolution of \eqref{H} on the bounded interval $[0,1/\sigma_0)$. Thus friction compresses the effective clock, although physical time runs forever. For $\sigma_0=1$ this is the system studied by Jin \cite{Jin2016}, and Theorem~\ref{existence-uniqueness} recovers his existence and uniqueness result as a special case. The transformation also explains the similarity between his generalized potential \cite[(1.3)]{Jin2016} and that of the homogeneous theory. By Proposition~\ref{potential-transfer}, the former is obtained from the latter by replacing $\tau$ with $A(t)=1-e^{-t}$. Since the homogeneous solution does not depend on $\sigma_0$, the same representation yields the vanishing-damping limit directly.

\begin{proposition}[Vanishing damping]\label{vanishing-damping}
Let $\alpha\equiv1$, let $\rho_0\in\M_{\mathrm{loc}}(\R)$ be nonnegative and $u_0\in L^\infty(\R,d\rho_0)$, and, for $\sigma_0>0$, let $(\rho^{\sigma_0},u^{\sigma_0})$ be the entropy solution of \eqref{NA} with $\sigma\equiv\sigma_0$. Let $(\tilde\rho,w)$ be the entropy solution of \eqref{H} with the same data.
Then, for every $\varphi\in C_c(\R)$ and every $T>0$,
$$
\sup_{t\in[0,T]}\Big|\int_\R\varphi\,d\rho^{\sigma_0}_t-\int_\R\varphi\,d\tilde\rho_t\Big|
+\sup_{t\in[0,T]}\Big|\int_\R\varphi\,u^{\sigma_0}\,d\rho^{\sigma_0}_t-\int_\R\varphi\,w\,d\tilde\rho_t\Big|
\longrightarrow0
\qquad\text{as }\sigma_0\to0^+.
$$
\end{proposition}

\begin{proof}
Let $A_{\sigma_0}(t)=(1-e^{-\sigma_0 t})/\sigma_0$. By Theorems~\ref{cauchy-equivalence} and~\ref{uniqueness}, $\rho^{\sigma_0}_t=\tilde\rho_{A_{\sigma_0}(t)}$ and $u^{\sigma_0}(\cdot,t)\rho^{\sigma_0}_t=e^{-\sigma_0 t}\,(w\tilde\rho)_{A_{\sigma_0}(t)}$ as measures. Since $1-\sigma_0 s\le e^{-\sigma_0 s}\le1$, we have $t-\sigma_0 t^2/2\le A_{\sigma_0}(t)\le t$, so $A_{\sigma_0}(t)\to t$ uniformly on $[0,T]$. For $\varphi\in C_c(\R)$, the functions $\tau\mapsto\int_\R\varphi\,d\tilde\rho_\tau$ and $\tau\mapsto\int_\R\varphi\,w\,d\tilde\rho_\tau$ are continuous by condition (a) of Definition~\ref{weak-solution}, hence uniformly continuous and bounded on $[0,T]$. Therefore $\int_\R\varphi\,d\tilde\rho_{A_{\sigma_0}(t)}\to\int_\R\varphi\,d\tilde\rho_t$ uniformly on $[0,T]$, and, since also $e^{-\sigma_0 t}\to1$ uniformly on $[0,T]$, $e^{-\sigma_0 t}\int_\R\varphi\,w\,d\tilde\rho_{A_{\sigma_0}(t)}\to\int_\R\varphi\,w\,d\tilde\rho_t$ uniformly on $[0,T]$.
\end{proof}

The vanishing-damping limit was studied by Tang \cite[Theorem~1.1]{Tang2025}, with $\sigma_0=1/\tau_{\mathrm r}$ for a relaxation time $\tau_{\mathrm r}\to\infty$, by passing to the limit in the generalized potential and its minimizers. Through the transformation, the limit follows from the continuity in time of the homogeneous solution alone, and Proposition~\ref{vanishing-damping} gives it for the density and the momentum, locally uniformly in time. Between the undamped and the constant-damping cases lies time-gradually-degenerate damping, a class covered by the Riemann theory of De~la~cruz and Neves \cite{DelaCruzNeves2025}, for instance $\alpha\equiv1$ and $\sigma(t)=\sigma_0/(1+t)$ with $\sigma_0>0$, which is the case $\ell=1$ of the family $\sigma_0(1+t)^{-\ell}$ of Section~\ref{preliminaries}. Then $S(t)=\sigma_0\log(1+t)$, so that $e^{-S(t)}=(1+t)^{-\sigma_0}$, and
$$
A(t)=
\begin{cases}
\dfrac{(1+t)^{1-\sigma_0}-1}{1-\sigma_0}, & \sigma_0\ne1,\\[2mm]
\log(1+t), & \sigma_0=1.
\end{cases}
$$
Hence $A_\infty=\infty$ for $\sigma_0\le1$ and $A_\infty=1/(\sigma_0-1)$ for $\sigma_0>1$. Whether the effective time remains bounded depends on integrability of $e^{-S}$, a dichotomy that has no counterpart in the constant-damping case. When $A_\infty<\infty$, the entire evolution of \eqref{NA} is a finite portion of the homogeneous evolution.
In particular, concentration that would occur in the homogeneous problem only after the effective time $A_\infty$ never occurs for \eqref{NA}.

For the cosmological coefficients of the introduction, $\alpha=1/a$ and $\sigma=a'/a$, with $a\in C^1([0,\infty))$ positive, nondecreasing, and $a(0)=1$, we have $S=\log a$, so that $e^{-S}=1/a$, and
$$
A(t)=\int_0^t\frac{ds}{a(s)^2}.
$$
For the Einstein--de Sitter scale factor $a(t)=(1+t/t_0)^{2/3}$ of the introduction, where $t_0>0$ is the time elapsed since the origin of the expansion at the initial instant $t=0$, one finds $A(t)=3t_0\big(1-(1+t/t_0)^{-1/3}\big)$ and $A_\infty=3t_0$. Thus, when the gravitational field generated by the departures from homogeneity is neglected, the whole evolution corresponds to the homogeneous evolution on the bounded interval $[0,3t_0)$, and only the concentrations that the homogeneous dynamics produces before the effective time $3t_0$ ever form. For these coefficients, $A$ and $w$ coincide, up to constant reference scales and the origin of time, with the time and the velocity of the supercomoving variables of Martel and Shapiro \cite[equations~(17), (18), and (22)]{MartelShapiro1998}, in which the drag terms of the cosmological fluid equations disappear. Their analysis concerns the fluid equations with the linked coefficients $\sigma=-\alpha'/\alpha$, whereas the present results concern measure-valued entropy solutions with independent coefficients.

\subsection{The Riemann problem}
\label{riemann-example}

The previous examples illustrate the transformation on the side of the coefficients. We now illustrate it on the side of the data, with Riemann initial data
\begin{equation}\label{riemann-data}
(\rho_0,u_0)(x)=
\begin{cases}
(\rho_-,u_-), & x<0,\\
(\rho_+,u_+), & x>0,
\end{cases}
\qquad \rho_\pm>0,\ u_\pm\in\R.
\end{equation}
The measure $\rho_0=\rho_-\mathbf 1_{(-\infty,0)}\,dx+\rho_+\mathbf 1_{(0,\infty)}\,dx$ is nonnegative and locally finite, and $u_0$ is bounded, so Theorem~\ref{existence-uniqueness} gives a unique entropy solution of \eqref{NA}--\eqref{riemann-data}. This case deserves an explicit treatment because an independent construction of the solution is available.
Indeed, De~la~cruz and Neves \cite{DelaCruzNeves2025} solved this Riemann problem directly. They used the velocity rescaling $e^{S(t)}u$ and the similarity variable $x/A(t)$ \cite[(32) and (33)]{DelaCruzNeves2025}, together with a vanishing-viscosity argument for the delta shock, but neither the measure-valued solution classes nor the uniqueness theory of the homogeneous system.

The Riemann problem for \eqref{H} with data \eqref{riemann-data} is well known, see for instance Sheng and Zhang \cite{ShengZhang1999} and Huang and Wang \cite{HuangWang2001}. Three cases occur, according to the sign of $u_+-u_-$. If $u_-<u_+$, the solution contains a vacuum region, namely
\begin{equation}\label{H-rarefaction}
(\rho,w)(x,\tau)=
\begin{cases}
(\rho_-,u_-), & x<u_-\tau,\\
\big(0,\,x/\tau\big), & u_-\tau<x<u_+\tau,\\
(\rho_+,u_+), & x>u_+\tau.
\end{cases}
\end{equation}
If $u_-=u_+=:u_*$, the two states are separated by the contact discontinuity $x=u_*\tau$, and neither vacuum nor concentration is created. If $u_->u_+$, the mass concentrates on a delta shock $x=s_H(\tau):=w_\delta\tau$, and the density and the momentum are given by
\begin{align}
\rho_\tau
&=\rho_-\mathbf 1_{(-\infty,s_H(\tau))}\,dx
+W_H(\tau)\delta_{s_H(\tau)}
+\rho_+\mathbf 1_{(s_H(\tau),\infty)}\,dx,
\label{H-shock-density}\\
w(\cdot,\tau)\rho_\tau
&=\rho_-u_-\mathbf 1_{(-\infty,s_H(\tau))}\,dx
+W_H(\tau)w_\delta\,\delta_{s_H(\tau)}
+\rho_+u_+\mathbf 1_{(s_H(\tau),\infty)}\,dx,
\label{H-shock-momentum}
\end{align}
where
$$
w_\delta=\frac{\sqrt{\rho_-}\,u_-+\sqrt{\rho_+}\,u_+}{\sqrt{\rho_-}+\sqrt{\rho_+}},
\qquad
W_H(\tau)=\sqrt{\rho_-\rho_+}\,(u_--u_+)\,\tau.
$$
The shock speed satisfies $u_+<w_\delta<u_-$, which is the admissibility condition of the delta shock.
These pairs are weak solutions of \eqref{H}, see Sheng and Zhang \cite{ShengZhang1999}, and they are entropy solutions in the sense of Definition~\ref{entropy-solution}.
Indeed, the velocity in \eqref{H-rarefaction} is continuous and piecewise affine with slopes $0$ and $1/\tau$, and in the delta-shock case the velocity equal to $u_-$, $w_\delta$, and $u_+$ on $x<s_H(\tau)$, at $x=s_H(\tau)$, and on $x>s_H(\tau)$, respectively, is nonincreasing, and for $u_-=u_+$ the velocity is constant, so \eqref{homogeneous-oleinik} holds in all three cases.
The energy condition is immediate in the first two cases. In the delta-shock case, the kinetic energy differs from $u_-^2\rho_-\mathbf 1_{(-\infty,s_H(\tau))}\,dx+u_+^2\rho_+\mathbf 1_{(s_H(\tau),\infty)}\,dx$, which converges weakly to $u_0^2\rho_0$, by the measure $W_H(\tau)w_\delta^2\delta_{s_H(\tau)}$, whose total mass tends to $0$ as $\tau\to0^+$.

By Theorem~\ref{cauchy-equivalence}, the entropy solution of \eqref{NA}--\eqref{riemann-data} is obtained by substituting $\tau=A(t)$ in \eqref{H-rarefaction} or in \eqref{H-shock-density}--\eqref{H-shock-momentum} and multiplying every velocity, namely $w$, $w_\delta$, and $x/\tau$, by $e^{-S(t)}$, while the density and the weight $W_H$ of the delta shock are evaluated at $\tau=A(t)$. For $u_-<u_+$, this gives
\begin{equation}\label{NA-rarefaction}
(\rho,u)(x,t)=
\begin{cases}
(\rho_-,u_-e^{-S(t)}), & x<u_-A(t),\\
\Big(0,\,\dfrac{x\,e^{-S(t)}}{A(t)}\Big), & u_-A(t)<x<u_+A(t),\\
(\rho_+,u_+e^{-S(t)}), & x>u_+A(t).
\end{cases}
\end{equation}
Since $\rho_t$ vanishes on the vacuum region, the weak formulation does not depend on the values of $u$ there, and the value $x\,e^{-S(t)}/A(t)$ is the admissible representative inherited from the construction of Theorem~\ref{HW}, and it attains equality in \eqref{intrinsic-oleinik}. For $u_-=u_+=u_*$, the two constant states are separated by the curve $x=u_*A(t)$, and the velocity on both sides is $u_*e^{-S(t)}$. For $u_->u_+$, let $s(t):=w_\delta A(t)$, $u_\delta(t):=w_\delta e^{-S(t)}$, and $W_{\mathrm{NA}}(t):=W_H(A(t))$. The density and the momentum are given by
\begin{align*}
\rho_t
&=\rho_-\mathbf 1_{(-\infty,s(t))}\,dx
+W_{\mathrm{NA}}(t)\delta_{s(t)}
+\rho_+\mathbf 1_{(s(t),\infty)}\,dx,\\
u(\cdot,t)\rho_t
&=\rho_-u_-e^{-S(t)}\mathbf 1_{(-\infty,s(t))}\,dx
+W_{\mathrm{NA}}(t)u_\delta(t)\,\delta_{s(t)}
+\rho_+u_+e^{-S(t)}\mathbf 1_{(s(t),\infty)}\,dx.
\end{align*}
These measures define an entropy solution of \eqref{NA}, since they are the image of \eqref{H-shock-density}--\eqref{H-shock-momentum} under the bijection of Theorem~\ref{entropy-correspondence}. In particular, the jump relations across the delta shock, which in the direct approach must be imposed separately, need no further verification. As a direct check, the delta shock moves with speed $s'(t)=w_\delta A'(t)=\alpha(t)u_\delta(t)$, in accordance with Theorem~\ref{characteristics-transfer}.

With $A$ given by \eqref{effective-time}, formula \eqref{NA-rarefaction} coincides with the solution given for $u_-<u_+$ in De~la~cruz and Neves \cite[Section~3.1]{DelaCruzNeves2025}, since the regions, the boundaries $x=u_\pm A(t)$, and the velocity $x\,e^{-S(t)}/A(t)$ in the vacuum region are the same.
For $u_->u_+$, the location, the velocity, and the weight of the delta shock,
$$
s(t)=w_\delta A(t),\qquad
u_\delta(t)=w_\delta e^{-S(t)},\qquad
W_{\mathrm{NA}}(t)=\sqrt{\rho_-\rho_+}\,(u_--u_+)A(t),
$$
also coincide with those of De~la~cruz and Neves \cite[Theorem~3.1]{DelaCruzNeves2025}, both for $\rho_-\ne\rho_+$ and for $\rho_-=\rho_+$. In particular, the Riemann solutions of De~la~cruz and Neves \cite{DelaCruzNeves2025}, in which the delta shock was selected by the overcompressive condition $\alpha(t)u_+e^{-S(t)}<s'(t)<\alpha(t)u_-e^{-S(t)}$ \cite[(44)]{DelaCruzNeves2025}, are entropy solutions in the sense of Definition~\ref{entropy-solution}, and Theorem~\ref{uniqueness} shows that no other entropy solution exists in the much larger class of measure-valued solutions. The earlier construction uses similarity reduction and viscosity, while the present argument transfers the homogeneous entropy theory. Both use the same coefficient-dependent velocity rescaling and effective-time expression.

\subsection{Immediate concentration}
\label{immediate-concentration}

The last example concerns the exceptional set $\Omega$ defined at the end of Section~\ref{strips}, formed by the initial points whose particles belong to a concentration at every positive time. We show that it may have full Lebesgue measure on an interval, although the curves of discontinuity are countable. By Theorem~\ref{cauchy-equivalence}, it suffices to describe the homogeneous solution. We use the system \eqref{H}, the generalized potential $G$ of Section~\ref{strips} with $\bar a=0$, which is allowed because $\rho_0$ has no atoms, the point $y_*^H$ of Proposition~\ref{potential-transfer}, and the forward generalized characteristics $x_H$ of Section~\ref{generalized-characteristics}.

For $y\in\R$ let $\mathrm{frac}(y):=y-\lfloor y\rfloor\in[0,1)$ be its fractional part. Let $\rho_0:=\mathbf 1_{[0,1)}\,dx$ and
$$
u_0(\eta):=-\sum_{k\ge1}k^{-2}\,\mathrm{frac}(2^k\eta),\qquad\eta\in[0,1),
$$
and put $s_n:=\sum_{k=1}^nk^{-2}2^k$ for $n\ge1$ and $s_0:=0$, with the convention $1/s_0:=+\infty$. The series converges uniformly, so $u_0$ is Borel and $-\pi^2/6\le u_0\le0$. Since $s_n\to\infty$, every $\tau>0$ belongs to exactly one interval $[1/s_n,1/s_{n-1})$ with $n\ge1$.
The function $u_0$ has infinite total variation on $(0,1)$. Indeed, one-sided limits can be taken term by term. If $\eta_0\in(0,1)$ is an odd multiple of $2^{-j}$, then $\mathrm{frac}(2^k\cdot)$ is continuous at $\eta_0$ for $k<j$ and jumps from $1$ to $0$ at $\eta_0$ for $k\ge j$, so
$$
u_0(\eta_0+0)-u_0(\eta_0-0)=\sum_{k\ge j}k^{-2}\ge\sum_{k\ge j}\frac{1}{k(k+1)}=\frac1j.
$$
There are $2^{j-1}$ such points, so the total variation of $u_0$ is at least $\sum_{j\ge1}2^{j-1}/j=\infty$.

For $n\ge1$, the dyadic intervals of level $n$ are $D=[\lambda,\lambda+2^{-n})$ with $\lambda\in2^{-n}\{0,1,\dots,2^n-1\}$. For $\tau>0$ put $\Theta^0_\tau(\eta):=\eta+\tau u_0(\eta)$, the free position at time $\tau$ of the particle issued from $\eta$, and let
$$
c_D(\tau):=2^n\int_D\Theta^0_\tau(\eta)\,d\eta
$$
be the barycenter of the free positions of the particles issued from $D$.

\begin{lemma}\label{dyadic-lemma}
Let $n\ge1$ and $\tau>0$.
\begin{enumerate}[label=(\roman*)]
\item If $\tau\ge1/s_n$, then, for every dyadic interval $D=[\lambda,\lambda+2^{-n})$ of level $n$, $\int_\lambda^y(\Theta^0_\tau-c_D(\tau))\,d\eta\ge0$ for $y\in[\lambda,\lambda+2^{-n}]$, with equality for $y=\lambda$ and $y=\lambda+2^{-n}$.
\item If $D=[\lambda,\lambda+2^{-n})$ and $D'=D+2^{-n}$ are consecutive dyadic intervals of level $n$ in $[0,1)$, then $c_{D'}(\tau)-c_D(\tau)\ge2^{-n}(1-\tau s_{n-1})$, with equality if and only if $D$ and $D'$ are the two halves of a dyadic interval of level $n-1$. In particular, $c_D(\tau)$ is strictly increasing in $D$ if $\tau<1/s_{n-1}$.
\end{enumerate}
\end{lemma}

\begin{proof}
Let $\eta\in D=[\lambda,\lambda+2^{-n})$. For $k\le n$, the interval $D$ is contained in a dyadic interval of length $2^{-k}$, so $\mathrm{frac}(2^k\eta)=\mathrm{frac}(2^k\lambda)+2^k(\eta-\lambda)$. Hence
$$
\Theta^0_\tau(\eta)=\lambda+\tau a_D+(1-\tau s_n)(\eta-\lambda)+\tau f_n(\eta),
$$
where
$$
a_D:=-\sum_{k\le n}k^{-2}\,\mathrm{frac}(2^k\lambda),
\qquad
f_n(\eta):=-\sum_{k>n}k^{-2}\,\mathrm{frac}(2^k\eta).
$$
For $k>n$, the interval $D$ is the union of $2^{k-n}$ periods of $\mathrm{frac}(2^k\cdot)$, each starting at a multiple of $2^{-k}$, so the mean of this function over $D$ is $\tfrac12$. Over $[\lambda,y]$ with $y\in D$, its mean is at most $\tfrac12$, since $[\lambda,y]$ consists of complete periods and of at most one incomplete period $[\beta,\beta+r2^{-k})$, $0<r<1$, over which the mean is $r/2$.

(i) It suffices to show that, for $y\in(\lambda,\lambda+2^{-n}]$, the mean of $\Theta^0_\tau$ over $[\lambda,y]$ is at least $c_D(\tau)$, its mean over $D$. The mean of $(1-\tau s_n)(\eta-\lambda)$ over $[\lambda,y]$ is $(1-\tau s_n)(y-\lambda)/2$, which is at least its mean $(1-\tau s_n)2^{-n-1}$ over $D$, because $1-\tau s_n\le0$. By the preceding paragraph, the mean of $f_n$ over $[\lambda,y]$ is at least $-\tfrac12\sum_{k>n}k^{-2}$, which is its mean over $D$. Equality at $y=\lambda+2^{-n}$ is the definition of $c_D(\tau)$.

(ii) Put $\lambda':=\lambda+2^{-n}$, and let $j\in\{1,\dots,n\}$ be such that $\lambda'$ is an odd multiple of $2^{-j}$. For $k<j$, $\mathrm{frac}(2^k\lambda')=\mathrm{frac}(2^k\lambda)+2^{k-n}$, and for $j\le k\le n$, $\mathrm{frac}(2^k\lambda')=\mathrm{frac}(2^k\lambda)+2^{k-n}-1$, since $2^k\lambda'$ is an integer exactly for $k\ge j$. Hence $a_{D'}-a_D=-2^{-n}s_n+\sum_{k=j}^nk^{-2}$. Since $f_n$ has the same mean on $D$ and on $D'$,
$$
c_{D'}(\tau)-c_D(\tau)=2^{-n}+\tau\big(a_{D'}-a_D\big)=2^{-n}(1-\tau s_n)+\tau\sum_{k=j}^nk^{-2}\ge2^{-n}(1-\tau s_n)+\tau n^{-2},
$$
and the right-hand side equals $2^{-n}(1-\tau s_{n-1})$, because $s_n-s_{n-1}=n^{-2}2^n$. Equality holds exactly when $j=n$, that is, when $D$ and $D'$ are the two halves of a dyadic interval of level $n-1$.
\end{proof}

Thus, for $n\ge2$, the barycenters of sibling intervals meet at $\tau=1/s_{n-1}$, while those of non-sibling intervals remain apart. For $n=1$, the convention $1/s_0=\infty$ means that the two atoms of level $1$ never meet.

\begin{proposition}\label{dyadic-solution}
Let $n\ge1$, $\tau\in[1/s_n,1/s_{n-1})$, and $x\in\R$, and let $\beta(x,\tau)$ be the total length of the dyadic intervals $D$ of level $n$ with $c_D(\tau)<x$. If $x$ is different from all the points $c_D(\tau)$, then $y_*^H(x,\tau)=\beta(x,\tau)$. Consequently, the density of the solution of Theorem~\ref{HW} with data $(\rho_0,u_0)$ is
$$
\rho_\tau=\sum_D2^{-n}\,\delta_{c_D(\tau)},
$$
where the sum runs over the $2^n$ dyadic intervals of level $n$ in $[0,1)$, and the points $c_D(\tau)$ are pairwise distinct.
\end{proposition}

\begin{proof}
Since $\bar a=0$, $G(y;x,\tau)=\int_0^y(\Theta^0_\tau-x)\,d\eta$ for $y\in[0,1]$, $G(y;x,\tau)=0$ for $y\le0$, and $G(y;x,\tau)=G(1;x,\tau)$ for $y\ge1$. Let $\bar\Theta^0_\tau$ be equal to $c_D(\tau)$ on each dyadic interval $D$ of level $n$. For $y\in[0,1]$,
$$
G(y;x,\tau)=F_1(y)+F_2(y),\qquad F_1(y):=\int_0^y\big(\Theta^0_\tau-\bar\Theta^0_\tau\big)\,d\eta,\qquad F_2(y):=\int_0^y\big(\bar\Theta^0_\tau-x\big)\,d\eta.
$$
By Lemma~\ref{dyadic-lemma}(i), $F_1\ge0$ on $[0,1]$, with equality at the endpoints of the dyadic intervals of level $n$. The function $F_2$ is piecewise linear, with slope $c_D(\tau)-x$ on $D$, and by Lemma~\ref{dyadic-lemma}(ii) these slopes are strictly increasing in $D$. In particular, the points $c_D(\tau)$ are pairwise distinct. If $x$ is different from all of them, the slope of $F_2$ is negative on the intervals $D$ with $c_D(\tau)<x$ and positive on the others. These intervals form an initial segment of $[0,1)$ of length $\beta:=\beta(x,\tau)$, so $F_2$ attains its minimum over $[0,1]$ only at $\beta$. Since $F_1(\beta)=0\le F_1$, the same holds for $G(\cdot;x,\tau)$. If $0<\beta<1$, then $F_2(\beta)<0=G(y;x,\tau)$ for $y\le0$ and $F_2(\beta)<F_2(1)=G(y;x,\tau)$ for $y\ge1$, so $\beta$ is the only point at which $G(\cdot;x,\tau)$ attains its infimum over $\R$, and $y_*^H(x,\tau)=\beta$. If $\beta\in\{0,1\}$, the infimum is attained on $(-\infty,0]$ or on $[1,\infty)$, and $y_*^H(x,\tau)$, which belongs to the support $[0,1]$ of $\rho_0$, equals $\beta$.
Since $\rho_0$ has no atoms, the two cases in \cite[(2.35)]{HuangWang2001} coincide, and the mass primitive of the solution of Theorem~\ref{HW} is $m_0(y_*^H(x,\tau))=\beta(x,\tau)$ at every $x$ different from the points $c_D(\tau)$. This function is constant between consecutive points $c_D(\tau)$ and jumps by $2^{-n}$ at each of them, which gives the formula for $\rho_\tau$.
\end{proof}

For $\tau\in[1/s_n,1/s_{n-1})$ the solution therefore consists of $2^n$ atoms of mass $2^{-n}$, and, for $n\ge2$, sibling atoms merge at $\tau=1/s_{n-1}$. The curves of discontinuity are the trajectories of the atoms, one for each dyadic interval, so they are countable. Nevertheless, the particle issued from every point of $(0,1)$ that is not a dyadic rational lies in a concentration at every positive time.

\begin{proposition}\label{dyadic-exceptional}
For $\tau>0$, let $\gamma_\tau:\R\to\R$ be nondecreasing maps such that, for every $n\ge1$, every $\tau\in[1/s_n,1/s_{n-1})$, and every dyadic interval $D$ of level $n$, $\gamma_\tau(\eta)=c_D(\tau)$ for $\eta$ in the interior of $D$. Then every $\eta_0\in(0,1)$ that is not a dyadic rational belongs to the set $\Omega$ of Section~\ref{strips} defined with these maps. The hypothesis holds for every nondecreasing extension $\gamma_\tau:\R\to\R$ of the forward generalized characteristics $x_H(\cdot,\tau)$ of Section~\ref{generalized-characteristics}, which are defined on $\supp\rho_0=[0,1]$, and for every family of nondecreasing maps $\gamma_\tau:\R\to\R$ with $(\gamma_\tau)_\#\rho_0=\rho_\tau$.
\end{proposition}

\begin{proof}
Let $\tau>0$, let $n\ge1$ be such that $\tau\in[1/s_n,1/s_{n-1})$, let $D=[\lambda,\lambda+2^{-n})$ be the dyadic interval of level $n$ containing $\eta_0$, and put $x:=c_D(\tau)$. Then $\eta_0$ lies in the interior of $D$, and $\gamma_\tau(\eta)=x$ for every $\eta\in(\lambda,\lambda+2^{-n})$. By monotonicity, every $\eta\in\R$ with $\gamma_\tau(\eta)<x$ satisfies $\eta\le\lambda$, and every $\eta\in\R$ with $\gamma_\tau(\eta)>x$ satisfies $\eta\ge\lambda+2^{-n}$. Hence $\xi^-(x,\tau)\le\lambda<\eta_0<\lambda+2^{-n}\le\xi^+(x,\tau)$. Thus $\eta_0\in\Omega$.
For the forward characteristics, \cite[(2.46)]{HuangWang2001} gives $x_H(\eta,\tau)=x$ whenever $y_*^H(x-0,\tau)<\eta<y_*^H(x+0,\tau)$ and $\eta\in\supp\rho_0$. By Proposition~\ref{dyadic-solution}, $y_*^H(c_D(\tau)-0,\tau)=\lambda$ and $y_*^H(c_D(\tau)+0,\tau)=\lambda+2^{-n}$, and $x_H(\cdot,\tau)$ is nondecreasing on $[0,1]$, which gives the hypothesis for every nondecreasing extension. Finally, let $\gamma_\tau$ be a nondecreasing map that transports $\rho_0$ onto the measure of Proposition~\ref{dyadic-solution}. The transport determines $\gamma_\tau$ only $\rho_0$-almost everywhere, so we use the preimages of half-lines. Since the points $c_E(\tau)$ are strictly increasing in $E$, the transport property and the monotonicity of $\gamma_\tau$ show that the sets $\{\eta\in[0,1):\gamma_\tau(\eta)<c_D(\tau)\}$ and $\{\eta\in[0,1):\gamma_\tau(\eta)\le c_D(\tau)\}$ are initial segments of $[0,1)$ of lengths $\lambda$ and $\lambda+2^{-n}$, so $\gamma_\tau=c_D(\tau)$ on the interior of $D$.
\end{proof}

The maps of Proposition~\ref{dyadic-exceptional} are not the curves with which Huang and Wang define $\Omega_s$. We now show that the same conclusion holds for those curves. Let $w$ be the velocity of Theorem~\ref{HW} for the data of this subsection, and let $M:=\|u_0\|_{L^\infty(d\rho_0)}=\pi^2/6$, so that $|w|\le M$. Following \cite[(3.4)]{HuangWang2001} with initial time $0$, let $w_\varepsilon:=w*j_\varepsilon$, where $j_\varepsilon$ is a standard mollifier, and, for $\eta\in\R$, let $\gamma^\varepsilon(\eta,\cdot)$ be the solution of
$$
\frac{dx}{d\tau}=w_\varepsilon(x,\tau),\qquad x(0)=\eta.
$$
Whether the convolution is taken in $x$ or in $(x,\tau)$, in the latter case with $w$ extended by $0$ to $\tau<0$, the function $w_\varepsilon$ is bounded by $M$ and Lipschitz in $x$, so this solution exists and is unique. By \cite[Lemma~3.1(1)]{HuangWang2001}, there is a sequence $\varepsilon_i\to0^+$ such that the limit $\gamma(\eta,\tau):=\lim\limits_{i\to\infty}\gamma^{\varepsilon_i}(\eta,\tau)$ exists for all $\eta\in\R$ and $\tau\ge0$. The set $\Omega_s$ is the set $\Omega$ of Section~\ref{strips} defined with $\gamma_\tau:=\gamma(\cdot,\tau)$.

\begin{proposition}\label{HW-exceptional}
For every choice of the sequence $\varepsilon_i$, the set $\Omega_s$ contains every point of $(0,1)$ that is not a dyadic rational, except for a countable set. In particular, $\Omega_s$ has Lebesgue measure at least $1$.
\end{proposition}

\begin{proof}
We first collect the properties of $\gamma$ that we use. Uniqueness for the approximate equations implies that the maps $\gamma^\varepsilon(\cdot,\tau)$ are nondecreasing, and $|w_\varepsilon|\le M$ gives $|\gamma^\varepsilon(\eta,\tau)-\gamma^\varepsilon(\eta,\tau')|\le M|\tau-\tau'|$. Hence $\gamma(\cdot,\tau)$ is nondecreasing, $\gamma(\eta,\cdot)$ is $M$-Lipschitz, $|\gamma(\eta,\tau)-\eta|\le M\tau$, and, for each $\eta$, $\gamma^{\varepsilon_i}(\eta,\cdot)\to \gamma(\eta,\cdot)$ uniformly on compact subsets of $[0,\infty)$. Let $\mathcal R$ be the set of points of $(0,1)$ that are not dyadic rationals.

Next we describe $w$ in the vacuum. Let $n\ge1$ and $\tau\in[1/s_n,1/s_{n-1})$. For a dyadic interval $D=[\lambda,\lambda+2^{-n})$ of level $n$, the definition of $c_D(\tau)$ gives
$$
c_D(\tau)=\lambda+2^{-n-1}+\tau\bar u_D,\qquad\bar u_D:=2^n\int_Du_0\,d\eta\in[-M,0].
$$
We call a point $(x,\tau)$ with $\tau>0$ a vacuum point if $x$ is different from all the points $c_D(\tau)$. Let $x\in\R$. With the decomposition $G=F_1+F_2$ on $[0,1]$ used in the proof of Proposition~\ref{dyadic-solution}, $F_1\ge0=F_1(\beta(x,\tau))$, and $F_2$ is strictly decreasing on $[0,\beta(x,\tau)]$, because its slopes $c_D(\tau)-x$ are negative there. Hence $G(y;x,\tau)>G(\beta(x,\tau);x,\tau)$ for $y\in[0,\beta(x,\tau))$. Since $G(\cdot;x,\tau)$ is continuous, every limit of a minimizing sequence is a minimizer, and since $y_*^H(x,\tau)$ is such a limit lying in $\supp\rho_0=[0,1]$, we obtain
\begin{equation}\label{ystar-lower}
y_*^H(x,\tau)\ge\beta(x,\tau)\qquad\text{for every }x\in\R,
\end{equation}
with equality at vacuum points by Proposition~\ref{dyadic-solution}.
Let $(x_1,\tau_1)$ be a vacuum point with $\tau_1\in[1/s_n,1/s_{n-1})$, and put $\beta:=\beta(x_1,\tau_1)$. We claim that
\begin{equation}\label{vacuum-velocity}
w(x_1,\tau_1)=
\begin{cases}
(x_1-\beta)/\tau_1, & |x_1-\beta|\le M\tau_1,\\
-M, & x_1-\beta<-M\tau_1,\\
M, & x_1-\beta>M\tau_1.
\end{cases}
\end{equation}
The last two cases are the definition of $w$, since $y_*^H(x_1,\tau_1)=\beta$. In the first case, \cite[(2.29)]{HuangWang2001} defines $w(x_1,\tau_1)$ as the derivative at $\tau_1$ of the forward generalized characteristic issued from $(x_1,\tau_1)$, whose position at a time $\tau_0>\tau_1$ is, by the construction at the beginning of the proof of \cite[Lemma~2.6]{HuangWang2001},
$$
\chi(\tau_0)=\sup\{x\in\R:p(x)<x_1\},\qquad
p(x):=x-\frac{x-y_*^H(x,\tau_0)}{\tau_0}(\tau_0-\tau_1).
$$
Put $r(\tau_0):=\beta+(x_1-\beta)\tau_0/\tau_1$. For $\tau_0>\tau_1$ close to $\tau_1$, the point $r(\tau_0)$ lies in the open interval between the two atoms, or the half-line beyond the extreme atom, that contains $x_1$ at time $\tau_1$, so $\beta(x,\tau_0)=\beta$ for $x$ in that interval near $r(\tau_0)$. For such $x$, $p(x)=\beta+(x-\beta)\tau_1/\tau_0$, which is less than $x_1$ for $x<r(\tau_0)$. For $x\ge r(\tau_0)$, \eqref{ystar-lower} and the monotonicity of $\beta(\cdot,\tau_0)$ give $y_*^H(x,\tau_0)\ge\beta$, hence
$$
p(x)=y_*^H(x,\tau_0)\Big(1-\frac{\tau_1}{\tau_0}\Big)+x\,\frac{\tau_1}{\tau_0}\ge\beta\Big(1-\frac{\tau_1}{\tau_0}\Big)+r(\tau_0)\,\frac{\tau_1}{\tau_0}=x_1.
$$
Therefore $\chi(\tau_0)=r(\tau_0)$, and the derivative at $\tau_1$ is $(x_1-\beta)/\tau_1$, which proves \eqref{vacuum-velocity}. In particular, on each connected component of the set of vacuum points, $\beta$ is constant and $w$ is continuous, and $w(\cdot,\tau)$ is Lipschitz with constant $1/\tau$ on each interval of vacuum.

We show that, for $\eta\in \mathcal R$ and every $\tau>0$, the point $\gamma(\eta,\tau)$ is an atom of $\rho_\tau$. Suppose that $(\gamma(\eta,\tau_b),\tau_b)$ is a vacuum point for some $\tau_b>0$. On every interval $[\delta,\tau_b]$ with $\delta>0$ the atoms move along finitely many continuous curves, so the set of $\tau\in(0,\tau_b]$ at which $\gamma(\eta,\tau)$ is an atom of $\rho_\tau$ is relatively closed in $(0,\tau_b]$. Let $\tau_a$ be its supremum, with $\tau_a:=0$ if the set is empty. Then $\tau_a<\tau_b$, and $(\gamma(\eta,\tau),\tau)$ is a vacuum point for every $\tau\in(\tau_a,\tau_b]$.
Within each interval $[1/s_n,1/s_{n-1})$ the atoms are strictly ordered by Lemma~\ref{dyadic-lemma}(ii), so, by continuity, $\gamma(\eta,\cdot)$ stays between the same two consecutive atoms, or to the left of the first one, or to the right of the last one, and $\beta(\gamma(\eta,\tau),\tau)$ is constant. At a time $1/s_{n-1}\in(\tau_a,\tau_b]$, two sibling atoms of level $n$ meet, since equality holds in Lemma~\ref{dyadic-lemma}(ii) for siblings, and they continue as the atom of their parent, whose barycenter is the mean of theirs. The curve $\gamma(\eta,\cdot)$ cannot lie between two siblings just before that time, since it would then meet their common atom. Two consecutive intervals of level $n$ that are not siblings have the same common endpoint as their parents. Hence there is a dyadic rational $\beta\in[0,1]$ with $\beta(\gamma(\eta,\tau),\tau)=\beta$ for every $\tau\in(\tau_a,\tau_b]$.

Assume first that $\tau_a=0$, and let $\tau$ be small, with $\tau\in[1/s_n,1/s_{n-1})$. If $0<\beta<1$, then $\gamma(\eta,\tau)$ lies between the atoms of $D=[\beta-2^{-n},\beta)$ and $D'=[\beta,\beta+2^{-n})$. Since $c_D(\tau)\ge\beta-2^{-n-1}-M\tau$ and $c_{D'}(\tau)\le\beta+2^{-n-1}$, we get $|\gamma(\eta,\tau)-\beta|\le2^{-n-1}+M\tau$, hence $|\eta-\beta|\le2^{-n-1}+2M\tau$. Letting $\tau\to0^+$, so that $n\to\infty$, gives $\eta=\beta$, which is impossible because $\eta$ is not a dyadic rational. If $\beta=0$, then $\gamma(\eta,\tau)<c_{D'}(\tau)\le2^{-n-1}$ with $D'=[0,2^{-n})$, so $\eta<2^{-n-1}+M\tau$, and letting $\tau\to0^+$ contradicts $\eta>0$. The case $\beta=1$ is symmetric.

Assume now that $\tau_a>0$, and let $n$ be such that $\tau_a\in[1/s_n,1/s_{n-1})$. By continuity, $\gamma(\eta,\tau_a)$ is one of the atoms adjacent to the vacuum interval that contains $\gamma(\eta,\tau)$ for $\tau$ slightly larger than $\tau_a$. We first show that $\gamma(\eta,\cdot)$ solves the characteristic equation away from $\tau_a$. Let $[\tau',\tau'']\subset(\tau_a,\tau_b]$. The compact set $\{(\gamma(\eta,\tau),\tau):\tau'\le\tau\le\tau''\}$ lies in the open set of vacuum points, and so does a compact neighborhood $K$ of it, on which $w$ is uniformly continuous. Hence $w_\varepsilon\to w$ uniformly on a smaller neighborhood of the curve, and the uniform convergence of $\gamma^{\varepsilon_i}(\eta,\cdot)$ on $[\tau',\tau'']$ allows us to pass to the limit in $\gamma^{\varepsilon_i}(\eta,\tau'')-\gamma^{\varepsilon_i}(\eta,\tau')=\int_{\tau'}^{\tau''}w_{\varepsilon_i}(\gamma^{\varepsilon_i}(\eta,\tau),\tau)\,d\tau$. Thus $\partial_\tau \gamma(\eta,\tau)=w(\gamma(\eta,\tau),\tau)$ for every $\tau\in(\tau_a,\tau_b)$.
Suppose that $\gamma(\eta,\tau_a)=c_D(\tau_a)$ with $D=[\beta-2^{-n},\beta)$, and put $d(\tau):=\gamma(\eta,\tau)-c_D(\tau)$. Then $d$ is Lipschitz, $d(\tau_a)=0$, and there is $\delta>0$ with $0<d<2^{-n-1}$ on $(\tau_a,\tau_a+\delta]$, and $d'(\tau)=w(\gamma(\eta,\tau),\tau)-\bar u_D$ there. Moreover, $\gamma(\eta,\tau)-\beta=d(\tau)-2^{-n-1}+\tau\bar u_D<0$, so, by \eqref{vacuum-velocity}, either $w(\gamma(\eta,\tau),\tau)=-M\le\bar u_D$ or
$$
w(\gamma(\eta,\tau),\tau)=\frac{\gamma(\eta,\tau)-\beta}{\tau}=\bar u_D+\frac{d(\tau)-2^{-n-1}}{\tau}<\bar u_D.
$$
Thus $d'\le0$ on $(\tau_a,\tau_a+\delta]$, and $d\le d(\tau_a)=0$ there by continuity, a contradiction. If $\gamma(\eta,\tau_a)=c_{D'}(\tau_a)$ with $D'=[\beta,\beta+2^{-n})$, the same argument applied to $e(\tau):=c_{D'}(\tau)-\gamma(\eta,\tau)$, for which $0<e<2^{-n-1}$ gives $\gamma(\eta,\tau)-\beta>-M\tau$, shows that $e'\le0$, again a contradiction. Hence $\gamma(\eta,\tau)$ is an atom of $\rho_\tau$ for every $\eta\in \mathcal R$ and $\tau>0$.

Consequently, for $\eta\in \mathcal R$, the curve $\gamma(\eta,\cdot)$ follows the atoms. Within an interval $[1/s_n,1/s_{n-1})$ the atoms move along finitely many continuous and pairwise disjoint curves, so $\gamma(\eta,\cdot)$ follows one of them, say that of $D$, and at time $1/s_{n-1}$ it continues along the atom of the parent of $D$. In particular, if $\eta,\eta'\in \mathcal R$ and $\gamma(\eta,\tau_1)=\gamma(\eta',\tau_1)$, then $\gamma(\eta,\tau)=\gamma(\eta',\tau)$ for every $\tau\ge\tau_1$.

For $\eta\in\R$ and $\tau>0$ let $J_\tau(\eta):=\{\eta'\in\R:\gamma(\eta',\tau)=\gamma(\eta,\tau)\}$, which is an interval containing $\eta$ because $\gamma(\cdot,\tau)$ is nondecreasing. With $x:=\gamma(\eta,\tau)$, the definitions of Section~\ref{strips} give $\xi^-(x,\tau)=\inf J_\tau(\eta)$ and $\xi^+(x,\tau)=\sup J_\tau(\eta)$. Hence $\eta\in\Omega_s$ whenever $J_\tau(\eta)$ is not reduced to a point for every $\tau>0$. Let $B_\tau$ be the set of $\eta\in \mathcal R$ with $J_\tau(\eta)=\{\eta\}$. Distinct points of $B_\tau$ have distinct images under $\gamma(\cdot,\tau)$, and these images are atoms of $\rho_\tau$, so $B_\tau$ is finite. If $\eta\in \mathcal R\setminus B_{\tau'}$ and $\tau>\tau'$, then $J_{\tau'}(\eta)$ is a nondegenerate interval, so it contains a point $\eta'\in \mathcal R$ different from $\eta$, and $\gamma(\eta',\tau)=\gamma(\eta,\tau)$ by the preceding paragraph. Thus $\eta\notin B_\tau$, that is, $B_\tau\subset B_{\tau'}$ for $\tau'<\tau$. Every $\eta\in \mathcal R\setminus\bigcup_{k\ge1}B_{1/k}$ therefore satisfies $J_\tau(\eta)\ne\{\eta\}$ for every $\tau>0$, and belongs to $\Omega_s$. Since $\bigcup_kB_{1/k}$ is countable, the proposition follows.
\end{proof}

Hence, in this example, the set $\Omega_s$, whose nullity is used in the passage to the initial time described in Section~\ref{strips}, has Lebesgue measure at least $1$, although the curves of discontinuity are countable.\\

{\bf Acknowledgments:}
The author R. De la cruz gratefully acknowledges Universidad Pedag\'ogica y Tecnol\'ogica de Colombia (UPTC) for providing the academic environment and dedicated time that made this research possible. The author R. Mondal gratefully acknowledges IIT Bombay for providing the excellent research environment.

{\bf Author Contributions:}
All authors contributed to the study conception and design.\\

{\bf Funding:}
The author W. Neves has received research grants from CNPq, Brazil through the grants 313005/2023-0, 403675/2025-1, and also by FAPERJ (Cientista do Nosso Estado), Brazil through the grant E-26/204.171/2024.
The author R. Mondal has received support from the Institute Postdoctoral Fellowship (F\&SS/Rect/33/2026) at the Indian Institute of Technology Bombay (IIT Bombay).

\section*{Declarations} 
{\bf Conflict of interest:} The authors declare that they have no known competing financial interests or personal relationships that could have appeared to influence the work reported in this paper.

\end{document}